\documentclass[brochure,12pt]{bourbaki}
\usepackage[matrix,arrow]{xy}
\usepackage{amssymb,amsfonts,amsmath,footnote,url}
\usepackage[french=guillemets]{csquotes}
\usepackage{tikz}
\usetikzlibrary{shapes,arrows}
\usepackage[Algorithme]{algorithm}
\usepackage[noend]{algpseudocode}
\usepackage[russian,francais]{babel}
\usepackage[utf8]{inputenc} 
\usepackage[T1, T2A]{fontenc}    
\newcommand{\El}{\mbox{\usefont{T2A}{\rmdefault}{m}{n}\CYRL}}
\DeclareMathOperator{\Sym}{Sym}
\DeclareMathOperator{\Iso}{Iso}
\DeclareMathOperator{\Aut}{Aut}
\DeclareMathOperator{\Alt}{Alt}
\renewcommand{\labelenumi}{(\alph{enumi})}
\renewcommand{\theenumi}{\alph{enumi}}

\newtheorem{exer}[defi]{Exercice}

\DeclareFontFamily{OT1}{rsfs}{}
\DeclareFontShape{OT1}{rsfs}{n}{it}{<-> rsfs10}{}
\DeclareMathAlphabet{\mathscr}{OT1}{rsfs}{n}{it}

\tikzstyle{decision} = [diamond, draw, fill=gray!20,
    text width=4.5em, text badly centered, node distance=3cm, inner sep=0pt]
\tikzstyle{block} = [rectangle, draw, fill=gray!20,
  text width=5em, text centered, rounded corners, minimum height=4em]
\tikzstyle{noblock} = [rectangle, draw, 
  text width=5em, text centered, rounded corners, minimum height=4em]
\tikzstyle{rect} = [rectangle, draw, text centered, minimum height=4em]
\tikzstyle{line} = [draw, very thick, color=black!50, -latex']
\tikzstyle{cloudin} = [draw, ellipse, text width = 5.75em, node distance=3cm,
  minimum height=5em]
\tikzstyle{cloudout} = [draw, ellipse,  text width = 4em, node distance=3cm,
  minimum height=2em]

\date{Janvier 2017}
\bbkannee{69\`eme ann\'ee, 2016-2017}
\bbknumero{1125}
\title{Isomorphismes de graphes en temps quasi-polynomial}
\subtitle{d'apr\`es Babai et Luks, Weisfeiler-Leman, \dots}
\author{Harald Andrés HELFGOTT}
\address{Universit\"at G\"ottingen\\
Mathematisches Institut\\
Bunsenstrasse 3-5\\
D-37073 G\"ottingen\\
Allemagne}
\email{helfgott@math.univ-paris-diderot.fr}

\begin{document}
\maketitle

\noindent{\em Résumé:}
Soient donnés deux graphes $\Gamma_1$, $\Gamma_2$ à $n$ sommets. Sont-ils
isomorphes?
S'ils le sont, l'ensemble des isomorphismes de $\Gamma_1$ à
$\Gamma_2$ peut être identifié avec une classe $H\cdot\pi$ du groupe
symétrique sur $n$ éléments. Comment trouver $\pi$ et des générateurs de $H$?

Le défi de donner un algorithme toujours efficace en réponse à ces questions est resté longtemps ouvert. Babai a récemment montré comment résoudre ces questions -- et d'autres qui y sont liées -- en temps quasi-polynomial, c'est-à-dire
en temps $\exp\left(O(\log n)^{O(1)}\right)$.
Sa stratégie est basée en partie sur l'algorithme de Luks (1980/82), qui a résolu le cas de graphes de degré borné.


\section{Introduction}

Soient $\mathbf{x}$, $\mathbf{y}$ deux chaînes de caractères, à savoir,
deux applications $\Omega\to \Sigma$, où $\Sigma$ (l'{\em alphabet}) et $\Omega$
(le {\em domaine}) sont des ensembles finis. Tout groupe de permutations\footnote{Pour nous, $G<S$ (ou $S>G$) veut dire {\og $G$ est un sous-groupe de $S$,
  pas forcement propre.\fg}}
$G<\Sym(\Omega)$ agit sur l'ensemble $\Sigma^\Omega$
des chaînes de domaine $\Omega$ sur un alphabet $\Sigma$.
Pour nous, {\em décrire un groupe $G$}, ou {\em être donné un groupe $G$},
voudra toujours dire \enquote{donner, voire être donné, un ensemble de
  générateurs de $G$}; {\em décrire une classe}
$H \pi$ voudra dire \enquote{donner un élément $\pi$ de la classe
et un ensemble de générateurs de $H$}. 

Le {\em problème de l'isomorphisme de chaînes} consiste à déterminer,
étant donnés $\mathbf{x}$, $\mathbf{y}$ et $G$,
s'il y a au moins un élément $\pi$ de $G$ qui envoie
$\mathbf{x}$ sur $\mathbf{y}$, et, si de tels éléments ({\em isomorphismes})
existent, à les décrire. Il est clair que l'ensemble des isomorphismes
$\Iso_G(\mathbf{x},\mathbf{y})$ forme
une classe $\Aut_G(\mathbf{x}) \pi$ du groupe
$\Aut_G(\mathbf{x})$ d'automorphismes de $\mathbf{x}$ dans~$G$,
c'est-à-dire du groupe consistant dans les éléments de $G$ qui envoient $\mathbf{x}$ sur
lui-même. 

Le défi consiste à donner un algorithme qui résolve le problème en
temps polynomial en la taille $n = |\Omega|$ de $\Omega$,
voire en temps raisonnable.
Par exemple, le temps employé pourrait être {\em quasi-polynomial} en $n$,
ce qui veut dire $\exp\left(O(\log n)^{O(1)}\right)$.
Ici,
comme toujours,
$O(f(n))$ désigne une quantité bornée par $C\cdot f(n)$, pour $n$
assez grand et $C>0$ une constante, et $O_\epsilon$ indique que la constante
$C$ dépend de $\epsilon$.

Une grande partie de la motivation pour le problème de l'isomorphisme de
chaînes vient du fait que le {\em problème de l'isomorphisme de graphes}
se réduit à lui. Ce problème consiste à déterminer si deux graphes finis
$\Gamma_1$ et $\Gamma_2$
sont isomorphes, et, s'ils le sont, à décrire la classe de leurs isomorphismes.
(Un {\em isomorphisme} $\pi:\Gamma_1\to \Gamma_2$ est
une bijection $\pi$ de l'ensemble de sommets de
$\Gamma_1$ vers celui de $\Gamma_2$ telle que $\pi(\Gamma_1) = \Gamma_2$.)
Une solution permettrait, par exemple, de trouver une molécule dans une
base de données.

Le problème de l'isomorphisme de graphes se réduit en temps polynomial
au problème de
l'isomorphisme de chaînes, de la façon suivante. Supposons sans perte
de généralité que $\Gamma_1$ et $\Gamma_2$ ont le même ensemble de sommets
$V$. Alors, nous pouvons définir $\Omega$ comme l'ensemble des paires
d'éléments de $V$ (ordonnés ou non ordonnés, suivant que nos graphes
sont orientés ou pas). La chaîne $\mathbf{x}_i$, $i=1,2$, est définie
comme suit: pour la paire $a=\{v_1,v_2\}$ (ou $a=(v_1,v_2)$, si nos graphes
sont orientés), la valeur de $\mathbf{x}_i(a)$ est $1$ s'il y a une arête
entre $v_1$ et $v_2$ en $\Gamma_1$, et $0$ dans le cas contraire. Soit $G$
l'image de l'homomorphisme $\iota:\Sym(V)\to \Sym(\Omega)$ définie par
$\sigma^\iota(\{v_1,v_2\}) = \{\sigma(v_1),\sigma(v_2)\}$, où
$\sigma^\iota = \iota(\sigma)$. Alors $\iota$ induit une bijection entre
la classe des isomorphismes de
$\Gamma_1$ à $\Gamma_2$ et la classe $\Iso_G(\mathbf{x}_1,\mathbf{x}_2)$.

\begin{theo}[Babai]
  Le problème de l'isomorphisme de chaînes $\Omega\to \Sigma$
  peut être résolu en temps quasi-polynomial en 
le nombre d'éléments du domaine $\Omega$.
\end{theo}

En novembre 2015, Babai a annoncé une solution en temps
quasipolynomial, avec un algorithme explicite.
La préparation de cet exposé m'a conduit à trouver une erreur non triviale
dans l'analyse du temps, mais Babai a réussi à le réparer en simplifiant
l'algorithme. La preuve est maintenant correcte.

\begin{coro}[Babai]
  Le problème de l'isomorphisme de graphes peut être résolu en temps
 quasi-polynomial en le nombre de sommets.
\end{coro}

Notre référence principale sera
\cite{Ba}; nous nous servirons aussi de la version courte \cite{Ba2}.
Nous essayerons d'examiner la preuve de la façon la plus détaillée possible
dans un exposé de ce format, en partie pour aider à éliminer tout doute
qui pourrait rester sur la forme actuelle du résultat.


La meilleure borne générale connue antérieurement pour le temps requis par
le problème de l'isomorphisme de graphes, due à Luks \cite{BKL},
était $\exp(O(\sqrt{n \log n}))$,


\begin{center}
  * * *
\end{center}

L'usage de la {\em canonicité} joue un rôle crucial dans la stratégie
de Babai. Comme dans la théorie de catégories, voire dans l'usage courant,
un choix est {\em canonique} s'il est fonctoriel. La situation typique
pour nous sera la suivante: un groupe $G<\Sym(\Omega)$ agit sur $\Omega$,
et donc sur $\Sigma^\Omega$; il agit aussi sur un autre ensemble $S$,
et donc aussi sur les applications $S\to \mathscr{C}$,
où $\mathscr{C}$ est un ensemble fini. 
Une application $S\to \mathscr{C}$ s'appelle un {\em coloriage};
l'ensemble $\mathscr{C}$
s'appelle l'ensemble de {\em couleurs}.
Un choix {\em canonique} (en relation à $G$) d'un coloriage de $\Omega$
pour chaque chaîne $\mathbf{x}\in \Sigma^\Omega$ est une application qui va de
$\Sigma^\Omega$ aux coloriages et qui commute avec l'action de $G$.


En particulier, un choix canonique peut être un outil pour détecter des
non-isomorphismes: si les coloriages $C(\mathbf{x})$ et
$C(\mathbf{y})$ induits canoniquement par
$\mathbf{x}$ et $\mathbf{y}$ ne sont pas isomorphes l'un à l'autre 
-- par exemple, s'ils ont un nombre différent d'éléments vermeils -- alors
$\mathbf{x}$ et $\mathbf{y}$ ne sont pas isomorphes l'un à l'autre.
Même quand il y a des isomorphismes dans $G$ qui envoient $C(\mathbf{x})$ sur
$C(\mathbf{y})$, la classe $\Iso_G(C(\mathbf{x}),C(\mathbf{y}))$
de tels isomorphismes sert à délimiter la classe d'isomorphismes
$\Iso_G(\mathbf{x},\mathbf{y})$
de $\mathbf{x}$ à $\mathbf{y}$, puisque cette dernière est
forcément un sous-ensemble
de $\Iso_G(C(\mathbf{x}),C(\mathbf{y}))$.

La preuve assimile aussi plusieurs idées développées lors d'approches
antérieures au problème. La première étape de la procédure consiste à
essayer de suivre ce qui est en essence
l'algorithme de Luks \cite{Lu}. Si cet algorithme s'arrête, c'est parce
qu'il s'est heurté contre un quotient
$H_1/H_2$ isomorphe à $\Alt(\Gamma)$, où $H_2\triangleleft H_1 < G$
et $\Gamma$ est plutôt grand.

Notre tâche majeure consiste à étudier ce qui se passe à ce moment-là.
La stratégie
principale sera de chercher à colorier $\Gamma$ d'une façon qui
dépend canoniquement de $\mathbf{x}$. Cela limitera les automorphismes
et isomorphismes possibles à considérer. Par exemple, si la moitié
de $\Gamma$ est coloriée en rouge et l'autre en noir, le groupe
d'automorphismes possibles se réduit à
$\Sym(|\Gamma|/2)\times \Sym(|\Gamma|/2)$. Un coloriage similaire induit
par $\mathbf{y}$ limite les isomorphismes aux applications qui alignent
les deux coloriages. Nous trouverons toujours des coloriages qui nous
aident, sauf quand certaines structures ont une très grande symétrie,
laquelle, en revanche, permettra une descente à $\Omega$ considérablement
plus petit. Cette double récursion -- réduction du groupe $H_1/H_2$ ou
descente à des chaînes considérablement plus courtes -- résoudra le problème.
\section{Fondements et travaux précédents}




En suivant l'usage courant pour les groupes de permutations, nous
écrirons $r^g$ pour l'élément $g(r)$ auquel $g\in \Sym(\Omega)$ envoie
$r\in \Omega$. \'Etant donnés une chaîne $\mathbf{x}:\Omega\to \Sigma$ et
un élément
$g\in \Sym(\Omega)$, nous définissons $\mathbf{x}^g:\Omega\to \Sigma$ par
$\mathbf{x}^g(r) = \mathbf{x}\left(r^{g^{-1}}\right)$.

Par contre, nous
écrivons $\Omega^k$ pour l'ensemble des $\vec{x}=(x_1,\dotsc,x_k)$
avec l'action à gauche donnée par
$(\phi(\vec{x}))_r = \vec{x}_{\phi(r)}$. L'idée est que ceci est défini non pas
seulement pour $\phi$ une permutation, mais pour toute application
$\phi:\{1,\dotsc,k\}\to \{1,\dotsc,k\}$, même non injective. Nous appelons
les éléments de $\Omega^k$ {\em tuples} plutôt que {\em chaînes}.

\subsection{Algorithmes de base}

\subsubsection{Schreier-Sims}

Plusieurs algorithmes essentiels se basent sur une idée de Schreier \cite{Sch}.
Il a remarqué que, pour tout sous-groupe 
$H$ d'un groupe $G$
et tout sous-ensemble $A\subset G$ qui engendre $G$ et contient des
représentants de toutes les classes de $H$ dans $G$,
\[A' = A A A^{-1}\cap H =
\left\{\sigma_1 \sigma_2 \sigma_3^{-1} : \sigma_i\in A\right\} \cap H\]
est un ensemble de générateurs de $H$.

L'étape suivante est celle de Sims \cite{Si1}, \cite{Si2}, qui a montré l'utilité
de travailler avec un groupe de permutations $G<\Sym(\Omega)$,
$\Omega = \{x_1,\dotsc, x_n\}$,
en
termes d'une {\em chaîne de stabilisateurs}
\[G =G_0> G_1 > G_2 > \dotsc >G_{n-1} = \{e\},\]
où $G_k = G_{(x_1,x_2,\dotsc,x_k)} =
\{g \in G: \forall 1\leq i\leq k\;\; x_i^g = x_i \}$ ({\em stabilisateur
  de points}).

L'algorithme de Schreier-Sims (Algorithme 1; description basée sur
\cite[\S 1.2]{Lu}) construit des ensembles $C_i$ de représentants
de $G_i/G_{i+1}$ tels que
$\cup_{i\leq j<n-1} C_j$ engendre $G_i$ pour tout $0\leq i<n-1$.
Le temps pris par l'algorithme
est $O(n^5 + n^3 |A|)$, où $A$ est l'ensemble de générateurs de $G$ qui
nous est donné: la fonction \textsc{Filtre} prend $O(n)$ de temps,
et tout $g$ pour lequel elle est appelée satisfait
$g \in A C \cup C A \cup C^2$, où $C$ est la valeur de
$\cup_i C_i$ à la fin de la procédure. Bien sûr, $|C|\leq n (n+1)/2$.

Grâce à l'algorithme lui-même, nous pourrons toujours supposer
que nos ensembles de générateurs sont de taille $O(n^2)$.
Le temps pris par l'algorithme est donc $O(n^5)$.\footnote{Nous supposons que l'ensemble
  de générateurs initial, spécifiant le groupe $G$ du problème,
  est de taille $O(n^C)$, $C$ une constante. Le temps pris par la
première utilisation de l'algorithme est donc
$O\left(n^{\max(5,3+C)}\right)$.}

\begin{algorithm}
  \caption{Schreier-Sims: construction d'ensembles $C_i$}\label{alg:schrsims}
  \begin{algorithmic}[1]
    \Function{SchreierSims}{$A$, $\vec{x}$}
    \Comment{$A$ engendre $G<\Sym(\{x_1,\dotsc,x_n\})$}
    \Ensure{$\cup_{i\leq j<n-1} C_j$ engendre $G_i$
    et $C_i\mapsto G_i/G_{i+1}$ est injectif\; $\forall i\in \{0,1,\dotsc,n-2\}$}
    \State{$C_i\gets \{e\}$ pour tout $i\in \{0,1,\dotsc,n-2\}$}
    \State{$B\gets A$}
    \While{$B\ne \emptyset$}
    \State{Choisir $g\in B$ arbitraire, et l'enlever de $B$}
    \State{$(i,\gamma) \gets \text{\textsc{Filtrer}($g$, $(C_i)$, $\vec{x}$)}$}
    \If{$\gamma\ne e$}
    \State{ajouter $\gamma$ à $C_i$}
    \State{$B\gets B \cup \bigcup_{j\leq i} C_j \gamma \cup
          \bigcup_{j\geq i} \gamma C_j$}  
    \EndIf
    \EndWhile
    \State \Return $(C_i)$ 
    \EndFunction
    \Function{Filtrer}{$g$, $(C_i)$, $\vec{x}$}
    \Comment{retourne $(i,\gamma)$ tel que   $\gamma\in G_i$,
      $g\in C_0 C_1\dotsb C_{i-1} \gamma$
    }
    \Require 
    $C_i\subset G_i$ et
  $C_i \to G_i/G_{i+1}$ injectif\; $\forall i\in \{0,1,\dotsc,n-2\}$
  \Ensure 
  $g\notin C_0 C_1 \dotsb C_i G_{i+1}$ sauf si $(i,\gamma) = (n-1,e)$
  \State{$\gamma \gets g$}
  \For{$i=0$   \textbf{jusqu'à} $n-2$}
  \If{$\exists h\in C_i$ tel quel $x_i^{h} = x_i^{g}$}
  \State{$\gamma\gets h^{-1} \gamma$} 
  \Else{}
  \State{\Return $(i,\gamma)$}
    \EndIf
    \EndFor 
    \State \Return $(n-1,e)$
    \EndFunction
  \end{algorithmic}
\end{algorithm}

Une fois les ensembles $C_i$ construits, il devient possible d'accomplir
plusieurs tâches essentielles rapidement.

\begin{exer}\label{ex:fhl}
  Montrer comment accomplir les tâches suivantes en temps polynomial,
  étant donné un groupe $G<\Sym(\Omega)$, $|\Omega|=n$:
  \begin{enumerate}
  \item   Déterminer si un élément $g\in \Sym(\Omega)$ est dans $G$.
  \item \'Etant donnés un homomorphisme $\phi: G\to \Sym(\Omega')$,
    $|\Omega'|\ll |\Omega|^{O(1)}$,
    et un sous-groupe $H<\Sym(\Omega')$, décrire $\phi^{-1}(H)$.
  \item\label{it:richt} \cite{FHL} Soit $H<G$ avec $\lbrack G:H\rbrack\ll n^{O(1)}$.
    \'Etant donné un test qui détermine en temps polynomial si un élément
    $g\in G$ appartient à $H$, décrire $H$. {\em Astuce:} travailler avec
    $G > H > H_1>H_2 >\dotsc$ à la place de $G =G_0> G_1 > G_2 > \dotsc$.
  \end{enumerate}
\end{exer}

 Ici, comme toujours, {\og décrire \fg} veut dire
 {\og trouver un ensemble de générateurs\fg}, et un groupe nous est {\og donné\fg} si
 un tel ensemble nous est donné.

L'algorithme de Schreier-Sims
décrit le stabilisateur de points $G_{(x_1,\dotsc,x_k)}$
 pour $x_1,\dotsc,x_k \in \Omega$ arbitraires. Par contre, nous ne pouvons
 pas demander allègrement un ensemble de générateurs d'un
 {\em stabilisateur d'ensemble} $G_{\{x_1,\dotsc,x_k\}} = \{g\in G: \{x_1^g,\dotsc,x_k^g\} =
 \{x_1,\dotsc,x_k\}\}$ pour $G$, $x_i$ arbitraires: faire ceci serait
 équivalent à résoudre le problème de l'isomorphisme lui-même.

 \subsubsection{Orbites et blocs}\label{subs:orbl}

 Soit donné, comme toujours,
 un groupe de permutations $G$ agissant sur un ensemble fini $\Omega$. 
 Le domaine $\Omega$ est l'union disjointe des {\em orbites}
 \mbox{$\{x^g: g\in G\}$} de~$G$. Ces orbites peuvent être
 déterminées en temps polynomial\footnote{Pour être précis:
   $O\left(|\Omega|^{O(1)} + |A| |\Omega|\right)$, où $A$ est la taille de l'ensemble de générateurs de $G$ qui nous est donné. Nous omettrons toute mention de cette taille par la suite, puisque, comme nous l'avons déjà dit, nous pouvons
 la garder toujours sous contrôle.}
 en $|\Omega|$. Ceci est un exercice simple. La tâche se réduit à celle -- simple elle aussi -- de trouver les composantes connexes d'un graphe.

 Supposons que l'action de $G$ soit transitive. (Il y a donc une seule orbite.)
 Un {\em bloc} de $G$ est un sous-ensemble $B\subset \Omega$,
 $B\notin \{\emptyset, \Omega\}$,
 tel que, pour $g,h\in G$ quelconques, soit $B^g = B^h$, soit $B^g\cap B^h =
 \emptyset$.
 La collection $\{B^g: g\in G\}$ ({\em système de blocs})
 pour $B$ donné partitionne $\Omega$. L'action de $G$ est {\em primitive}
 s'il n'y a pas de blocs de taille $>1$; autrement, elle s'appelle
 {\em imprimitive}.
 Un système de blocs est {\em minimal}\footnote{Pour paraphraser \cite[\S 1.1]{Lu}: il faut avouer qu'un
   tel système pourrait s'appeler plutôt maximal. La taille des blocs est
   maximale, leur nombre est minimal.}
 si l'action de $G$ sur lui est primitive.
 
 Voyons comment déterminer si l'action de $G$ est primitive,
 et, s'il ne l'est pas, comment
 trouver un système de blocs de taille $>1$. En itérant la procédure, nous
 obtiendrons un système de blocs minimal en temps polynomial.
 (Nous suivons \cite{Lu}, qui cite \cite{Si1}.)

 Pour $a, b\in \Omega$
 distincts, soit $\Gamma$ le graphe avec $\Omega$ comme son ensemble de sommets et
 l'orbite $\{\{a,b\}^g : g\in G\}$ comme son ensemble d'arêtes.
 La composante connexe qui contient $a$ et $b$ est le bloc le plus petit
 qui contient $a$ et $b$. (Si $\Gamma$ est connexe, alors le {\og bloc\fg} est
 $\Omega$.) L'action de $G$ est imprimitive ssi $\Gamma$ est non connexe pour
 un $a$ arbitraire et au moins un $b$; dans ce cas-là, nous obtenons un bloc
 qui contient $a$ et $b$, et donc tout un système de blocs de taille $>1$. 

 Un dernier mot: si $G<\Sym(\Omega)$, nous disons que $G$ est {\em transitif},
 voire {\em primitif}, si son action sur $\Omega$ l'est.
 
 \subsection{Luks: le cas de groupes avec facteurs d'ordre borné}\label{subs:luks}

 Luks a montré comment résoudre le problème de l'isomorphisme de graphes
 en temps polynomial dans le cas spécial de graphes
 de degré borné. (Le {\em degré}, ou {\em valence}, d'un sommet dans
 un graphe non orienté est le nombre d'arêtes qui le contiennent.)
 Il réduit ceci au problème de décrire le groupe d'automorphismes
 de chaînes dans le cas
 d'un groupe $G$ tel que tout facteur de composition de $G$
 -- c'est-à-dire, tout quotient dans une suite principale (Jordan-H\"older)
 de $G$ -- est borné. Le processus de réduction, élégant et loin d'être
 trivial, ne nous concerne pas ici. Voyons plutôt comment Luks résout ce
 cas du problème de l'isomorphisme de chaînes.

 Nous suivrons la notation de \cite{Ba}, même si les idées viennent de \cite{Lu}.
 \begin{defi}
   Soient $K\subset \Sym(\Omega)$ et $\Delta\subset \Omega$ (la {\og fenêtre\fg}).
   L'ensemble d'{\em isomorphismes partiels} $\Iso_K^\Delta$ est 
   \[\Iso_K^\Delta(\mathbf{x},\mathbf{y}) = \{\tau\in K: \mathbf{x}(x) =
\mathbf{y}(x^\tau)\;\;\; \forall x\in \Delta\}.\]
   L'ensemble d'{\em automorphismes partiels} $\Aut_K^\Delta(\mathbf{x})$
   est égal à $\Iso_K^\Delta(\mathbf{x},
   \mathbf{x})$.
 \end{defi}
 $\Iso_K^\Delta$ est donc l'ensemble de toutes les permutations $g\in K$ qui
 envoient $\mathbf{x}$ sur $\mathbf{y}$ -- au moins à en juger par
 ce qui peut se voir par la fenêtre $\Delta$. Nous travaillerons en général
 avec $K$ de la forme $H \pi$, où $H$ laisse $\Delta$ invariante
 (en tant qu'ensemble).
 
 Il est clair que, pour $K,K_1,K_2\subset \Sym(\Omega)$ et $\sigma\in
 \Sym(\Omega)$,
 \begin{equation}\label{eq:udu1}
 \Iso_{K \sigma}^\Delta(\mathbf{x},\mathbf{y}) =
 \Iso_K^\Delta\left(\mathbf{x},\mathbf{y}^{\sigma^{-1}}\right) \sigma,\end{equation}
 \begin{equation}\label{eq:udu2}\Iso_{K_1\cup K_2}^\Delta(\mathbf{x},\mathbf{y}) =
 \Iso_{K_1}^\Delta(\mathbf{x},\mathbf{y}) \cup
 \Iso_{K_2}^\Delta(\mathbf{x},\mathbf{y}).
 \end{equation}
 
 Il est aussi clair que, si $G$ est un sous-groupe de $\Sym(\Omega)$ et $\Delta$ est
 invariant sous $G$, alors $\Aut_G(\mathbf{x})$ est un sous-groupe de $G$, et, pour tout $\sigma\in \Sym(\Omega)$,
 $\Iso_{G \sigma}(\mathbf{x},\mathbf{y})$ est soit vide,
 soit une classe à
 droite de la forme $\Aut_G(\mathbf{x}) \tau$, $\tau \in \Sym(\Omega)$.
 Soient $\Delta_1, \Delta_2\subset \Omega$, $\Delta_1$ invariant sous $G$.
 Pour $G' = \Aut_G(\mathbf{x})$ et $\sigma, \tau$ tels que
 $\Iso_{G \sigma}^{\Delta_1}(\mathbf{x},\mathbf{y}) = G' \tau$,
 \begin{equation}\label{eq:udu3}
   \Iso_{G \sigma}^{\Delta_1\cup \Delta_2}(\mathbf{x},\mathbf{y}) =
 \Iso_{G' \tau}^{\Delta_2}(\mathbf{x},\mathbf{y}) =
 \Iso_{G'}^{\Delta_2}\left(\mathbf{x},\mathbf{y}^{\tau^{-1}}\right) \tau,\end{equation}
 où la deuxième équation est une application de (\ref{eq:udu1}). 
Babai appelle (\ref{eq:udu3}) la {\em règle de la chaîne}.


L'énoncé suivant n'utilise pas la classification de groupes finis simples.
\begin{theo}[\cite{BCP}\footnote{\`A vrai dire, \cite[Thm 1.1]{BCP} est plus général que ceci; par exemple, des
facteurs abéliens arbitraires (non bornés) sont admis. Cela donne une généralisation du Théorème \ref{thm:luxor}.}]\label{thm:bcp}
  Soit $G<\Sym(\Omega)$ un groupe primitif. Soit $n = |\Omega|$.
  Si tout facteur de composition de $G$ est d'ordre $\leq k$, alors
  $|G|\leq n^{O_k(1)}$.
\end{theo}
Ici, comme d'habitude, $O_k(1)$ désigne une quantité qui dépend seulement de $k$.

\begin{theo}[Luks \cite{Lu}]\label{thm:luxor}
  Soient $\Omega$ un ensemble fini et $\mathbf{x}, \mathbf{y}:\Omega\to
  \Sigma$ deux chaînes.
  Soit donné un groupe $G<\Sym(\Omega)$ tel que
  tout facteur de composition de $G$ est d'ordre $\leq k$.
  Il est possible de déterminer
  $\Iso_G(\mathbf{x},\mathbf{y})$ en temps polynomial en $n = |\Omega|$.
\end{theo}
  \noindent{\sc Preuve} ---
  {\bf Cas 1: $G$ non transitif.}
  Soit $\Delta_1\subsetneq \Omega$,
  $\Delta_1\ne \emptyset$, $\Delta_1$ stable sous l'action de $G$.
  Définissons $\Delta_2 = \Omega\setminus \Delta_1$.
  Alors, par (\ref{eq:udu3}), il suffit
  de calculer $\Iso_G^{\Delta_1}(\mathbf{x},\mathbf{y})$ (égal
  à une classe que nous notons $G' \tau$) et $\Iso_{G'}^{\Delta_2}(\mathbf{x},\mathbf{y}')$ pour $\mathbf{y}' = \mathbf{y}^{\tau^{-1}}$.
  Or, pour déterminer $\Iso_G^{\Delta_1}(\mathbf{x},\mathbf{y})$, nous
  déterminons, de façon récursive,
  $\Iso_G\left(\mathbf{x}|_{\Delta_1},\mathbf{y}|_{\Delta_1}\right)$, puis, par Schreier-Sims, le stabilisateur de points $G_{(\Delta_1)}$.
  De la même manière, déterminer $\Iso_{G'}^{\Delta_2}(\mathbf{x},\mathbf{y}')$ pour $\mathbf{y}' = \mathbf{y}^{\tau^{-1}}$ se réduit à déterminer
  le groupe d'isomorphismes (dans un groupe $G'$)
  entre deux chaînes de longueur $|\Delta_2|$. Comme
  $|\Delta_1| + |\Delta_2| = n$ et Schreier-Sims prend du temps $O(n^5)$,
  tout va bien. (La comptabilité est laissée au lecteur.)

  {\bf Cas 2: $G$ transitif.} Soit $N$ le stabilisateur d'un système
  de blocs minimal pour $G$; donc, $G/N$ est primitif. Par le Théorème
  \ref{thm:bcp}, $|G/N|\leq m^{O_k(1)}$, où $m$ est le nombre de blocs.
  Or, pour $\sigma_1,\dotsc,\sigma_{\ell}$ ($\ell = |G/N|$)
tels que $G = \cup_{1\leq i\leq \ell} N \sigma_i$,
  \begin{equation}\label{eq:rulu}
    \Iso_G(\mathbf{x},\mathbf{y}) = \Iso_{\cup_i N \sigma_i}(\mathbf{x},
  \mathbf{y}) = \bigcup_{1\leq i\leq \ell} \Iso_{N \sigma_i}(\mathbf{x},
  \mathbf{y}) = \bigcup_{1\leq i\leq \ell} \Iso_{N}(\mathbf{x},
  \mathbf{y}^{\sigma_i^{-1}}) \sigma_i
  \end{equation}
  par (\ref{eq:udu1}) et (\ref{eq:udu2}). Comme les orbites de $N$ sont
  contenues dans les blocs, qui sont de taille $n/m$,
  déterminer $\Iso_N(\mathbf{x},\mathbf{y}_i)$
  ($\mathbf{y}_i = \mathbf{y}^{\sigma_i^{-1}}$)
  se réduit -- par la règle (\ref{eq:udu3}) --
  à déterminer les groupes d'isomorphismes de $m$ paires
  de chaînes de longueur $n/m$.
  Nous avons donc réduit le problème à la solution de
  $\ell\cdot m = m^{O_k(1)}$
  problèmes pour des chaînes de longueur $n/m$.
  
  Le pas final consiste à faire l'union de classes
   en (\ref{eq:rulu}). Nous avons une description de
  chaque $\Iso_N(\mathbf{x},\mathbf{y}_i)$, soit comme l'ensemble vide,
  soit comme une classe à droite $H \tau_i$
  du groupe $H=\Aut_N(\mathbf{x})$, dont nous avons
  trouvé une description, c'est-à-dire un ensemble de générateurs $A$. Alors
  \[\begin{aligned}   \Iso_G(\mathbf{x},\mathbf{y}) &=
\bigcup_{1\leq i\leq \ell} \Iso_{N}(\mathbf{x},\mathbf{y}_i) \sigma_i =
\bigcup_{1\leq i\leq \ell} H \tau_i \sigma_i\\
&= \left\langle A \cup \left\{
\tau_i \sigma_i (\tau_1 \sigma_1)^{-1}: 1\leq i\leq \ell\right\}\right\rangle
  \tau_1 \sigma_1.
  \end{aligned}\]
\qed

Nous aurions pu éviter quelques appels à Schreier-Sims en travaillant toujours
avec des isomorphismes partiels, mais cela a peu d'importance
qualitative.

\subsection{Relations, partitions, configurations}\label{subs:secf1}

Soit $\mathscr{C}$ ({\og couleurs\fg}) un ensemble fini que nous pouvons supposer
ordonné (disons, de rouge à violet). Une {\em relation $k$-aire} sur un
ensemble fini $\Gamma$ est un sous-ensemble $R\subset \Gamma^k$. Une
{\em structure (relationnelle) $k$-aire} est une paire
$\mathfrak{X} = (\Gamma, (R_i)_{i\in \mathscr{C}})$,
où, pour chaque $i\in \mathscr{C}$, $R_i$ est une relation $k$-aire sur
$\Gamma$. Si les $R_i$ sont tous non vides et partitionnent $\Gamma^k$,
nous disons que $\mathfrak{X}$ est une {\em structure de partition $k$-aire}.
Dans ce cas-là, nous pouvons décrire $\mathfrak{X}$ par une fonction
$c:\Gamma^k\to \mathscr{C}$ qui assigne à chaque $\vec{x}\in \Gamma^k$ l'indice $i$ de la relation $R_i$ à laquelle il appartient. Nous disons que
$c(\vec{x})$ est la couleur de $\vec{x}$.

Un isomorphisme entre deux structures $k$-aires $\mathfrak{X} =
(\Gamma, (R_i)_{i\in \mathscr{C}})$ et $\mathfrak{X}' =
(\Gamma', (R_i')_{i\in \mathscr{C}})$ est une bijection $\Gamma\to \Gamma'$
qui envoie $R_i$ à $R_i'$ pour chaque $i$.
Il est possible de construire un foncteur $F_1$ qui envoie
chaque structure $k$-aire $\mathfrak{X}$ sur $\Gamma$
à une structure de partition $k$-aire
$F_1(\mathfrak{X})$ sur $\Gamma$; qui plus est, $\Iso(\mathfrak{X},\mathfrak{Y})
=\Iso(F_1(\mathfrak{X}),F_1(\mathfrak{Y}))$.
La procédure est plutôt triviale; nous la détaillons (Algorithme
\ref{alg:struraff}) pour montrer ce qu'{\em indexer} veut dire. Cela nous
permet de ne pas utiliser plus de $\min\left(|\Gamma|^k,2^{|\mathscr{C}|}\right)$ couleurs, où $n= |\Omega|$, tout en gardant
leur signification en termes des couleurs originales $\mathscr{C}$. 
Le temps pris pour calculer $F_1(\mathfrak{X})$ est $O(|\mathscr{C}| |\Gamma|^{O(k)})$.
Nous ne nous occupons pas des détails d'implémentation de la collection
de tuples $\mathscr{I}$,
mais il peut s'agir tout simplement d'une liste ordonnée lexicographiquement;
dans ce cas, $|\Gamma|^{O(k)}$ est $|\Gamma|^{2k}$. (Dans la réalité,
$\mathscr{I}$ serait implémentée avec du
{\em hachage}, ce qui n'est que l'art de bien organiser une bibliothèque.)

\begin{algorithm}
  \caption{Raffinement d'une structure de relation. Indexeur.}\label{alg:struraff}
  \begin{algorithmic}[1]
    \Function{$F_1$}{$\Gamma$,$k$,$\mathscr{C}$,$(R_i)_{i\in \mathscr{C}}$}
    \State $\mathscr{I} \gets \emptyset$
    \For{$\vec{x}\in \Gamma^k$}
    \State $\mathit{a} \gets \{i\in \mathscr{C}: \vec{x}\in R_i\}$
    \State $c(\vec{x}) \gets \text{\textsc{Indexeur}}(\mathscr{I},\mathit{a})$
    \EndFor
    \State \Return $(\mathscr{I}, c)$ \Comment{retourne $c:\Gamma^k\to \mathscr{C}'$}\\
    \Comment{$\mathscr{C}'$ est l'ensemble d'indices de $\mathscr{I}$;
      $\mathscr{I}$
      explique $\mathscr{C}'$ en termes de $\mathscr{C}$}
    \EndFunction
    \Function{Indexeur}{$\mathscr{I}$,$\mathit{a}$}
    \Comment{$\mathscr{I}$ est une collection modifiable}
    \If{$\mathit{a}$ n'est pas dans $\mathscr{I}$}
    \State ajouter $\mathit{a}$ à $\mathscr{I}$
    \EndIf
    \State \Return indice de $\mathit{a}$ dans $\mathscr{I}$
        \EndFunction
  \end{algorithmic}
\end{algorithm}

Un élément $\vec{x}\in \Gamma^k$ définit une relation d'équivalence
$\rho(\vec{x})$ sur $\{1,\dotsc,k\}$: $i\sim j$ ssi $x_i = x_j$.
Le monoïde $\mathfrak{M}(S)$ ($S$ un ensemble) consiste en les
applications $S\to S$, avec la composition comme opération.
 \begin{defi}\label{def:ukuru}
   Une structure de partition $k$-aire $\mathfrak{X} = (\Gamma,c)$ est dite
   {\em configuration $k$-aire} si
   \begin{enumerate}
   \item\label{it:klara} Pour tous $\vec{x}, \vec{y}\in \Gamma^k$, si
     $c(\vec{x}) = c(\vec{y})$, alors $\rho(\vec{x}) = \rho(\vec{y})$.
   \item\label{it:klarb} Il y a un homomorphisme de monoïdes
     $\eta: \mathfrak{M}(\{1,\dotsc,k\}) \to \mathfrak{M}(\mathscr{C})$
     tel que, pour tout $\tau\in \mathfrak{M}(\{1,\dotsc,k\})$,
     $c\left(\tau(\vec{x})\right) = \tau^\eta(c(\vec{x}))$ pour tout
     $\vec{x}\in \Gamma^k$.
     \end{enumerate}
 \end{defi}
   Alors, par exemple, pour $k=2$, (\ref{it:klara}) veut dire que la couleur
   de $\vec{x}=(x_1,x_2)$ {\og sait\fg} si $x_1=x_2$ ou pas, dans le sens où,
   si nous connaissons $c(\vec{x})$, alors nous savons si $x_1=x_2$ ou pas.
   De la même façon,
   (\ref{it:klarb}) nous indique que la couleur de $\vec{x}$ connaît
   les couleurs de $(x_2,x_1)$, $(x_1,x_1)$ et $(x_2,x_2)$.

Nous pouvons définir un foncteur $F_2$ qui envoie
chaque structure de partition
\mbox{$k$-aire $\mathfrak{X}$} sur $\Gamma$ à une configuration
$k$-aire; comme pour $F_1$, le fait que $F_2(\mathfrak{X})$ est un
raffinement de $\mathfrak{X}$ implique que
$\Iso(\mathfrak{X},\mathfrak{Y}) = \Iso(F_2(\mathfrak{X}),F_2(\mathfrak{Y}))$.
La procédure pour calculer $F_2$ est très similaire à celle pour calculer
$F_1$ (Algorithme \ref{alg:struraff}). Au lieu d'assigner à
$\vec{x}$ la couleur $\{i\in \mathscr{C}: \vec{x}\in R_i\}$, nous lui assignons
la couleur $\left(\rho(\vec{x}),
(c(\phi(\vec{x})))_{\phi \in \mathfrak{M}(\{1,\dotsc,k\})}\right)$.

Il est aisé de voir que $F_2(\mathfrak{X})$ est le raffinement le plus
grossier d'une structure de partition
$\mathfrak{X}$ qui est une configuration, de la même manière que
$F_1(\mathfrak{X})$ est le raffinement le plus grossier d'une structure
$\mathfrak{X}$ qui est
une structure de partition.

\begin{defi}
  Soit $\mathfrak{X} = (\Gamma,c)$, $c:\Gamma^k\to \mathscr{C}$,
  une structure de partition $k$-aire.
  Pour $1\leq l\leq k$, nous définissons $c^{(l)}:\Gamma^l\to
  \mathscr{C}$ comme suit:
  \[c^{(l)}(\vec{x}) = c(x_1,x_2,\dotsc,x_l, x_l, \dotsc x_l).\]
  La structure de partition $l$-aire
  $\mathfrak{X}^{(l)} = \left(\Gamma,c^{(l)}\right)$ est dite
  le {\em ($l$-)squelette} de $\mathfrak{X}$.
\end{defi}
La chaîne vide sera viride.

\begin{exer}
Tout squelette d'une configuration est une configuration.
\end{exer}
Ici le fait que l'axiome (\ref{it:klarb}) dans la définition de {\em configuration}
soit valable même pour~$\eta$ non injectif est crucial.

Pour $\mathfrak{X} = (\Gamma,c)$ une structure de partition
et $\Gamma'\subset \Gamma$, la {\em sous-structure induite}
$\mathfrak{X}\lbrack \Gamma'\rbrack$ est la structure
$(\Gamma',c|_{\Gamma'})$ définie par la restriction de $c$ à $\Gamma'$.
Il est clair que, si $\mathfrak{X}$ est une configuration, alors
$\mathfrak{X}\lbrack \Gamma'\rbrack$ l'est aussi.

\begin{center}
  * * *
\end{center}

Il ne faut pas confondre une structure de partition ({\em partition structure})
avec ce que nous appellerons un {\em découpage} ({\em colored partition}).
Un découpage d'un ensemble $\Gamma$
est un coloriage de $\Gamma$ supplémenté d'une
partition de chaque classe de couleur. (Une {\em classe de couleur}
est l'ensemble de sommets d'une couleur donnée.) Un découpage est dit
{\em admissible} si chaque ensemble $B$ dans chaque partition est de taille $\geq 2$. Pour $\alpha<1$,
un {\em $\alpha$-découpage} est un découpage admissible tel que
$|B|\leq \alpha |\Gamma|$ pour chaque $B$.

Un découpage est une structure plus fine que le coloriage qu'il raffine,
mais moins fine que la structure que nous obtiendrions si nous donnions à chaque
élément de chaque partition une couleur différente. Un automorphisme ou
isomorphisme d'un découpage doit préserver les couleurs de celui-ci, mais
pourrait permuter les ensembles de la même taille qui appartiennent à
la partition d'une couleur.
Comme les ensembles de taille différente ne peuvent, évidemment, être
permutés, il est clair que
nous pouvons supposer sans perte de généralité que toute couleur est
partitionnée en ensembles de la même taille. Nous ajoutons ceci à la
définition de $\alpha$-découpage à partir de maintenant.

\subsection{Configurations cohérentes $k$-aires}\label{subs:concoh}

Pour $\vec{x}\in \Gamma^k$, $z\in \Gamma$ et $1\leq i\leq k$, nous
définissons $\vec{x}^i(z) \in \Gamma^k$ comme suit:
\[\vec{x}^i(z) = (x_1, x_2, \dotsc,x_{i-1},z,x_{i+1},\dotsc, x_k).\]
\begin{defi}\label{def:coh}
  Une {\em configuration cohérente $k$-aire} $\mathfrak{X}=(\Gamma,c)$
  est une configuration $k$-aire ayant la propriété suivante: il y a une
  fonction $\gamma:\mathscr{C}^k\times \mathscr{C} \to \mathbb{Z}_{\geq 0}$ telle
  que, pour $\vec{k}\in \mathscr{C}^k$ et $j\in \mathscr{C}$ arbitraires et
  tout $\vec{x}\in \Gamma^k$ tel que $c(\vec{x})=j$,
  \[|\{z\in \Gamma: c(\vec{x}^i(z)) = k_i\; \forall 1\leq i\leq k\}| = \gamma(\vec{k},j).\]
  Les valeurs $\gamma(\vec{k},j)$ sont appelées {\em nombres d'intersection}
  de $\mathfrak{X}$.
\end{defi}
Une configuration cohérente est dite {\em classique} si $k=2$.

\begin{rema}
  Les configurations cohérentes classiques ont été
  introduites par Higman \cite{Hi}. Les premiers exemples étaient
  du type {\em schurien}: une configuration est
  {\em schurienne} si elle est la partition
  de $\Gamma^2$ dans ses orbites ({\og orbitales\fg}) sous l'action d'un groupe $G<\Sym(\Gamma)$.
\end{rema}

\begin{defi}
  Si une configuration cohérente classique n'a que deux couleurs, une pour
  $\{(x,x) : x\in \Gamma\}$ et l'autre pour son complément, la configuration
  est dite {\em une clique}, ou {\em triviale}.
\end{defi}

\begin{exer}
Tout squelette d'une configuration cohérente est cohérent.
\end{exer}
Encore une fois, l'axiome (\ref{it:klarb}) des configurations joue un rôle clé.

\begin{exer}\label{ex:rescoh}
  Soient $\mathfrak{X} = (\Gamma,c)$ une configuration cohérente et
  $\Gamma'\subset \Gamma$ une classe de couleurs en relation au coloriage
  induit par $c$ sur $\Gamma$. Alors la sous-structure
  induite $\mathfrak{X}\lbrack \Gamma'\rbrack$ est une configuration cohérente.
\end{exer}
Ici, c'est un cas spécial de (\ref{it:klarb}) qu'il faut utiliser: la couleur 
$c(x_1,\dotsc,x_n)$ {\og connaît\fg} les couleurs $c(x_1),\dotsc,c(x_n)$,
puisque $c(x_i)=c(x_i,\dotsc,x_i)$.

Soient $0\leq l<k$ et $\vec{x} \in \Gamma^l$. Nous colorions
$\Gamma^{k-l}$ comme suit: pour $\vec{y}\in \Gamma^{k-l}$,
\[c_{\vec{x}}(\vec{y}) = c(\vec{x} \vec{y}).\]
  En résulte une structure de partition $(k-l)$-aire
  $\mathfrak{X}_{\vec{x}} = \left(\Gamma,c_{\vec{x}}\right)$.
  \begin{exer}
    Soit $\mathfrak{X} = (\Gamma,c)$ une structure de partition;
    soit $\vec{x}\in \Gamma^l$, $0\leq l <k$. Alors
    \begin{enumerate}
      \item
    $c_{\vec{x}}$ est un raffinement du coloriage du squelette $\mathfrak{X}^{(k-l)}$. \item Si $\mathfrak{X}$ est cohérente, $\mathfrak{X}_{\vec{x}}$
        l'est aussi.
    \end{enumerate}
  \end{exer}
  Il est clair que, de plus, $\mathfrak{X}_{\vec{x}}$ est
  {\em canonique en relation à $\vec{x}$}, ce qui veut dire que
  $\mathfrak{X}\to \mathfrak{X}_{\vec{x}}$
  commute avec l'action sur $\Gamma$ du stabilisateur dans $\Sym(\Gamma)$
  des points $x_1,\dotsc,x_l$.

  \begin{defi}\label{def:uniprim}
    Une configuration cohérente $(\Gamma,c)$
    est dite {\em homogène} si la couleur
    $c(x,x,\dotsc,x)$ de tout sommet $x\in \Gamma$ est la même.
    Une configuration
    cohérente classique est dite {\em primitive} si elle est homogène
    et les graphes $\mathscr{G}_r = \{(x,y): x,y\in \Gamma, c(x,y)=r\}$
    (pour toute couleur $r$ telle que $c(x,y)= r$ pour au moins une paire
    $(x,y)$ avec $x\ne y$) sont tous connexes. Elle est
    dite {\em uniprimitive} si elle est primitive et non triviale.
  \end{defi}
  Nous n'avons pas besoin de préciser si ces graphes son connexes dans
  le sens propre (à savoir, il y a un chemin de tout sommet à tout autre,
  respectant l'orientation) ou dans le sens faible (sans compter l'orientation):
  le fait que $(\Gamma,c)$ soit cohérente, classique et homogène implique que
  $d^+_r(x)=|\{y\in \Gamma: (x,y)\in \mathscr{G}_r\}|$
  est indépendant de $x$ (pourquoi?),
  ce qui implique que toute composante faiblement connexe de
  $\mathscr{G}_r$ est connexe (exercice).

  \begin{exer}\label{ex:noclique}
    Soit $\mathfrak{X}=(\Gamma,c)$ une configuration cohérente classique uniprimitive.
    Il n'y a aucun ensemble $B\subset \Gamma$, $|B|>|\Gamma|/2$, tel que
    la restriction de $\mathfrak{X}$ à $B$ soit une clique.
  \end{exer}
  \noindent{\sc Solution} --
  Si les arêtes de la grande clique sont sensées être blanches,
  soit {\em noir} une autre couleur d'arêtes de
  $\mathfrak{X}$, et soit $\mathscr{G} = \mathscr{G}_{\text{noir}}$.
  Or, pour
  un graphe orienté birégulier\footnote{Voir la définition du \S \ref{subs:hypdes}.} $\mathscr{G}$ non vide avec $\Gamma$ comme ensemble de sommets, il est impossible qu'il y ait un ensemble $B\subset \Gamma$,
  $|B|>|G|/2$, tel que  la réduction du graphe à $B$ soit vide (pourquoi?).
  \qed 
  
  \begin{exer}\label{ex:samesize}
    Soit $(\Gamma,c)$ une configuration cohérente classique homogène.
      \begin{enumerate}
      \item\label{it:bibr1} Soit $r_0,\dotsc,r_k$ une séquence de couleurs.
        Alors, si $x_0,x_k\in \Gamma$ sont tels que $c\left(x_0,x_k\right)=r_0$,
        le nombre de $x_1,\dotsc,x_{k-1}\in \Gamma$ tels que
        $c\left(x_{i-1},x_i\right) = r_i$ pour tout $1\leq i\leq k$ dépend
        seulement de $r_0,\dotsc,r_k$.
      \item\label{it:bibr2}
        Pour toute couleur $r$, toute composante connexe de
    $\mathscr{G}_r$ est de la même taille.
    \end{enumerate}
  \end{exer}
  \noindent{\sc Solution (esquisse)} --- En (\ref{it:bibr1}), le cas
  $k=2$ vaut par la définition de {\og cohérent\fg}; prouvez les cas $k>2$
  par induction. Pour prouver (\ref{it:bibr2}), utilisez (\ref{it:bibr1}).

  \subsection{Le raffinement canonique $k$-aire à la façon de Weisfeiler-Leman}\label{subs:wl}

  Définissons un foncteur $F_3$ qui envoie une configuration
  $\mathfrak{X}=(\Gamma,c)$ à une configuration cohérente
  $F_3(\mathfrak{X})= (\Gamma,c')$. Comme $F_3(\mathfrak{X})$ sera
  un raffinement de $\mathfrak{X}$, nous aurons
$\Iso(\mathfrak{X},\mathfrak{Y}) = \Iso(F_3(\mathfrak{X}),F_3(\mathfrak{Y}))$.

  L'algorithme \ref{alg:weislem}, qui calcule $F_3$, est basé sur une idée
  de Weisfeiler et Leman\footnote{Aussi appelé Lehman, mais \cite{Ba} indique que le deuxième auteur préférait {\em Leman}. Deux transformations naturelles $L\to \El$, $\El\to L$ peuvent ne pas être l'inverse l'une de l'autre.}
  \cite{WL}.
  Il s'agit d'itérer une procédure de raffinement. Si, dans une itération,
  aucun raffinement ne se produit -- c'est-à-dire, si les classes d'équivalence
  du nouveau coloriage $\mathscr{C}_i$ sont les mêmes que celles de
  l'ancien coloriage $\mathscr{C}_{i-1}$ -- alors, (a) aucun raffinement
  ne se produira dans le futur, (b) le coloriage $\mathscr{C}_{i-1}$ est
  déjà cohérent.

  Si le coloriage $\mathscr{C}=\mathscr{C}_0$ 
  a $r$ couleurs différentes du début,
  il est clair qu'il ne peut être raffiné que
  $|\Gamma|^k - r$ fois. Alors, $|\Gamma|^k-r$ itérations sont suffisantes
pour produire une configuration cohérente. 
En particulier, si l'indexation est faite en temps logarithmique, et le
vecteur dans le pas \ref{sta:thisline} de l'algorithme \ref{alg:weislem}
est représenté comme un vecteur creux (puisque son nombre d'entrées non-nulles
est au plus $|\Gamma|$), le temps pris par l'algorithme est
$O(k^2 |\Gamma|^{2k+1} \log |\Gamma|)$. (En outre,
\cite[\S 2.8.3]{Ba} affirme une borne plus forte.)

  Les algorithmes de type Weisfeiler-Leman étaient autrefois regardés comme
  une approche plausible au problème de l'isomorphisme de graphes. Depuis
  \cite{CFI}, \cite{EvP}, il est clair qu'il ne se suffisent pas
  à eux-mêmes. Ils sont quand même un outil précieux. La version
  $k$-aire ici est due à Babai-Mathon \cite{Ba3} et Immerman-Lander
  \cite{ImL}. 
  
\begin{algorithm}
  \caption{Weisfeiler-Leman pour les configurations $k$-aires.}\label{alg:weislem}
  \begin{algorithmic}[1]
    \Function{WeisfeilerLeman}{$\Gamma$,\; $k$,\; $c:\Gamma^k\to \mathscr{C}$}
    \State $\mathscr{C}_0\gets \mathscr{C}$;
    $c_0\gets c$; $i_0\gets |\Gamma^k|-|c(\Gamma^k)|$
    \For{$i=1$   \textbf{jusqu'à} $i_0$}
    \State $\mathscr{I}_i\gets \emptyset$
    \For $\vec{x} \in \Gamma^k$
    \State\label{sta:thisline} $\nu\gets
    \left(c_{i-1}(\vec{x}),
    \left(\left|\left\{z\in \Gamma: c_{i-1}(\vec{x}^j(z)) = r_j\;\;\;
    \forall\; 1\leq j\leq k\right\}\right|\right)_{\vec{r}\in \mathscr{C}_{i-1}^k}
    \right)$
    \State $c_i(\vec{x}) = \text{\textsc{Indexeur}}\left(\mathscr{I}_i,
    \nu\right)$\Comment{\textsc{Indexeur} est
      comme dans l'algorithme \ref{alg:struraff}}
    \EndFor
    \State $\mathscr{C}_i\gets \text{indices de $\mathscr{I}_{i}$}$    
    \EndFor
    \State \Return $\left(c_{n-i_0}:\Gamma^k\to \mathscr{C}_{n-i_0},
    (\mathscr{I}_i)_{1\leq i
      \leq n-i_0}\right)$ \Comment{$(\mathscr{I}_i)$ donne du sens à $\mathscr{C}_{n-r}$}
    \EndFunction
  
  \end{algorithmic}
\end{algorithm}

\subsection{Graphes, hypergraphes et {\em designs} en blocs}\label{subs:hypdes}
Nous savons déjà qu'un {\em graphe} est une paire $(V,A)$, où $V$
est un ensemble ({\og sommets\fg}) et $A$ est une collection de
paires d'éléments de $V$ (voire de sous-ensembles de $V$ avec deux éléments,
si le graphe est non orienté). Un graphe non orienté est dit {\em régulier}
si le degré de tout sommet est le même; un graphe
orienté est dit {\em birégulier} si le {\em degré sortant}
$d^+(v)=|\{w\in V: (v,w)\in A\}|$ et le {\em degré entrant}
$d^-(v)=|\{w\in V: (v,w)\in A\}|$ sont
indépendants de $v$. (Pour $V$ fini, ils sont forcément la même constante.)

Un {\em graphe biparti} est un triplet
$(V_1,V_2;A)$ avec $A\subset V_1\times V_2$. Un graphe biparti
est {\em semirégulier} si le degré\footnote{
Nous omettons les mots
{\og entrant\fg} et {\og sortant\fg}, puisqu'il est évident qu'il s'agit du degré
entrant dans le cas de $v_1$ et du degré sortant dans le cas de $v_2$.}
$d^+(v_1)$ est indépendant
de $v_1\in V_1$, et le degré $d^-(v_2)$ est indépendant de $v_2\in V_2$. 

\begin{exer}\label{ex:biparhom}
      Soit $\mathfrak{X}=(\Gamma,c)$ une configuration cohérente classique homogène.
  \begin{enumerate}
  \item\label{it:bar1}
    Soient $C_1$, $C_2$ deux classes de couleur, et soit {\em vert} une
couleur d'arêtes en $C_1\times C_2$. Alors,
  le graphe biparti $(C_1,C_2;\mathscr{G}_{\text{vert}})$ est semirégulier.
\item\label{it:bara2}
  Soit $y\in \Gamma$, et $L_i(y) = \{x\in \Gamma: c(x,y) = i\}$. Soient
  {\em lin}, {\em bis} et {\em terre} trois couleurs d'arêtes. Alors,
  pour $L_1 = L_{\text{lin}}(y)$ et $L_2 = L_{\text{bis}}(y)$,
  le graphe biparti $(L_1,L_2; R_{\text{terre}} \cap (L_1\times L_2))$
  est semirégulier.
\end{enumerate}
\end{exer}

\begin{exer}\label{ex:biparcomp}
      Soit $\mathfrak{X}$ une configuration cohérente classique homogène.
      Soient $C_1$, $C_2$ deux classes de couleur. Soient {\em vert} une couleur
      d'arêtes en $C_1\times C_2$ et 
    {\em rouge} une couleur d'arêtes en $C_2\times C_2$. Soient
    $B_1,\dotsc,B_m$ les composantes connexes de $\mathscr{G}_{\text{rouge}}$
    en $C_2$. Définissons le graphe biparti $Y = (C_1,\{1,\dotsc,m\}; D)$
    comme suit: $(x,y)\in D$ ssi $(x,y)\in \mathscr{G}_{\text{vert}}$ pour au moins un
    $y\in B_i$. Alors $Y$ est semirégulier.
\end{exer}
\noindent{\sc Solution} ---
Notez que, pour $y\in B_i$ et $x\in C_1$, $(x,y')$ est vert pour au moins
un $y'\in B_i$ ssi il existe $x_0=x, x_1,\dotsc, x_m$ tels que $(x_i,x_{i+1})$ est rouge pour $0\leq i<m$ et $(x_m,y)$ est vert. Concluez par l'exercice
\ref{ex:samesize}\ref{it:bibr1} que tous les sommets en $\{1,\dotsc,m\}$ ont
le même degré en $Y$.

De façon analogue, montrez que, pour $x\in V_1$ et $y\in B_i$
tels que $(x,y)$ est rouge, le nombre de $z\in B_i$ tels que $(x,z)$
est rouge ne dépend pas de $x$, $y$ ou $i$. Notons ce nombre~$q$.
Alors, le degré de tout $v\in C_1$ est son degré en $X$, divisé par $q$.
Par (\ref{it:bar1}), il ne dépend donc pas de $v$. \qed

Un graphe biparti est {\em complet} (en tant que graphe biparti) si
$A = V_1\times V_2$. Un graphe biparti qui n'est ni vide ni complet est
appelé {\em non trivial}.

Un {\em hypergraphe} $\mathscr{H} = (V,\mathscr{A})$ consiste en un ensemble
$V$ ({\og sommets\fg}) et une collection $\mathscr{A}$ de sous-ensembles de $V$
({\og arêtes\fg}), peut-être avec 
des sous-ensembles répétés. Un hypergraphe est dit {\em $u$-uniforme} si
$|A|=u$ pour tout $A\in \mathscr{A}$. Il est dit {\em régulier de degré~$r$}
si tout $v\in V$ appartient à exactement $r$ ensembles $A$ dans $\mathscr{A}$.

L'hypergraphe $u$-uniforme complet sur $V$ est $(V,\{A\subset V: |A|=u\})$,
où chaque ensemble $A$ est compté une fois. Un
{\em coloriage des arêtes} de l'hypergraphe complet est une application
de $\{A\subset V: |A|=u\}$ à un ensemble fini $\mathscr{C}$.

Un {\em block design équilibré} (BDE) de paramètres $(v,u,\lambda)$
est un hypergraphe avec $|V|=v$ sommets, 
$u$-uniforme et régulier de degré $r\geq 1$, tel que toute paire 
$\{v_1,v_2\}$ de sommets distincts est contenue dans
exactement $\lambda\geq 1$ arêtes ({\og blocks\fg}).
Un block design {\em dégénéré} a la même définition, mais avec $\lambda=0$,
et la condition additionnelle d'être un hypergraphe régulier. (La régularité
peut être déduite de la définition si $\lambda\geq 1$.) Un block design
est {\em incomplet} si $u<v$.
Notons $b$ le nombre $|\mathscr{A}|$ d'arêtes d'un BDE.
\begin{prop}[Inégalité de Fisher\footnote{Si, R. A. Fisher, le
      statisticien.
Ici {\em design} vient d'{\em experimental design}.} \cite{F}]\label{prop:fishy}
  Pour tout block design équilibré incomplet, $b\geq v$.
\end{prop}
Il est aisé de voir que cette inégalité est vraie même pour les designs
dégénérés.

Les blocks designs admettent une généralisation. Un
{\em design $t$-$(v,u,\lambda)$} est un hypergraphe $(V,\mathscr{A})$
$u$-uniforme avec $v=|V|$ sommets tel que tout $T\subset V$ de taille $t$
est contenu dans exactement $\lambda$ arêtes. Ici $t\geq 2$ et
$\lambda\geq 1$.
Nous écrivons toujours
$b=|\mathscr{A}|$.
\begin{prop}[\cite{RChW}]\label{prop:RW}
  Pour tout design $t$-$(v,u,\lambda)$ et tout $s\leq \min(t/2,v-u)$, nous
  avons $b\geq \binom{v}{s}$.
\end{prop}

\subsection{Schémas de Johnson}\label{sec:schjoh}

Un {\em schéma d'association} est une configuration cohérente
classique $(\Gamma,c:\Gamma^2\to \mathscr{C})$ telle que
$c(x,y)=c(y,x)\;\; \forall x,y\in \Gamma$.
(Il s'agit donc d'un sens du mot {\em schéma} qui n'a rien à voir avec les {\em schémas} de la géométrie algébrique.)

Soient $s\geq 2$ et $r\geq 2 s + 1$. Un {\em schéma de Johnson}
$\mathscr{J}(r,s) = (\Gamma,c)$ est donné par
\[
  \Gamma = \mathscr{S}_s(\Lambda) =
  \{S\subset \Lambda: |S|=s\},\;\;\;\;\;\;\;\;\;\;\;\;
  c(S_1,S_2) = |S_1\setminus (S_1\cap S_2)|,\]
  où $\Lambda$ est un ensemble à $r$ éléments. La relation $R_i$
  est bien sûr l'ensemble \[R_i = \{(S_1,S_2): c(S_1,S_2) =i\}.\]
  
  Notons que nous avons défini implicitement un foncteur de la catégorie
  d'ensembles $\Lambda$ avec $|\Lambda|=r$ à la catégorie de schémas de
  Johnson. Ceci est un foncteur plein; autrement dit, les seuls
  automorphismes de $\mathscr{J}(r,s)$ sont ceux qui sont induits par
  $\Sym(\Lambda)$.

\subsection{Identification de groupes et de schémas}\label{subs:idgroup}

Il est une chose de démontrer que deux groupes $G$, $H$ sont isomorphes,
et une autre de construire un isomorphisme $\phi$ de façon explicite
entre eux. Cette dernière tâche
implique, au moins, de donner les images $\phi(g_1),\dotsc,\phi(g_r)$
de générateurs $g_1,\dotsc,g_r$ de $G$.

Voyons un cas particulier qui nous sera crucial. Nous aurons un groupe
de permutation $G<\Sym(\Gamma)$, et nous saurons qu'il est isomorphe au
groupe abstrait $\Alt_m$. Comment construire un isomorphisme?

Si $m$ n'est pas trop petit en relation à $n = |\Gamma|$, il est connu
que $G$ doit être isomorphe à un groupe de permutation de la forme
$\Alt_m^{(k)}$, qui n'est autre que le groupe $\Alt_m$ agissant
sur l'ensemble $\mathscr{S}_k(\Lambda_0) =
\{S\subset \Lambda_0: |S|=k\}$ à $\binom{m}{k}$
éléments, où $\Lambda_0$ est un ensemble à $m$ éléments.\footnote{
  Babai nomme les groupes $\Alt_m^{(k)}$ {\em groupes de Johnson}, par
analogie avec les schémas de Johnson. Puisque $\Alt_m^{(k)}$ n'est
  qu'un déguisement de $\Alt_m$, ne faudrait-il pas appeler ce dernier
  {\em groupe de Ramerrez}?} En d'autres termes, 
il existe une bijection $\iota_0:\Gamma\to \mathscr{S}_k(\Lambda_0)$
et un isomorphisme $\phi_0:G\to \Alt(\Lambda_0)$ tels que
\[\iota_0\left(\omega^g\right) = \iota_0(\omega)^{\phi_0(g)}.\]
Le problème consiste à construire $\iota:\Gamma\to \mathscr{S}_k(\Lambda)$
et $\phi:G\to \Alt(\Lambda)$,
calculables en temps polynomial, avec ces mêmes propriétés.

Nous suivons \cite{BLS}. Soient $\Upsilon\subset \Gamma\times \Gamma$
l'orbitale la plus petite de $G$ (hors la diagonale
$(\{\omega,\omega\}: \omega\in \Gamma\}$); soit $\Delta\subset
\Gamma\times \Gamma$ l'orbitale la plus grande. Nous supposerons
que $m>(k+1)^2-2$, ce qui revient à dire que $n$ n'est pas trop grand en
relation à $m$.\footnote{Si $m\leq (k+1)^2-2$, alors $n$ est si grand que
  $m! = n^{O(\log n)}$. En ce cas, nous pouvons enlever le groupe $G$
  (c'est-à-dire, dans l'application qui nous intéressera,
  un quotient $G/N$) de façon brutale, comme dans 
  le cas 2 de la preuve du théorème \ref{thm:luxor} (Luks). Nous pourrions
  aussi nous passer de la supposition $m>(k+1)^2-2$ au coût de quelques
  complications en ce qui suit. En particulier, $\phi(\Delta)$ ne serait pas
  $R_k$ comme dans (\ref{eq:ludod}), sinon un autre $R_j$.}
  Alors,
\begin{equation}\label{eq:ludod}\begin{aligned}
\phi(\Upsilon) &= R_1 = \{(S_1,S_2)\in \mathscr{S}_k(\Lambda_0): |S_1\cap S_2|=k-1\},\\
\phi(\Delta) &= R_k = \{(S_1,S_2)\in \mathscr{S}_k(\Lambda_0): S_1\cap S_2=\emptyset\}.
\end{aligned}\end{equation}

Définissons, pour $(x,y)\in \Upsilon$,
\[
B(x,y) = \{z\in \Gamma: (x,z)\notin \Delta, (y,z)\in \Delta\}.\]
Ceci est l'ensemble de tous les $z$ tels que $\iota_0(z)$ intersecte
$\iota_0(x)$ mais pas $\iota_0(y)$. Soit
\[C(x,y) = \Gamma \setminus  \bigcup_{z\in B(x,y)} \{r: (z,r) \in \Delta(z)\}.\]
Alors \[\begin{aligned}
\iota_0&(C(x,y))\\ &= \{S\in \mathscr{S}_k(\Lambda_0): 
S\cap S' \ne \emptyset\;\; \forall S'\in \mathscr{S}_k(\Lambda_0)
\; \text{t.q. $S'\cap \iota_0(x)\ne \emptyset$, $S'\cap \iota_0(y) = \emptyset$}
\}\\ &= \{S\in \mathscr{S}_k(\Lambda_0): i\in S\},\end{aligned}\]
où $i$ est l'élément de $\iota_0(x)$ qui n'est pas dans $\iota_0(y)$.

Soit $\Lambda$ la collection $\{C(x,y): (x,y)\in \Upsilon\}$, sans
multiplicités. Nous pouvons calculer et comparer $C(x,y)$ pour $(x,y)$
donné, et calculer et indexer $\Lambda$, tout
en temps polynomial. 
Nous calculons, aussi en temps polynomial, l'action de $G$ sur
$\Lambda$ induite par l'action de $G$ sur $\Upsilon$. Ceci définit
$\phi:G\to \Alt(\Lambda)$.

Il y a une
bijection naturelle $j:\Lambda\to \Lambda_0$ qui commute avec l'action de $G$:
elle envoie $C(x,y)$ à $i$, où $i$ est
l'élément de $\Lambda_0$ tel que $\iota_0(C(x,y)) =
\{S\in \mathscr{S}_k(\Lambda_0): i\in S\}$.
Il est clair que, pour $\omega\in \Gamma$, $\omega\in C(x,y)$ ssi
$j(C(x,y))\in \iota_0(\omega)$.
Ainsi, nous obtenons la bijection
$\iota:\Gamma\to \mathscr{S}_k(\Lambda)$, donnée par
\[\iota(\omega) = \{\gamma\in \Lambda: \omega \in \gamma\}.\]
Celle-ci satisfait $\iota\left(\omega^g\right) = \iota(\omega)^{\phi(g)}$.

Les applications $\phi$, $\iota$ sont donc celles que nous désirions; nous
avons construit un isomorphisme explicite entre $G$ et $\Alt(\Lambda)$.
Notons que cette même procédure nous permet de construire un isomorphisme
explicite entre, d'un côté,
un schéma d'association (\S \ref{sec:schjoh}) qu'on
sait être isomorphe à un schéma de Johnson $\mathscr{J}(m,k)$, et, de
l'autre côté,
ce même schéma.

\newpage
\section{La procédure principale}
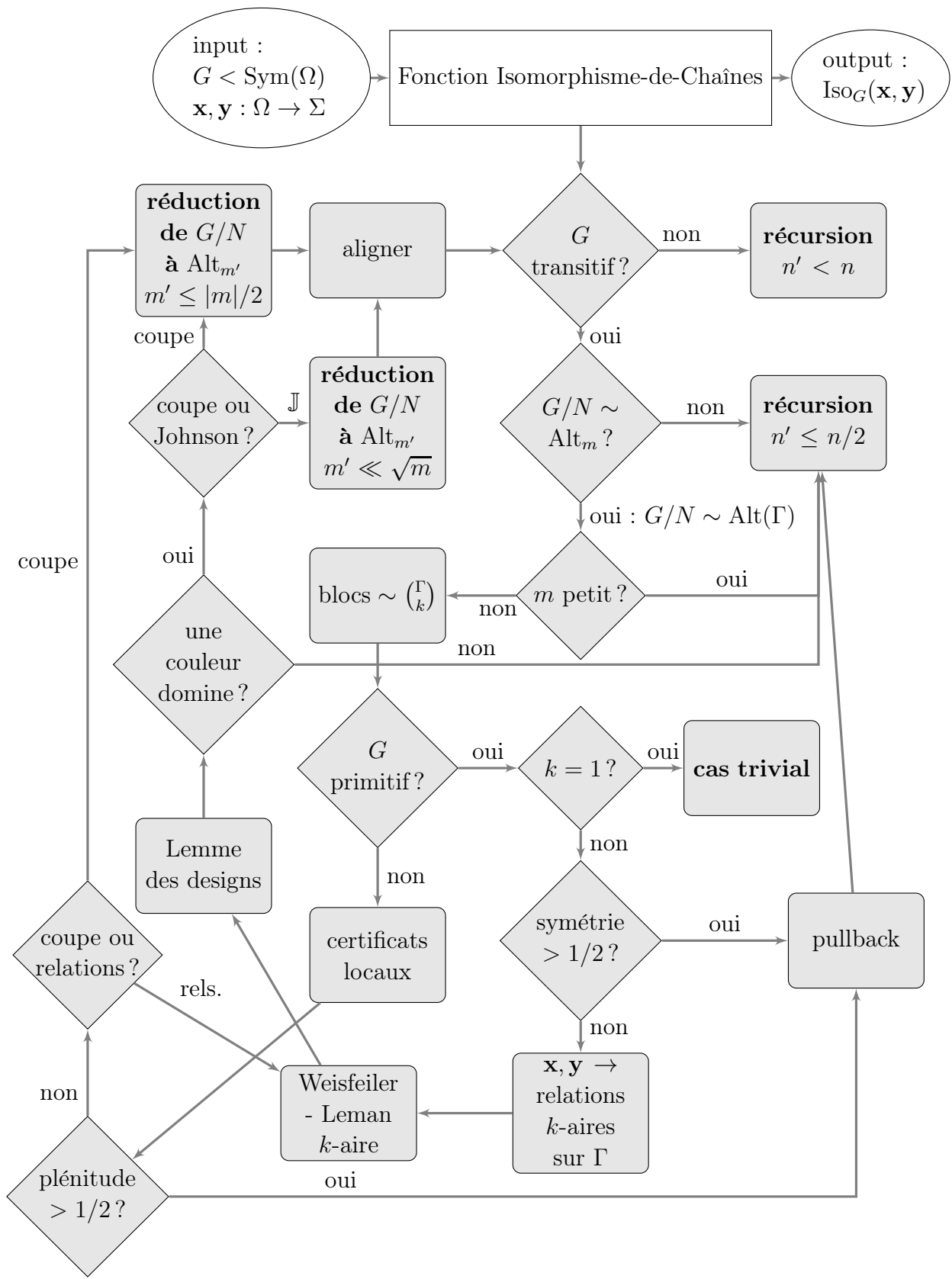
\label{sec:procprin}
\begin{center}
\begin{tikzpicture}[scale=2, node distance = 2cm, auto]
  \node [rect] (init) {Fonction Isomorphisme-de-Chaînes};
    \node [cloudin, left of=init, node distance = 5.5cm] (input) {
    input:\\ $G<\Sym(\Omega)$\\ $\mathbf{x},\mathbf{y}:\Omega\to \Sigma$};
    \node [cloudout, right of=init, node distance = 5cm] (output) {output: $\Iso_G(\mathbf{x},\mathbf{y})$};
    \node [decision, below of=init] (trans) {$G$ transitif?};
    \node [decision, below of=trans] (idtop) {$G/N \sim \Alt_m$?};
    \node [block, right of=idtop, node distance = 4.1cm] (luks) {\bf récursion
      $n'\leq n/2$};
        \node [block, right of=trans, node distance = 4.1cm] (petrec) {\bf récursion
      $n'<n$};
    \node [block, left of = trans, node distance = 3.5cm] (align) {aligner};
    \node [decision, below of=idtop, node distance = 3cm] (gampetit) {$m$ petit?};
    \node [block, left of=gampetit, node distance = 3.5cm] (struc)
          {$\text{blocs} \sim \binom{\Gamma}{k}$};
    \node [decision, below of=struc] (primi)
          {$G$ primitif?};
          \node [decision, right of = primi, node distance = 3.5cm] (k1) {$k=1$?};
          \node [block, right of =k1, node distance = 2.95cm] (triv) {\bf cas trivial};
          \node [decision, below of = k1] (sym) {symétrie $>1/2$?};
                    \node [block, right of =sym, node distance = 4.75cm] (pull) {pullback};
                    \node [block, below of = sym, node distance = 3cm] (xrel) {$\mathbf{x},
                    \mathbf{y}\to$ relations $k$-aires sur $\Gamma$};
                    \node [block, left of =xrel, node distance = 4cm] (wl) {Weisfeiler - Leman $k$-aire};
                    \node [block, below of = primi, node distance = 3.25cm] (cert) {certificats locaux};
  \node [decision, left of = cert, node distance = 5cm] (splirel) {coupe ou relations?};
                    \node [decision , below of = splirel, node distance = 4.2cm] (certsym) {plénitude $>1/2$?};
                    \node [block, left of = idtop, node distance = 3.5cm] (redsqrt) {\bf réduction de $G/N$ à $\Alt_{m'}$\\ $m'\ll \sqrt{m}$};
                    \node [decision, left of = redsqrt, node distance = 3cm] (spj) {coupe ou {\textsc Johnson}?};
                    \node [block, above of = spj, node distance = 3cm] (redhalf)
                          {\bf réduction de $G/N$ à $\Alt_{m'}$\\ $m'\leq |m|/2$};
                          \node [decision, below of = spj, node distance = 4.2cm] (coldom) {une couleur domine?};
                          \node [block, below of = coldom, node distance = 3.5cm] (design) {Lemme des designs};

    \path[line] (input) -- (init);
    \path[line] (init) -- (output);      
    \path [line] (init) -- (trans);
        \path [line] (trans) -- node [near start, color=black] {non} (petrec);
    \path [line] (trans) -- node [,color=black] {oui} (idtop);
    \path [line] (idtop) -- node [, color=black] {non} (luks);
    \path [line] (idtop) -- node [, color=black] {oui:
      $G/N \sim \Alt(\Gamma)$} (gampetit);
    \path [line] (gampetit) -| node [near start, color=black] {oui} (luks);
    \path [line] (gampetit) -- node [near start, color=black] {non}(struc);
    \path [line] (struc) -- (primi);
    \path [line] (align) -- (trans);
    \path [line] (primi) -- node [, color=black] {oui} (k1);
    \path [line] (k1) -- node [, color=black] {oui} (triv);
    \path [line] (k1) -- node [, color = black] {non} (sym);
    \path [line] (sym) -- node [, color = black] {oui} (pull);
    \path [line] (pull) -- (luks);
    \path [line] (sym) -- node [near start, color = black] {non} (xrel);
    \path [line] (xrel) -- (wl);
    \path [line] (primi) -- node [, color = black] {non} (cert);
    \path [line] (cert) -- (certsym);
    \path [line] (certsym) -- node [near start, color=black] {non} (splirel);
    \path [line] (splirel) -- node [near start, color = black] {
      rels.} (wl);
    \path [line] (certsym) -| node [very near start, color = black] {
      oui} (pull);
    \path [line] (wl) -- (design);
    \path [line] (design) -- (coldom);
    \path [line] (coldom) -| node [near start, color=black] {non\;\;\;\;\;\;\;\;\;\;\;\;\;\;\;\;\;\;\;\;\;\;\;\;} (luks);
    \path [line] (coldom) -- node [near start, color=black] {oui} (spj);
    \path [line] (redhalf) -- (align);
    \path [line] (redsqrt) -- (align);
    \path [line] (spj) -- node [near start, color=black] {coupe} (redhalf);
    \path [line] (spj) -- node [, color= black] {$\mathbb J$} (redsqrt);
    \path [line] (splirel) |- node [near start, color=black]  {coupe} (redhalf);
\end{tikzpicture}
\end{center}
\subsection{Premiers pas: récursion à la façon de Luks}

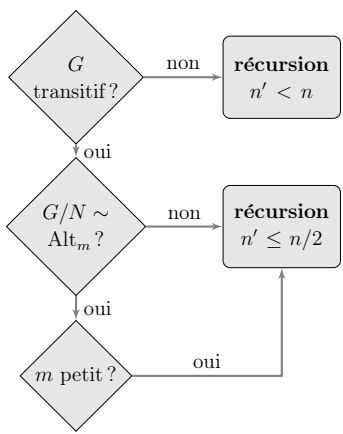
\begin{center}
  \scalebox{0.66}{
\begin{tikzpicture}[scale=2, node distance = 2cm, auto]
    \node [decision, below of=init] (trans) {$G$ transitif?};
    \node [decision, below of=trans, node distance = 3cm] (idtop) {$G/N \sim \Alt_m$?};
           \node [block, right of=trans, node distance = 4.1cm] (petrec) {\bf récursion
      $n'<n$};
    \node [block, right of=idtop, node distance = 4.1cm] (luks) {\bf récursion
      $n'\leq n/2$};
    \node [decision, below of=idtop, node distance = 3cm] (gampetit) {$m$ petit?};
    \path [line] (trans) -- node [, color=black] {non} (petrec);
    \path [line] (trans) -- node [,color=black] {oui} (idtop);
        \path [line] (idtop) -- node [,color=black] {oui} (gampetit);
    \path [line] (idtop) -- node [, color=black] {non} (luks);
    \path [line] (gampetit) -| node [near start, color=black] {oui} (luks);
\end{tikzpicture}
}
\end{center}

Les premiers pas de la procédure sont ceux de la preuve du Théorème
\ref{thm:luxor} (Luks). 
En particulier, si $G<\Sym(\Omega)$ n'est pas transitif, nous procédons exactement
comme dans le cas non transitif de la preuve du Théorème \ref{thm:luxor}.
Bien qu'il soit possible que $n=|\Omega|$ ne décroisse que très légèrement,
la récursion marche, puisque son coût est aussi très léger dans ce cas: nous
n'avons qu'à subdiviser le problème selon les orbites de $G$.

Supposons que $G$ soit transitif. Nous savons que nous pouvons trouver
rapidement un système de blocs minimal $R = \{B_i: 1\leq i\leq r\}$,
$B_i\subset \Omega$ (\S \ref{subs:orbl}). Par Schreier-Sims, nous trouvons
aussi, en temps polynomial, le sous-groupe $N\triangleleft G$ des éléments
$g\in G$ tels que $B_i^g = B_i$ pour tout $i$. Le groupe $H = G/N$ agit sur $R$.

Au lieu du Théorème \ref{thm:bcp}
\cite{BCP}, nous utiliserons une conséquence de la Classification
des Groupes Finis Simples (CGFS). Elle a été dérivée pour la première fois par
Cameron, puis raffinée par Mar\'oti. 

\begin{theo}[\cite{Cam}, \cite{Ma}]\label{thm:cam}
  Soit $H< \Sym(R)$ un groupe primitif, où $|R| = r$ est plus
  grand qu'une constante absolue. Alors, soit\footnote{Pour nous, $\log_2$
    désigne le logarithme en base $2$, et non pas le logarithme itéré
    $\log \log$.} 
  \begin{enumerate}
   \item\label{it:rila} $|H| < r^{1 + \log_2 r}$, soit
   \item\label{it:rilb}
    il y a un $M\triangleleft H$
    tel que $R$ se subdivise\footnote{
      L'énoncé dans \cite{Cam}, \cite{Ma} est plus fort: il décrit toute l'action de $H$ sur $R$. \`A vrai dire, le groupe $M$ est isomorphe, en tant que groupe de permutation, à $(\Alt_m^{(k)})^s$, $s\geq 1$. Nous avons $r =
\binom{m}{k}^s$.}
    en
    un système de $\binom{m}{k}$ blocs sur lequel $M$ agit comme un groupe 
    $\Alt_m^{(k)}$, $m\geq 5$. En plus,
    $\lbrack H: M\rbrack \leq r$.
  \end{enumerate}
\end{theo}
La borne $\lbrack H:M\rbrack\leq r$ se déduit de $m>2$, $|H|\geq r^{1+\log_2 r}$, $|H|\leq m!^s s!$, $m^s\leq r$ et
      $\lbrack H:M\rbrack\leq 2^s s!$, où $s\geq 1$ est un paramètre dans
      Cameron-Mar\'oti.

Il est possible \cite{BLS} de trouver en temps polynomial
le sous-groupe normal $M$ et les blocs de l'action de $M$.
Nous avons déjà vu au \S \ref{subs:idgroup} comment identifier explicitement
l'action de $M$ avec celle de $\Alt_m^{(k)}$.

Par ailleurs, l'algorithme de Schreier-Sims nous permet de calculer
$|H|$ en temps polynomial, et donc nous dit aussi si nous sommes dans le
cas (\ref{it:rila}). Si c'est le cas, nous procédons
comme dans le cas transitif de la preuve du Théorème \ref{thm:luxor}.
Nous réduisons ainsi le problème à $r^{1+\log_2 r}$
instances du problème pour des chaînes de longueur $\leq n/r$.

Si nous sommes dans le cas (\ref{it:rilb}) nous commençons toujours par
réduire le problème à $\lbrack H:M\rbrack$ instances du problème avec
$M$ à la place de $H$: par l'équation (\ref{eq:udu2})
et comme dans l'équation (\ref{eq:rulu}),
\[\Iso_H(\mathbf{x},\mathbf{y}) = \bigcup_{\sigma \in S} \Iso_M\left(
\mathbf{x}, \mathbf{y}^{\sigma^{-1}}\right) \sigma,\]
où $S$ est un système de représentants des classes de $M$ dans $H$.

Si $m\leq C \log n$, où $C$ est une constante,
\[|M| = \frac{m!}{2} < m^m \leq m^{C \log n}\leq (m')^{C \log n},\]
où $m' = \binom{m}{k}$.
Donc, ici comme dans le cas (\ref{it:rila}), nous nous permettons de
procéder comme dans le cas transitif de la preuve du Théorème \ref{thm:luxor}.
Nous obtenons une réduction à $\leq r\cdot (m')^{C \log n} = (m')^{O(\log n)}$
instances du problème pour des chaînes de longueur $n/m'$.
Ceci est tout à fait consistant avec l'objectif d'avoir une solution en
temps quasi-polynomial en $n$ (ou même en temps $n^{O(\log n)}$).

Il reste à savoir que faire si nous sommes dans le cas suivant:
il y a un isomorphisme $\phi:G/N \to \Alt(\Gamma)$, $|\Gamma|> C \log n$,
$C$ une constante.
(Ici nous avons déjà (i) remplacé
$G$ par la préimage de $M$ dans la réduction $G\to G/N$, et, après cela, (ii)
remplacé $N$ par le stabilisateur des blocs dans la partie (\ref{it:rilb})
du Théorème~\ref{thm:cam}.) Ce cas nous occupera pour le reste de
l'article.

\begin{center}
  * * *
\end{center}

Babai indique comment enlever la dépendance de CGFS à cette étape.
Soient $G$ et $N$ comme avant, avec $G$ transitif. Alors $G/N$ est
un groupe primitif agissant sur l'ensemble de blocs $R$.

Si
un groupe de permutations sur un ensemble $R$ est tel que son action sur
l'ensemble des paires d'éléments distincts de $R$ est transitive, le
groupe est dit {\em doublement transitif}. Or, un résultat de Pyber
qui ne dépend pas de CGFS \cite{Py} nous dit qu'un tel groupe
est soit $\Alt(R)$, soit $\Sym(R)$, soit d'ordre
$\leq |R|^{O(\log^2 |R|)}$.

Si $G/N$ est $\Alt(R)$ ou $\Sym(R)$, nous sommes dans le
cas que nous discuterons d'ici jusqu'à la fin. Si $G/N$ est doublement
transitif, mais n'est \'egal ni \`a $\Alt(R)$ ni \`a $\Sym(R)$,
nous pouvons procéder comme dans le
cas transitif de la preuve du Théorème \ref{thm:luxor}, puisque
$|G/N|\leq r^{O(\log^2 r)}$, $r = |R|\leq n$.
(Babai propose aussi un traitement alternatif,
même plus efficace et élémentaire.)

Supposons donc que $G/N$ n'est pas doublement transitif. Alors
la configuration cohérente schurienne (\S \ref{subs:concoh}) qu'elle
induit n'est pas une clique. En conséquence, nous pouvons donner
cette configuration à la procédure
\textsc{Coupe-ou-Johnson} (\S \ref{sec:coujoh}), et reprendre le fil de
l'argument à ce point-là. 

\section{La structure de l'action de $\Alt$}

\subsection{Stabilisateurs, orbites et quotients alternants}\label{subs:stab}

Nous aurons besoin de plusieurs résultats sur les épimorphismes $G\to
\Alt_k$. Ils joueront un rôle crucial dans la méthode des certificats
locaux (\S \ref{sec:certloc}). Dans la version originale \cite{Ba}, ils ont
aussi été utilisés dans
le rôle joué par \cite{BLS} dans cet exposé.

\begin{lemm}\label{lem:mustaf}
  Soit $G<\Sym(\Omega)$ primitif. Soit $\phi:G\to \Alt_k$ un épimorphisme
  avec $k>\max(8,2+\log_2 |\Omega|)$. Alors $\phi$ est un isomorphisme.
\end{lemm}
Prouver ce lemme est à peu près un exercice en théorie des groupes finis; 
il faut utiliser \cite[Prop.~1.22]{BaPS} pour le cas de socle abélien et
la conjecture de Schreier pour le cas de socle non abélien.
La conjecture de Schreier est
un théorème, mais un théorème dont la preuve dépend, à son tour, de
CGFS. 

Par contre, Pyber \cite{Py2} a donné une preuve du Lemme \ref{lem:mustaf}
qui n'utilise pas CGFS, avec une condition plus stricte:
$k> \max(C,(\log |\Omega|)^5)$, $C$ constante.
La dépendance de CGFS a donc été complètement enlevée de la preuve du
théorème principal.

\begin{defi}
  Soit $G<\Sym(\Omega)$. Soit $\phi:G\to \Sym_k$ un homomorphisme dont
  l'image contient $\Alt_k$. Alors $x\in \Omega$ est dit {\em atteint}
    si
  $\phi(G_x)$ ne contient pas $\Alt_k$.
\end{defi}

\begin{lemm}\label{lem:preaffect} Soit $G<\Sym(\Omega)$.
  Soit $\phi:G\to \Alt_k$ un épimorphisme avec
  $k>\max(8,2+\log_2 n_0)$, où
  $n_0$ est la taille de la plus grande orbite de $G$.
  \begin{enumerate}
  \item\label{it:transif}
    Si $G$ est transitif, tout $x\in \Omega$ est atteint.
  \item\label{it:atr}
    Au moins un $x\in \Omega$ est atteint.
  \end{enumerate}
\end{lemm}
\noindent{\sc Preuve (esquisse)} ---
(\ref{it:transif})  Ceci d\'ecoule immédiatement du Lemme \ref{lem:mustaf} si $G$ est primitif,
ou si $K<\ker(\phi)$ pour $K$ le stabilisateur d'un système de blocs
minimal. Il reste le cas de $\phi:K\to \Alt_k$ surjectif. En général:
\begin{quote} {\sc Lemme.}---
Pour $K_i$ arbitraires, $K<K_1\times\dots \times K_s$ et un épimorphisme
\mbox{$\phi:K\to S$,} $S$ simple, il doit y avoir un $i$ tel que $\phi$ se factorise
comme suit~: \mbox{$K\to K_i \stackrel{\psi}{\to} S$,} $\psi$ un épimorphisme.\end{quote}
En utilisant ce lemme pour les
restrictions $K_i$ de $K$ aux orbites de $K$, nous passons à une orbite $K_i$,
et procédons par induction.

(\ref{it:atr}) Soient $\Omega_1,\dotsc,\Omega_m$ les orbites de $G$, et
soit $G_i = G|_{\Omega_i}$ la restriction de $G$ à $\Omega_i$. Par le Lemme 
en (\ref{it:transif}),
il doit y avoir un $i$ tel que $\phi$ se factorise en
$G\to G_i \stackrel{\psi}{\to} \Alt_k$, $\psi$ un épimorphisme.
Alors, par (\ref{it:transif}), $(G_x)^\psi = ((G_i)_x)^\psi \ne \Alt_k$
pour tout $x\in \Omega_i$.
\qed

La proposition suivante jouera un rôle crucial au \S \ref{sec:casimp}. 
\begin{prop}\label{prop:atinl}
  Soient $G<\Sym(\Omega)$ transitif
  et $\phi:G\to \Alt_k$ un épimorphisme. Soit
  $U\subset \Omega$ l'ensemble des éléments non atteints. 
  \begin{enumerate}
  \item\label{it:unaffstab} Supposons que $k\geq \max(8,2+\log_2 n_0)$,
    où $n_0$ est la taille de la plus grande orbite de $G$.
    Alors $(G_{(U)})^\phi = \Alt_k$.
  \item\label{it:afforb} Supposons que $k\geq 5$. Si $\Delta$ est une orbite
    de $G$ qui contient des éléments atteints, alors chaque orbite de
    $\ker(\phi)$ contenue dans $\Delta$ est de longueur $\leq |\Delta|/k$.
  \end{enumerate}
\end{prop}
Rappelons que $G_{(U)} = \{g\in G: x^g = x \; \forall x\in U\}$
(stabilisateur de points).

\noindent{\sc Preuve} ---
(\ref{it:unaffstab}) Il est facile de voir que
$G$ fixe $U$ en tant qu'ensemble. Alors,
$G_{(U)}\triangleleft G$, et donc $(G_{(U)})^\phi \triangleleft G^\phi$.
Or, $G^\phi = \Alt_k$.
Supposons que $(G_{(U)})^\phi = \{e\}$. Alors $\phi$ se factorise comme suit :
\[G\to G|_U \stackrel{\psi}{\to} \Alt_k,\]
puisque $G_{(U)}$ est le noyau de $G\to G|_U$. 
Ici $\psi$ est un épimorphisme, et donc, par le Lemme
\ref{lem:preaffect} (\ref{it:afforb}), il existe un $x\in U$ tel que
$((G|_U)_x)^\psi \ne \Alt_k$. Or $((G|_U)_x)^\psi = (G_x)^\phi = \Alt_k$,
parce que $x$ est dans $U$, c'est-à-dire non atteint. Contradiction.

(\ref{it:afforb}) Comme $\Delta$ contient des éléments atteints et  est une orbite de $G$,
tout élément de $\Delta$ est atteint.
Soit $N=\ker(\phi)$, $x\in \Delta$. 
La longueur
de l'orbite $x^N$ est
\[\begin{aligned}
\left|x^N\right| &= \lbrack N:N_x\rbrack = \lbrack N:(N\cap G_x)\rbrack
= \lbrack N G_x : G_x\rbrack = \frac{\lbrack G:G_x\rbrack}{\lbrack G : N
  G_x\rbrack}\\ &= \frac{|\Delta|}{\lbrack G^\phi: (G_x)^\phi\rbrack}
= \frac{|\Delta|}{\lbrack \Alt_k : (G_x)^\phi\rbrack}.\end{aligned}\]
Or, tout sous-groupe propre de $\Alt_k$ est d'indice $\geq k$. Donc
$\left|x^N\right|\leq |\Delta|/k$.
\qed

\subsection{Le cas de grande symétrie}\label{sec:grasym}

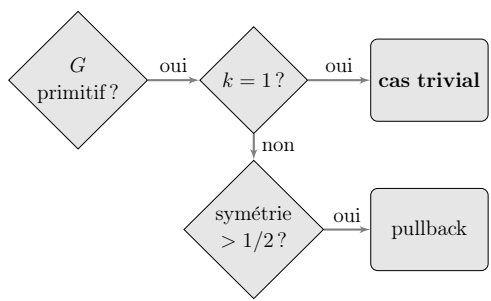
\begin{center}
    \scalebox{0.66}{
    \begin{tikzpicture}[scale=2.5, node distance = 2cm, auto]
         \node [decision] (primi)
          {$G$ primitif?};
      \node [decision, right of = primi, node distance = 3.5cm] (k1) {$k=1$?};
      \node [block, right of =k1, node distance = 3.5cm] (triv) {\bf cas trivial};
      \node [decision, below of = k1] (sym) {symétrie $>1/2$?};
      \node [block, right of =sym, node distance = 3.5cm] (pull) {pullback};
          \path [line] (primi) -- node [, color=black] {oui} (k1);
      \path [line] (k1) -- node [, color=black] {oui} (triv);
      \path [line] (k1) -- node [, color = black] {non} (sym);
      \path [line] (sym) -- node [, color = black] {oui} (pull);
    \end{tikzpicture}
  }
\end{center}

Considérons le cas de $G$ primitif. Nous pouvons supposer que
$G$ est isomorphe en tant que groupe de permutation à $\Alt_m^{(k)}$,
puisque nous avons déjà éliminé les autres cas au \S \ref{sec:procprin}
(peut-être en passant à un groupe non primitif $M$; le cas non primitif
sera traité au \S \ref{sec:casimp}).
Comme nous l'avons vu au
\S \ref{subs:idgroup}, nous pouvons construire une bijection $\iota$
entre $\Omega$ et l'ensemble $\mathscr{S}_k(\Gamma)$ des sous-ensembles avec $k$ éléments d'un ensemble $\Gamma$. Cette bijection induit
un isomorphisme $\phi:G\to \Alt(\Gamma)$.

Si $k=1$, alors $\Omega$ est en bijection avec $\Gamma$, et $G\sim \Alt_n =
\Alt_m$.
Nous sommes donc dans le cas trivial: le groupe $\Aut_G(\mathbf{x})$
consiste en les éléments de $\Alt_n$ qui permutent les lettres de
$\mathbf{x}$ de la même couleur, et $\Iso_G(\mathbf{x},\mathbf{y})$
est non vide ssi $\mathbf{x}$ et $\mathbf{y}$ ont exactement le même
nombre de lettres de chaque couleur -- où, si aucune lettre n'est répétée ni
en $\mathbf{x}$ ni en $\mathbf{y}$, nous ajoutons la condition que la
permutation de $\{1,\dotsc,n\}$ qui induit $\mathbf{x}\mapsto \mathbf{y}$
soit dans $\Alt_n$.

Alors, soit $G$ primitif, $k>1$. 

Deux éléments $\gamma_1,\gamma_2\in \Gamma$ sont des {\em jumeaux}
par rapport à un objet si la transposition $(\gamma_1 \gamma_2)$ le
laisse invariant. 
Il est clair que les jumeaux forment des classes d'équivalence, et que,
pour toute telle classe d'équivalence $C$, tout $\Sym(C)$ laisse l'objet
invariant. Notre objet sera la chaîne $\mathbf{x}$ (ou $\mathbf{y}$):
$\gamma_1$, $\gamma_2$ sont des jumeaux par rapport à $\mathbf{x}$ si, pour tout $i\in \Omega$,
$\mathbf{x}(i) = \mathbf{x}(\tau^{\phi^{-1}}(i))$, où
$\tau = (\gamma_1 \gamma_2)$.

Nous pouvons donc déterminer facilement (et en temps polynomial)
les classes d'équivalence en
$\Gamma$ (dites {\em classes de jumeaux}),
et vérifier s'il y a une classe d'équivalence $C$
de taille $>|\Gamma|/2$.
Examinons cette  possibilité puisque nous devrons l'exclure après. 

La classe $C$ de taille $>|\Gamma|/2$ est évidemment unique et donc canonique. Si
$\mathbf{x}$ a une telle classe et $\mathbf{y}$ ne l'a pas, ou
si les deux ont de telles classes, mais de tailles différentes, alors
$\mathbf{x}$ et $\mathbf{y}$ ne sont pas isomorphes.

Si $\mathbf{x}$, $\mathbf{y}$ ont des classes de jumeaux
$C_{\mathbf{x}}$, $C_{\mathbf{y}}$ de la même taille $>|\Gamma|/2$,
nous choisissons $\sigma\in \Alt(\Gamma)$ tel que
$C_{\mathbf{x}} = \left(C_{\mathbf{y}}\right)^\sigma$.
(Nous supposons $m>1$.)
En remplaçant $\mathbf{y}$
par $\mathbf{y}^{\sigma'}$, où $\sigma' = \phi^{-1}\left(\sigma^{-1}\right)$, nous réduisons notre problème au cas
$C_{\mathbf{x}} = C_{\mathbf{y}}$. (Voilà l'exemple le plus simple de
ce que Babai appelle {\em aligner}; nous avons {\em aligné}
$C_{\mathbf{x}}$ et $C_{\mathbf{y}}$.)

Alors, soit $C = C_\mathbf{x} = C_{\mathbf{y}}$. 
La partition
$\{C,\Gamma\setminus C\}$ de $\Gamma$ induit une partition
$\{\Omega_j\}_{0\leq j\leq k}$ de $\Omega$: $\omega\in \Omega_j$
  ssi $\psi(\omega)$ contient $k-j$ éléments de $C$ et $j$ éléments
  de $\Gamma\setminus C$.
  Il est ais\'e de montrer que $\alpha^{k-j} (1-\alpha)^j \binom{k}{j}
  <1/2$ pour $\alpha\in (1/2,1\rbrack$, $1\leq j\leq k$; donc,
  $|\Omega_j| < n/2$ pour $1\leq j\leq k$.

  Nous avons réduit notre problème à celui de déterminer
  $\Iso_H(\mathbf{x},\mathbf{y})$, où
  $H = \phi^{-1}\left(\Alt(\Gamma)_C\right)$. Ici le besoin de prendre
  un stabilisateur d'ensemble (à savoir,
  $\Alt(\Gamma)_C$) ne pose aucun souci: nous engendrons $H$
  en prenant des préimages $\phi^{-1}(h_1),\dotsc,\phi^{-1}(h_5)$
  de deux générateurs $h_1$, $h_2$ de $\Alt(C)<\Alt(\Gamma)$, deux générateurs
  $h_3$, $h_4$ de $\Alt(\Gamma\setminus C)<\Alt(\Gamma)$ et
  un élément $h_5\in \Alt(\Gamma)$ de la forme $(\gamma_1 \gamma_2)
  (\gamma_3 \gamma_4)$, où $\gamma_1,\gamma_2\in C$, $\gamma_3,\gamma_4\in
  \Gamma\setminus C$.
  (Si $|\Gamma|<8$, le nombre de générateurs est moindre, et la
  discussion se simplifie.) 
  Notre problème se réduit à celui de déterminer
  $\Iso_{H'}(\mathbf{x},\mathbf{y}')$ pour $\mathbf{y}' =
  \mathbf{y}$ et $\mathbf{y}' = \mathbf{y}^{h_5}$,
  où $H' = \phi^{-1}(\Alt(C)\times
  \Alt(\Gamma\setminus C)) = \phi^{-1}(\langle h_1,\dotsc,h_4\rangle)$.
  
  Comme $C$ est une classe de jumeaux pour $\mathbf{x}$, tout élément de
  $\phi^{-1}(\Alt(C))$ laisse $\mathbf{x}$ invariant.
  Si $\mathbf{x}|_{\Omega_0} \ne \mathbf{y}|_{\Omega_0}$,
  alors $\Iso_{H'}(\mathbf{x},\mathbf{y}) = \emptyset$.

  Soit alors $\mathbf{x}|_{\Omega_0} = \mathbf{y}|_{\Omega_0}$.
  Nous avons réduit notre problème \`a celui de déterminer
  $\Iso_{H'|_{\Omega'}}(\mathbf{x}|_{\Omega'},\mathbf{y}|_{\Omega'})$,
  où $\Omega' = \Omega  \setminus \Omega_0$.
  Rappelons que $H'|_{\Omega'}$ agit sur $\Omega'$ avec des orbites
  de longueur $|\Omega_i|<n/2$. Nous procédons donc comme dans le
  cas non transitif de la méthode de Luks (preuve du Thm.~\ref{thm:luxor}).

\section{Des chaînes aux schémas de Johnson}\label{subs:chasche}

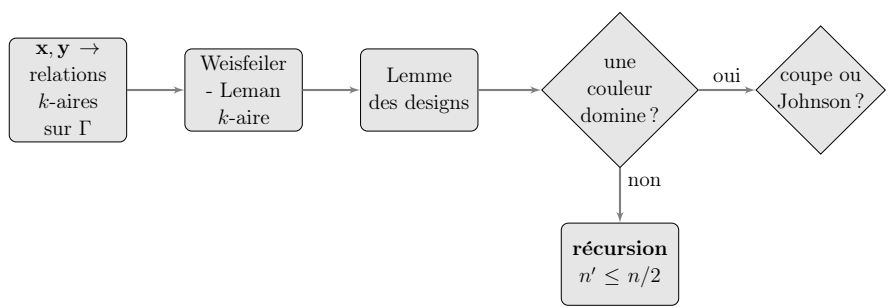
\begin{center}
      \scalebox{0.66}{
\begin{tikzpicture}[scale=2, node distance = 2cm, auto]
    \node [block, node distance = 3cm] (xrel) {$\mathbf{x},
    \mathbf{y}\to$ relations $k$-aires sur $\Gamma$};
    \node [block, right of =xrel, node distance = 3.5cm] (wl) {Weisfeiler - Leman $k$-aire};
    \node [block, right of = wl, node distance = 3.5cm] (design) {Lemme des designs};
    \node [decision, right of = design, node distance = 4cm] (coldom) {une couleur domine?};
    \node [decision, right of = coldom, node distance = 4cm] (spj) {coupe ou {\textsc Johnson}?};
    \node [block, below of = coldom, node distance = 3.5cm] (luks) {\bf récursion $n'\leq n/2$};
    \path [line] (xrel) -- (wl);
    \path [line] (wl) -- (design);
    \path [line] (design) -- (coldom);
    \path [line] (coldom) -- node [, color = black] {oui} (spj);
        \path [line] (coldom) -- node [near start, color = black] {non} (luks);
    
\end{tikzpicture}
}
\end{center}

Discutons maintenant le cas de $G$ primitif et, plus précisément, $G$
isomorphe à $\Alt_m^{(k)}$, $k\geq 2$. Maintenant nous pouvons supposer que nos
chaînes $\mathbf{x}, \mathbf{y}$ n'ont pas de classes de jumeaux de
taille $>m/2$. Les outils principaux que nous
développerons (Lemme des designs, coupe-ou-Johnson) nous seront
utiles, voire essentiels, aussi dans le cas de $G$ imprimitif.

Nous avons une bijection entre les éléments de $\Omega$
et $\{S\subset \Gamma: |S|=k\}$.
Pour $\mathbf{x}:\Omega\to \Sigma$ donné, nous avons donc une
structure relationnelle $\mathfrak{X}=(\Gamma, (R_i)_{i\in \Sigma})$ $k$-aire sur $\Gamma$:
$(x_1,\dotsc,x_k)\in R_i$ si $x_1,\dotsc,x_k$ sont tous différents et
$\mathbf{x}(\omega) = i$, où $\omega$ est l'élément de $\Omega$ qui
correspond à $\{x_1,\dotsc,x_k\}$.

Nous appliquons à $\mathfrak{X}$ le foncteur $F_1$ (\S \ref{subs:secf1}), qui fait
d'elle une structure de partition, puis le foncteur
$F_2$ (encore \S \ref{subs:secf1}), qui nous donne
une configuration $k$-aire,
et, finalement, le foncteur $F_3$ défini par Weisfeiler-Leman $k$-aire
(\S \ref{subs:wl}). Nous obtenons ainsi un raffinement
$F_3(F_2(F_1(\mathfrak{X}))) =
(\Gamma,c_\mathbf{x}:\Omega^k\to \mathscr{C})$
qui est une configuration cohérente $k$-aire.

Comme $F_1$, $F_2$, $F_3$ sont des foncteurs, l'assignation de $c_\mathbf{x}$
à $\mathbf{x}$ est canonique. Elle nous sera donc utile: si
$c_{\mathbf{x}}$ et $c_{\mathbf{y}}$ ne sont pas isomorphes sous l'action
de $\Alt_m$, alors $\mathbf{x}$ et $\mathbf{y}$ ne sont pas isomorphes
sous l'action de $\Alt_m^{(k)}$ non plus.

Nous obtiendrons une configuration cohérente classique de façon canonique
à partir de $c_\mathbf{x}$ (Lemme des {\em designs}). Soit cette nouvelle
configuration sera non triviale, soit nous obtiendrons un coloriage canonique
sans couleur dominante, ce qui nous permettra immédiatement de réduire
le problème à un certain nombre de problèmes pour des chaînes plus courtes,
comme dans l'algorithme de Luks.

Supposons, alors, que nous disposons d'une configuration cohérente classique
non triviale assignée de façon canonique à $\mathbf{x}$. La procédure
\textsc{Coupe-ou-Johnson} nous donnera l'un ou l'autre de ces deux
résultats: soit
un {\em découpage} canonique de $\Gamma$, soit un {\em schéma de Johnson}
plongé de façon canonique dans $\Gamma$.
Dans un cas comme dans l'autre, avoir une telle structure canonique
limite fortement l'ensemble d'isomorphismes et automorphismes possibles.
Nous pourrons réduire $G$ \`a un sous-groupe $\sim \Alt_{m'}$, avec
$m'\leq m/2$, dans le cas du découpage,
ou $m'\ll \sqrt{m}$, dans le cas de Johnson. Déjà $m'\leq m/2$ est
suffisante pour une récursion réussie.


\subsection{Lemme des designs}\label{sec:design}

\'Etant donnés une configuration $\mathfrak{X}=(\Gamma,c:\Gamma^k\to\mathscr{C})$
et un paramètre $1/2\leq \alpha<1$, une couleur $i$
est dite {\em $\alpha$-dominante} si $c(\gamma,\dotsc,\gamma) = i$
pour $\geq \alpha |\Gamma|$ valeurs de $\gamma\in \Gamma$. La classe
de couleurs $\{\gamma\in \Gamma: c(\gamma,\dotsc,\gamma) = i\}$
est, elle aussi, dite {\em dominante}. Par contre, si, pour toute
couleur $i$, la classe $\{\gamma\in \Gamma: c(\gamma,\dotsc,\gamma) = i\}$
est de taille $<\alpha |\Gamma|$, le coloriage est dit 
un {\em $\alpha$-coloriage}.

Comme avant, deux éléments $\gamma_1,\gamma_2\in \Gamma$ sont des {\em jumeaux}
par rapport à une structure~$\mathfrak{X}$ (ici, une configuration cohérente sur $\Gamma$)
si $(\gamma_1 \gamma_2) \in \Aut(\mathfrak{X})$.
\begin{prop}[Lemme des designs]\label{prop:design}
  Soit $\mathfrak{X}=(\Gamma, c:\Gamma^k\to \mathscr{C})$
  une configuration cohérente $k$-aire, où
  $2\leq k\leq |\Gamma|/2$. Soit $1/2\leq \alpha < 1$.
  Supposons qu'il n'y a aucune classe de jumeaux
  dans $\Gamma$ avec $>\alpha |\Gamma|$ éléments.

  Alors, au moins une des options suivantes est vraie:
  \begin{enumerate}
  \item\label{it:loro1} 
    il existe $x_1,\dotsc,x_\ell\in \Gamma$,
    $0\leq \ell<k$, tels que
    $\mathfrak{X}_{\vec{x}}^{(1)}$ n'a pas de couleur $\alpha$-dominante;
  \item\label{it:loro2} 
    il existe $x_1,\dotsc,x_\ell\in \Gamma$,
    $0\leq \ell<k-1$, tels que
    $\mathfrak{X}_{\vec{x}}^{(1)}$ a une couleur $\alpha$-dominante $C$
    et $(\mathfrak{X}_{\vec{x}})^{(2)}\lbrack C\rbrack$ n'est pas une clique.
  \end{enumerate}
  \end{prop}
La notation a \'et\'e définie dans les sections \ref{subs:secf1}
-- \ref{subs:concoh} . En particulier, le $1$-squelette
$\mathfrak{X}_{\vec{x}}^{(1)}$ est tout simplement un coloriage de $\Gamma$.

\begin{lemm}[Lemme de la grande clique]
  Soit $\mathfrak{X} = (\Gamma,c)$ une configuration cohérente classique.
  Soit $C\subset \Gamma$ une classe de couleurs avec
  $|C|\geq |\Gamma|/2$. Si $\mathfrak{X}\lbrack C\rbrack$ est une clique, alors
  $C$ est une classe de jumeaux.
\end{lemm}
\noindent{\sc Preuve} --- Supposons que $C$ n'est pas une
classe de jumeaux. Il y a donc un $x\in \Gamma$ et une couleur
(disons, {\em azur}) telle que $c(x,y)$ est de cette couleur pour au moins un
$y\in C$ mais pas pour tous. Comme $\mathfrak{X}\lbrack C\rbrack$
est une clique, $x\notin C$. Appelons la couleur de $C$ {\em carmin},
et celle de $x$ {\em bronze}. Soit $B\subset \Gamma$ l'ensemble des éléments
de couleur bronze.

Il s'agit de construire un block design équilibré (\S \ref{subs:hypdes})
qui contredise l'inégalité de Fisher (Prop.~\ref{prop:fishy}).
Définissons $A_b = \{y\in \Gamma : c(b y) = \text{azur}\}$ pour $b\in B$.
Comme $x\in B$ et  $c(x y) = \text{azur}$ pour au moins un $y\in C$,
et $c(x y)$ connaît la couleur de $y$, tous les éléments de $A_b$ sont
carmin.

Par la cohérence de $\mathfrak{X}$ et la
définition des nombres d'intersection (Def.~\ref{def:coh}),
\[|A_b| = \gamma(\text{azur},\text{azur}^{-1},\text{bronze}),\]
et donc $|A_b|$ ne dépend pas de $b$. Comme nous l'avons dit au d\'ebut,
$1\leq |A_x|< |C|$; donc, $1\leq |A_b|<C$ pour tout $b\in B$.

Montrez de façon similaire que, pour $v\in C$, la taille de
$\{b\in B: v\in A_b\} = \{b\in B: c(b v) = \text{azur}\}$ ne dépend pas
de $b$. Comme $\mathfrak{X}\lbrack C\rbrack$ est une clique,
$c(v,v')$ est de la même couleur pour tous $v,v'\in C$, $v=v'$; appelons
cette couleur {\em doré}. Montrez que
\[\{b\in B: v,v'\in A_b\} = \gamma(\text{azur},\text{azur}^{-1},\text{doré}).\]
Alors $(C,\{A_b\}_{b\in B})$ est un block design équilibré incomplet.

En conséquence, par l'inégalité de Fisher, $|B|\geq |C|$. Or, nous savons
que
$|C|>|\Gamma|/2$, $B,C\subset \Gamma$ et $B\cap C = \emptyset$. Contradiction.
\qed


\noindent{\sc Preuve du Lemme des designs} (Prop. \ref{prop:design}) ---
Supposons que
 pour chaque $\vec{x} \in \Omega^\ell$,
$0\leq \ell < k$,
$C_{\vec{x}}$ a une couleur $\alpha$-dominante $C(\vec{x})$, et, en plus,
si $\ell<k-1$, $(\mathfrak{X}_{\vec{x}})^{(2)}\lbrack C\rbrack$
est une clique. Nous arriverons à une contradiction.

Soit $C = C(\text{vide})$. Comme $|C|>\alpha |\Gamma|$,
$C$ est trop grande pour être un ensemble de jumeaux. Donc il existe
$u,v\in C$, $u\ne v$, tels que $\tau = (u v) \notin \Aut(\mathfrak{X})$.
Soit $\vec{y}$ de longueur minimale $r$ entre les chaînes satisfaisant
$c(\vec{y}^\tau) \ne c(\vec{y})$. Par cette minimalité et les règles
dans la définition \ref{def:ukuru}, $y_1,\dotsc y_r$ sont tous distincts.
En les permutant, nous pouvons assurer que
$u,v\notin \{y_1,y_2,\dotsc,y_{r-2}\}$, et, sans perte de généralité, que
soit (i) $y_{r-1}\ne u,v$, $y_r=u$, soit (ii) $y_{r-1} = u$, $y_r=v$. Dans le
cas (i), nous choisissons $\vec{x} = y_1,\dotsc,y_{r-1}$, $\ell=r-1$,
et voyons que $c_{\vec{x}}(u) \ne c_{\vec{x}}(v)$; dans le cas (ii), nous
choisissons $\vec{x} = y_1,\dotsc,y_{r-2}$, $\ell=r-2$, et obtenons
$c_{\vec{x}}(u,v) \ne c_{x}(v,u)$.
Nous aurons donc une contradiction avec notre supposition
une fois que nous aurons prouvé que $u,v\in C(\vec{x})$. 

Le fait que $u,v\in C(\vec{x})$ s'ensuivra immédiatement
de l'égalité $C(\vec{x}) = C\setminus \{x_1,\dotsc,x_\ell\}$; cette égalité, à son tour, 
se déduit par itération du fait que, pour $\vec{y}$ de longueur $\leq k-2$
et $\vec{x} = \vec{y} z$, $z\in \Omega$,
\begin{equation}\label{eq:kurut}
  C(\vec{x}) = C(\vec{y})\setminus \{z\}.\end{equation}
Pourquoi (\ref{eq:kurut}) est-il vrai? Nous sommes en train de supposer
que $\mathfrak{X}_{\vec{y}}^{(2)}\lbrack C(\vec{y})\rbrack$ est une clique,
et que $|C(\vec{y})|>\alpha |\Gamma|\geq |\Gamma|/2$.
Donc, par le lemme de la grande
  clique, tous les éléments de $C(\vec{y})$ sont des jumeaux en
  $\mathfrak{X}_{\vec{y}}^{(2)}$. En particulier, pour
  $u\in C(\vec{y})\setminus \{z\}$, $c_{\vec{x}}(u) = c_{\vec{y}}(z u)$
  ne dépend pas de $u$. Puisque le coloriage de sommets en $C_{\vec{x}}$
  est un raffinement de celui en $C_{\vec{y}}$ (par la deuxième règle de la
  définition~\ref{def:ukuru}), il s'ensuit que, soit
  $C(\vec{x}) = C(\vec{y})\setminus \{z\}$, soit
  $C(\vec{x})\subset \Gamma \setminus C(\vec{y})$, soit
  $C(\vec{x}) = \{z\}$. Comme
  $|C(\vec{x})|, |C(\vec{y})| > \alpha |\Gamma| \geq |\Gamma|/2$,
  les deux dernières possibilités sont exclues.
\qed


Nous appliquons le Lemme des designs (avec $\alpha=1/2$) à la configuration cohérente $k$-aire $\mathfrak{X}' = F_3(F_2(F_1(\mathfrak{X})))$, où
$\mathfrak{X}$ est donnée par $\mathbf{x}$
de la façon décrite au début de la section. Nous parcourons tous les
tuples possibles $\vec{x} = (x_1,\dotsc,x_{\ell}) \in \Gamma^\ell$,
$0\leq \ell<k$,
jusqu'à trouver un tuple pour lequel la première ou la deuxième conclusion
du Lemme des designs est vraie.

Si la première conclusion est vraie, nous
définissons $c_{\mathfrak{X}} = \mathfrak{X}_{\vec{x}}^{(1)}$ et
sautons à la section \ref{subs:cacoudo}.
Si la deuxième conclusion est vraie, nous passons au \S \ref{sec:coujoh},
ayant défini $\mathfrak{X}'' = \mathfrak{X}_{\vec{x}}^{(2)}\lbrack C\rbrack$,
où $C$ est la couleur $\alpha$-dominante de $\mathfrak{X}_{\vec{x}}^{(1)}$.

\subsection{Coupe ou Johnson}\label{sec:coujoh}

Nous avons une configuration classique cohérente homogène non triviale $\mathfrak{X}'' = (\Gamma,c)$. (Nous rappelons que ceci est un coloriage $c$
du graphe complet sur $\Gamma$ tel que (a) les sommets ont leur couleur
propre ({\og couleur diagonale\fg}),
(b) les arêtes $(x,y)$, $x\ne y$,
ne sont pas toutes de la même couleur, (c) la couleur $c(x,y)$
de l'arête $(x,y)$ détermine $c(y,x)$,
et (d) l'axiome de cohérence (\ref{def:coh}) se vérifie.)
Nous voudrions trouver des structures qui dépendent canoniquement
de $\mathfrak{X}''$ et qui contraignent son groupe d'automorphismes.

Il est raisonnable de s'attendre \`a ce que de telles structures existent: par le
Théorème \ref{thm:cam}, si le groupe
d'automorphismes est transitif, soit il est imprimitif (et donc il laisse
une partition invariante), soit il est près d'être $\Alt_m^{(k)}$, $k\geq 2$,
(qui laisse invariant un schéma de Johnson), soit il est petit
(et donc le stabilisateur de quelques points aura des orbites petites, et
ainsi nous donnera un coloriage sans couleur dominante). Le défi est de trouver
de telles structures, et de le faire canoniquement.

Si $\mathfrak{X}''$ n'est pas primitif (Déf.~\ref{def:uniprim}),
la tâche est plutôt facile: soit $r$ la couleur non diagonale la plus rouge
telle que le graphe  $\mathscr{G}_r = \{(x,y): x,y\in \Gamma,\; c(x,y)=r\}$
est connexe; par l'exercice \ref{ex:samesize}, ceci donne une partition
de $\Gamma$ dans des ensembles de la même taille $\leq |\Gamma|/2$. 

\begin{theo}[Coupe ou Johnson]\label{thm:coupjohn}
  Soit $\mathfrak{X} = (\Gamma,c)$ une configuration classique cohérente
  uniprimitive. Soit $2/3\leq \alpha< 1$.
  En temps $|\Gamma|^{O(1)}$, nous
  pouvons trouver
  \begin{itemize}
  \item soit un $\alpha$-découpage de $\Gamma$,
  \item soit un schéma de Johnson plongé sur
    $\Gamma_0\subset \Gamma$, $|\Gamma_0|\geq \alpha |\Gamma|$,
  \end{itemize}
  et un sous-groupe $H<\Sym(\Gamma)$ avec \[\lbrack \Sym(\Gamma):H\rbrack
  = |\Gamma|^{O(\log |\Gamma|)}\]
  tel que le découpage, voire le schéma,
  est canonique en relation à $H$.
  \end{theo}
Le groupe $H$ sera défini comme un stabilisateur de points
en $\Gamma$.
La valeur $2/3$ dans l'énoncé
est assez arbitraire; toute valeur $>1/2$ serait valable. Une valeur proche
à $1/2$ affecterait les constantes implicites.

\noindent{\sc Preuve } ---
Choisissons un $x\in \Gamma$ arbitraire. Donnons à chaque $y\in \Gamma$
la couleur de $c(x,y)$. Ce coloriage est canonique en relation à
$G_x$.  
S'il n'y a aucune classe de couleur $C$ de taille
$>\alpha |\Gamma|$, la partition triviale (non-partition)
de chaque classe nous donne un $\alpha$-découpage de $\Gamma$, et nous
avons fini.


Supposons, par contre,
qu'il y ait une classe de couleur -- disons, $C_{\text{lin}}$ --
de taille
$> \alpha |\Gamma|$. Comme $\alpha n>n/2$,
la relation $R_{\text{lin}}$ de cette couleur
est non orientée ($c(y,z) = \text{lin}$ ssi $c(z,y) = \text{lin}$).
Le complément de $R_{\text{lin}}$ (ou de toute autre relation)
est de diamètre $2$ (exercice). Soient $x,z\in \Gamma$ tels que
$c(x,z)=\text{lin}$, et soit $y\in \Gamma$ tel que
$c(x,y), c(z,y)\ne \text{lin}$. Appelons $c(x,y)$ {\em bis} et
$c(z,y)$ {\em terre}.

Considérons le graphe biparti $(V_1,V_2;A)$ avec sommets
$V_1 = C_{\text{lin}}$, $V_2 = C_{\text{bis}}$ et arêtes $R_{\text{terre}}\cap
(V_1\times V_2)$.
Le graphe est non vide par définition et semirégulier par l'exercice \ref{ex:biparhom}\ref{it:bara2}.
Par homogénéité et cohérence, le nombre de $y$ tels que
$c(y,w)$ est d'une couleur donnée $c_0$ est indépendant de $w$. Donc,  il est
toujours $\leq (1-\alpha) n < n/2$ pour $c_0\ne \text{lin}$.
Appliquant ceci à $c_0 = \text{terre}$ et $V_2$, nous voyons que le degré
$|\{v_1\in V_1: (v_1,v_2)\in A\}|$ est $<n/2$, et donc, comme $|V_1|>n/2$,
le graphe n'est pas complet. 

Nous appliquons donc la Proposition \ref{prop:bicoup} à $(V_1,V_2;A)$ avec
$\beta = \alpha |\Gamma|/|V_1|$. Notons que \mbox{$|V_2|\leq \beta |V_1|$.}
\qed

Nous travaillerons donc avec un graphe biparti $(V_1,V_2;A)$.
La stratégie sera d'essayer, soit de rendre $V_2$ plus petit (par au moins
un facteur constant), soit de trouver des structures en lui.
Soit ces structures nous permettront de réduire $V_2$ quand même, soit
elles nous aideront à découper $V_1$, ou à trouver un schéma de
Johnson assez grand sur~$V_1$.

Tout d'abord,
nous devrons borner la symétrie en $V_1$, c'est-à-dire réduire, voire
éliminer les jumeaux. Il y a deux raisons à ceci.
\begin{itemize}
\item Même si nous découvrions une structure assez riche en $V_2$, cela
  impliquerait peu ou rien sur $V_1$ si beaucoup d'éléments de $V_1$ se
  connectent à $V_2$ de la même façon.
\item  Si $V_2$ est petit,
  nous colorierons chaque sommet de $V_1$ par son ensemble de voisins en $V_2$.
Ceci nous donnera un coloriage canonique en relation à $G_{(V_2)}$.
Or, dans ce coloriage,
deux sommets en $V_1$ auront la même couleur ssi ils sont des jumeaux;
donc, si aucune classe de jumeaux en $V_1$ n'a $>\alpha |V_1|$ 
éléments, nous aurons un $\alpha$-coloriage.
\end{itemize}

\begin{exer}\label{ex:gita}
  Soit $(V_1,V_2;A)$ un graphe biparti semirégulier et non trivial. Alors,
  aucune classe de jumeaux en $V_1$ n'a plus de $|V_1|/2$ éléments.
\end{exer}
\noindent{\sc Solution} ---
  Nous assurons que $|A|\leq |V_1| |V_2|/2$ en prenant le complément
  s'il est
  nécessaire. Soient $d_2$ le degré des sommets en $V_2$ et $S$ une classe de jumeaux en $V_1$.
Montrez que $d_2\geq |S|$, et donc  $|A|\geq |S| |V_2|$. 
  \qed
  \begin{exer}\label{ex:reduc2}
    Soit $(V_1,V_2;A)$ un graphe biparti sans jumeaux en $V_1$.
    Soient $V_2 = C_1\cup C_2$, $C_1\cup C_2=\emptyset$.
    Montrez que, pour au moins un $i=1,2$,
    il n'y a aucune classe de $\geq |V_1|/2+1$ jumeaux en $V_1$ 
    dans le graphe $(V_1,C_i;A\cap (V_1\times C_i))$.
  \end{exer}
  \begin{exer}\label{ex:twindet}
    Soit $\mathfrak{X}=(\Gamma,c)$ une configuration cohérente. Soient
    $C_1$, $C_2$ deux classes de couleurs en $\Gamma$. Soit {\em brun} une
    couleur d'arêtes en $A\times B$. Alors, pour $x,y\in C_1$,
    la couleur $c(x,y)$ détermine si $x$ et $y$ sont des jumeaux dans
    le graphe biparti $(C_1,C_2;\mathscr{G}_{\text{brun}})$.
\end{exer}
  
\begin{prop}[Coupe ou Johnson biparti, ou {\og Una partita a poker\fg}]\label{prop:bicoup}
  Soit $X = (V_1,V_2;A)$ un graphe biparti 
  avec $|V_2|<\beta |V_1|$, où $2/3\leq \beta<1$, et tel qu'aucune classe
  de jumeaux en $V_1$ n'ait plus de $2|V_1|/3$ éléments.
  Alors, nous pouvons trouver, en temps $|V_1|^{O(1)}$,
  \begin{itemize}
  \item soit un $\beta$-découpage de $V_1$,
  \item soit un schéma de Johnson plongé sur
    $V_0\subset V_1$, $|V_0|\geq \beta |V_1|$,
  \end{itemize}
  et un sous-groupe $H<G$, $G = \Sym(V_1)\times \Sym(V_2)$, avec
  \[\lbrack G:H\rbrack = |V_1|^{O(\log |V_1|)}\]
  tel que le découpage, voire le schéma,
  est canonique en relation à $H$.
\end{prop}
La condition sur les classes de jumeaux ici 
était remplie (même avec $1/2$ à la place de $2/3$) à la fin de la preuve
du Thm.~\ref{thm:coupjohn}, grâce à l'exercice \ref{ex:gita}.

En ce qui concerne le temps de la procédure, nous expliciterons quelques
détails qui pourraient ne pas être évidents. Ce qui sera le détail le
plus délicat est l'indice $\lbrack G:H\rbrack$. Le groupe $H$ sera défini comme un
stabilisateur de points; nous devons bien contrôler le nombre de points
que nous stabilisons. 

Esquissons la stratégie générale de la preuve. Ce que nous voulons est
une réduction à la Proposition \ref{prop:cohcoup},
``Coupe-ou-Johnson cohérent''. Nous pouvons produire une configuration
cohérente classique sur $V_1\cup V_2$ à partir du graphe $X$,
tout simplement en utilisant Weisfeiler-Leman. Ce qui demande de la ruse est de
garantir que la restriction $\mathfrak{X}\lbrack C_2\rbrack$ à la classe
de couleurs dominante (s'il y a une) soit non triviale.

Pour obtenir une configuration non-triviale sur $C_2$, nous noterons que le
graphe $X$ induit lui-même une relation $d$-aire sur $C_2$, où $d$ est au plus
le degré de la majorité d'éléments de $V_1$ (si telle chose existe; sinon,
les degrés nous donnent une partition de $V_1$). Si la relation est triviale,
dans le sens de contenir toutes les $d$-tuples d'éléments distincts dans $C_2$,
nous obtenons un schéma de Johnson.  Si elle est non triviale mais contient
beaucoup de jumeaux, elle nous donne une manière de
descendre à un $C_2$ plus petit. S'il n'y a pas beaucoup de jumeaux, nous
utilisons le Lemme des Designs (supplémenté par un lemme standard sur les designs) pour obtenir une configuration cohérente classique non-triviale sur $C_2$,
ce qui était à trouver.

\noindent{\sc Preuve } ---
Si $|V_1|\leq c$, où $c$ est une constante, nous colorions chaque $v\in V_1$
par lui-même. Ce coloriage est canonique en relation à $H=\{e\}$; autrement dit, il n'est pas canonique du tout. Peu importe: 
trivialement, $|G|\leq (c!)^2 \leq |V_1|^{O(\log |V_1|)}$. Nous pouvons donc
supposer que $|V_1|>c$.

Si $|V_2|\leq (6 \log |V_1|)^{3/2}$ (disons), alors, par la discussion
ci-dessus, nous obtenons un $(2/3)$-coloriage de $V_1$
(et donc: un $(2/3)$-découpage de $V_1$). Ce coloriage est canonique
en relation à un $H$ d'indice
\[|V_2|! \leq |V_2|^{|V_2|} \leq (6 \log |V_1|)^{\frac{3}{2}
   (6 \log |V_1|)^{\frac{3}{2}}} \ll |V_1|^{\log |V_1|}.\]
Nous pouvons donc supposer que $|V_2|> (6 \log |V_1|)^{3/2}$.

Notre première tâche est d'éliminer les jumeaux.
Nous divisons $V_1$ dans ses classes de jumeaux et colorions chaque
$v\in V_1$ par son nombre de jumeaux et par son degré dans le graphe
$(V_1,V_2;A)$. Nous obtenons
un $\beta$-découpage de $V_1$, sauf s'il y a un entier $d_1$ tel que
l'ensemble $V_1'$ des sommets $v$ sans jumeaux et de degré $d_1$ est de taille
$|V_1'|> \beta |V_1|$. Supposons dorénavant que cela est le cas. Comme $|V_1'|
> |V_2|$ et qu'il n'y a pas de jumeaux en $V_1'$, nous voyons que $1<d_1<|V_2|-1$;
nous pouvons supposer que $d_1\leq |V_2|/2$ en remplaçant $A$ par son complément,
si nécessaire.

Soit $\mathscr{H} = (V_2,\mathscr{A})$ l'hypergraphe dont les arêtes sont les
voisinages en $(V_1,V_2;A)$ des sommets dans $V_1'$.
(Elles sont toutes contenues en $V_2$.) L'hypergraphe
est $d_1$-uniforme.
Comme il n'y a pas de jumeaux dans $V_1'$, il n'y
a pas d'arêtes identiques. Si $\mathscr{H}$ est l'hypergraphe
complet $d_1$-uniforme, alors $V_1'$ peut être identifié avec le schéma
de Johnson $\mathscr{S}_{d_1}(V_2)$. {\em
  (Scoppia in un pianto angoscioso e abbraccia la
  testa di Johnson.)}

Supposons alors que $\mathscr{H}$ n'est pas complet.
Nous voudrions avoir un coloriage canonique sur $V_2^d$
pour un $d\ll l$, $l=(\log |V_1'|)/\log |V_2|$, tel que les
éléments de $V_2$ ne soient pas tous jumeaux.
Si $d_1\leq 6 \lceil l\rceil$, nous définissons
$d=d_1$ et colorions $\{(v_1,\dotsc,v_d)\in V_2^d: \{v_1,\dotsc,v_d\}\in \mathscr{H}\}$ 
en écarlate, et tout le reste en gris.

Supposons, par contre, que $d_1> 6 \lceil l\rceil$. Soit
$d = 6 \lceil l \rceil$. Nous colorions $\vec{v}=(v_1,\dotsc,v_d)$ en gris
si les $v_i$ ne sont pas tous distincts; dans le cas contraire, nous
donnons à $\vec{v}$ la couleur
\begin{equation}\label{eq:juj}
  |\{H\in \mathscr{H}: \{v_1,\dotsc,v_d\} \subset H\}|.\end{equation}
Cette opération de coloriage peut être faite en temps de l'ordre de 
\[|V_1|\cdot \binom{d_1}{d} \leq |V_1|\cdot |V_2|^d =
|V_1| \cdot |V_2|^{6 \left\lceil \frac{\log |V_1'|}{\log |V_2|}
  \right\rceil} = |V_1|\cdot |V_1'|^{O(1)} = |V_1|^{O(1)}.\]

Si les tuples avec $v_1,\dotsc,v_d$ distincts n'avaient pas tous
la même couleur $\lambda$,
nous aurions un design $d-(|V_2|,d_1,\lambda)$ avec $|V_1'|$ arêtes. Donc,
par la Proposition~\ref{prop:RW}, $|V_1'|\geq \binom{|V_2|}{s}$ pour $s =
3\lceil l\rceil$. Comme $|V_2|\geq (6 \log |V_1'|)^{3/2}$ et $|V_1'|$ peut être supposé plus grand qu'une constante, 
\[\binom{|V_2|}{s} \geq 
\left(\frac{|V_2|}{s}\right)^s > \left(\frac{|V_2|}{6 l}\right)^s > 
\left(|V_2|^{1/3}\right)^{3 l}
= |V_2|^{\frac{\log |V_1'|}{\log |V_2|}}
= |V_1'|,\]
ce qui donne une contradiction.

Donc, pour $d_1$ arbitraire,
les tuples avec $v_1,\dotsc,v_d$ distincts
n'ont pas tous la même couleur; en d'autres termes,
les éléments de $V_2$ ne sont pas tous jumeaux en relation à notre
nouvelle structure $d$-aire.
S'il y a une classe $S$ de jumeaux de taille $> |V_2|/2$,
alors, par l'exercice \ref{ex:reduc2}, au moins un des deux graphes
$(V_1',S;A\cap (V_1'\times S))$, $(V_1',V_2\setminus S; A\cap (V_1'\times
(V_2\setminus S)))$
n'a aucune classe de $> |V_1'|/2+1$ jumeaux dans $V_1$. Comme $2 |V_1'|/3\geq
|V_1'|/2+1$ pour
$|V_1|\geq 8$, nous appliquons la Proposition~\ref{prop:bicoup}
elle-même à un de ces deux graphes (disons, celui sur $V_1'\times S$ si
les deux sont valables), et terminons. (Peut-être que la taille de $V_2$ est
descendue seulement à $|V_2|-1$, mais tous nos choix ont été canoniques
-- des non-choix, si l'on veut -- donc gratuits. Nous n'avons perdu
que du temps; pour être précis, $|V_1|^{O(1)}$ de temps, ce qui est acceptable.)

Alors, nous avons un coloriage de $V_2^d$ en relation auquel il n'y a aucune
classe de jumeaux en $V_2$ de taille $> |V_2|/2$. Nous appliquons les
foncteurs $F_1$, $F_2$, $F_3$ (Weisfeiler-Leman) à ce coloriage. Puis nous
utilisons le Lemme des designs (Prop.~\ref{prop:design}) avec $\alpha=2/3$.
Nous trouvons les éléments $x_1,\dotsc,x_\ell\in V_2$ ($\ell=d-1$ ou $\ell=d-2$)
dans l'énoncé de la Proposition~\ref{prop:design} par
force brute, en temps proportionnel à $|V_2|^d = |V_1|^{O(1)}$. Nous les fixons, et 
nous imposons que
$H$ fixe $x_1,\dotsc,x_\ell$, ce qui a un coût de $|V_1|^{O(1)}$, dans le
sens où \[\lbrack G: G_{x_1,\dotsc,x_\ell}\rbrack \leq |V_2|^d
= |V_1|^{O(1)}.\]




Si nous sommes dans le premier cas du Lemme de designs
(pas de couleur dominante), nous cueillons les classes de couleur,
en commençant par la plus rouge (interprétez la quantité en (\ref{eq:juj}) comme une longueur
d'onde), jusqu'à avoir une union des classes $S\subset V_2$ avec
$ |V_2|/3<|S|\leq 2 |V_2|/3$. (Ceci marche s'il n'y a aucune classe de taille
$> |V_2|/3$; si telles classes existent, nous définissons $S$ comme la classe la
plus grande de ce type.)
Nous appliquons l'exercice \ref{ex:reduc2},
et obtenons un graphe $(V_1',V_2',A\cap (V_1'\cap V_2'))$
remplissant les conditions de notre Proposition~\ref{prop:bicoup}
avec $V_2' = S$ ou $V_2' = V_2\setminus S$, et donc $|V_2'|\leq \alpha |V_2|$.
Donc, nous appliquons la Proposition~\ref{prop:bicoup} à ce graphe;
la récursion marche. (Il est important ici que $|V_2'|\leq \alpha |V_2|$,
puisque nous avons déjà encouru un coût considérable ($|V_1|^{O(1)}$)
  dans l'indice.)

  Restons donc dans le deuxième cas du Lemme des designs: nous avons
  un coloriage de $V_2$ avec
  une classe de couleurs $C\subset V_2$ telle que
  $|C|\geq 2 |V_2|/3$, et une
  configuration cohérente homogène
  classique $\mathfrak{Y}$ non triviale sur $C$. 

  Nous définissons un graphe avec
  des sommets $V_1'\cup V_2$, où $V_1'$ est de couleur nacrée et
  $V_2$ vient d'être colorié par le Lemme des Designs; les arêtes
  seront non pas seulement celles en
$A\cap (V_1'\times V_2)$ (coloriées en noir) mais aussi les arêtes
entre les éléments de $V_2$, dans les couleurs données par $\mathfrak{Y}$.
Nous appliquons les raffinements 
$F_1$, $F_2$ et $F_3$ (Weisfeiler-Leman) à ce graphe, et obtenons
une configuration cohérente $\mathfrak{X}$.


La configuration $\mathfrak{X}\lbrack V_2\rbrack$ est un raffinement de $\mathfrak{Y}$.
Si elle n'a pas de couleur $\alpha$-dominante, nous réduisons notre problème à
celui pour $(V_1',V_2',A\cap (V_1'\cap V_2'))$, $|V_2'|\leq \alpha |V_2|$,
comme avant; nous pouvons appliquer la proposition \ref{prop:bicoup}
\`a un tel graphe sans changer $\beta$ parce que
\[|V_2'|\leq \frac{2}{3} |V_2|\leq \beta |V_2| < \beta |V_1'|.\]
La r\'ecursion marche ici aussi parce que $|V_2'|\leq 2 |V_2|/3$:
il est important que $V_2'$ soit plus petit que $V_2$ par un facteur constant,
puisque le coût entraîn\'e jusqu'\`a maintenant  
dans l'index $\lbrack G:H\rbrack$ est d\'ej\`a considérable ($|V_1|^{O(1)}$).

Supposons donc qu'il y a une classe de couleurs $(2/3)$-dominantes
$C_2$ dans $\mathfrak{X}\lbrack V_2\rbrack$. Elle doit être un sous-ensemble de $C$
car $2/3+2/3>1$.
La restriction de $\mathfrak{X}\lbrack C_2\rbrack$ n'est pas une clique: si elle l'était,
la restriction de $\mathfrak{Y}$ à $C_2$ l'aurait été aussi, et cela est
impossible par l'exercice~\ref{ex:noclique}.

Nous pouvons supposer qu'il existe une classe de couleurs
$C_1\subset V_1'$ en $\mathfrak{X}_1$ qui
satisfait $|C_1|> \beta |V_1|$; sinon,
nous avons un $\beta$-coloriage de $V_1$, et pouvons finir.
Le fait que $|C_1|> \beta |V_1|$ implique que $|C_1|> |V_2|\geq |C_2|$.


Nous pouvons supposer aussi que les arêtes de $\mathfrak{X}$ en
$C_1\times C_2$ ne sont pas toutes de la même couleur. Si elles l'étaient,
il y aurait
une classe de $\geq |C_1|> \beta |V_1| > |V_1|/2+1$ jumeaux en $V_1$ dans le
graphe $(V_1,C_2;A\cap (V_1\times C_2))$, dont
$\mathfrak{X}\lbrack V_1\times C_2\rbrack$ est un raffinement.
Dans ce cas, par l'exercice~\ref{ex:reduc2},
nous aurions une réduction à
$(V_1,V_2\setminus C_2;A\cap (V_1\times (V_2\setminus C_2)))$, et nous
pourrions finir en utilisant la Proposition \ref{prop:bicoup} de façon
récursive.

Ainsi, nous avons tout réduit à la Proposition \ref{prop:cohcoup}:
nous l'appliquons à $\mathfrak{X}\lbrack C_1\cup C_2\rbrack$.
Nous obtenons, soit un $(1/2)$-découpage de $C_1$, soit un graphe
biparti $(W_1,W_2; A')$, $W_1\subset C_1$, $W_2\subset C_2$, avec
$|W_1|\geq |C_1|/2$, $|W_2|\leq |C_2|/2 \leq |V_2|/2$, tel qu'aucune classe
de jumeaux en $W_1$ n'a plus que $|W_1|/2$ éléments. Nous pouvons supposer
que $|W_1|>\beta |V_1|$, parce que, dans le cas contraire, nous avons obtenu
un $\beta$-découpage de $W_1$.
Alors, $|W_2|< |W_1|/2$. Nous pouvons, alors, faire
de la récursion: nous appliquons la Proposition \ref{prop:bicoup} avec
$(W_1,W_2;A')$ à la place de $(V_1,V_2;A)$.

La récursion se finit après pas plus que $O(\log |V_2|)$ pas
puisque $|W_2|\leq |V_2|/2$. Si la taille de $W_1$ (ou de $V_1$)
décroît en dessous de $\beta |V_1|$ (pour la valeur originale de $|V_1|$),
alors nous avons obtenu un $\beta$-découpage de $V_1$.
\qed

Comme nous l'avons vu, Coupe ou Johnson biparti utilise Coupe ou Johnson
cohérent. \`A son tour, Coupe ou Johnson cohérent se réduira à
Coupe ou Johnson biparti pour un graphe biparti $(V_1,V_2;A)$
avec $V_2$ de taille au plus une moitié de la taille du $V_2$ original.

\begin{prop}[Coupe ou Johnson cohérent]\label{prop:cohcoup}
  Soit $\mathfrak{X} = (C_1\cup C_2;c)$ une configuration cohérente avec
  des classes de couleurs de sommets $C_1$, $C_2$, où $|C_1|>|C_2|$.
  Supposons que ni $c|_{C_1\times C_2}$ ni $c|_{C_2\times C_2}$ est une fonction
  constante.

  Alors, nous pouvons trouver, en temps $|C_1|^{O(1)}$, soit
  \begin{itemize}
  \item un $(1/2)$-découpage de $C_1$, ou
  \item un graphe biparti $(V_1,V_2;A)$, $V_i\subset C_i$,
    $|V_1|\geq |C_1|/2$, $|V_2|\leq |C_2|/2$, tel que toute classe de jumeaux
    en $V_1$ contient au plus $|V_1|/2$ éléments,
  \end{itemize}
  et un élément $y\in C_2$, tel que le découpage, voire le graphe biparti,
  est canonique en relation à $G_y$, où $G=\Sym(C_1)\times \Sym(C_2)$.
\end{prop}
Il va de soi que dire que $c|_{C_2\times C_2}$ est constant équivaut à dire
que $\mathfrak{X}\lbrack C_2\rbrack$ est une clique.
\noindent{\sc Preuve } ---

Si la restriction $\mathfrak{X}\lbrack C_1\rbrack$ était une clique, alors,
par cohérence, pour toute couleur en $C_1\times C_2$ -- pourpre, disons --
les voisinages dans $(C_1,V_2;\mathscr{G}_{\text{pourpre}})$
des sommets en $C_2$ nous donneraient un block design
équilibré (et peut-être dégénéré) sur $C_1$. Le design est incomplet parce
que $c$ n'est pas monochrome sur $C_1\times C_2$.
L'inégalité de Fisher nous donne que $|C_2|\geq |C_1|$, en contradiction
avec nos suppositions.
Donc, $\mathfrak{X}_1\lbrack C_1\rbrack$ n'est pas une clique.

Si $\mathfrak{X}\lbrack C_1\rbrack$ n'est pas primitive, la plus rouge
de ses relations non connexes nous donne un $(1/2)$-découpage canonique
de $V_1$, par l'exercice \ref{ex:samesize}. Nous pouvons donc supposer
que $\mathfrak{X}\lbrack C_1\rbrack$ est primitif.


Nous avons deux cas à considérer: $\mathfrak{X}\lbrack C_2\rbrack$ primitive et
$\mathfrak{X}\lbrack C_2\rbrack$ imprimitive.

Supposons d'abord que $\mathfrak{X}\lbrack C_2\rbrack$ est imprimitive.
%
La relation non connexe la plus rouge dans $\mathfrak{X}\lbrack C_2\rbrack$
nous donne une partition de
$C_2$ dans des ensembles $B_1,\dotsc,B_m$, $m\geq 2$, tous
  de la même taille $\geq 2$. Nous avons donc trouvé une structure en $C_2$,
  et nous l'utiliserons, soit pour découper $C_1$, soit pour réduire
  $|C_2|$ par un facteur constant. Le premier pas consiste
  à montrer qu'il n'y a pas de jumeaux dans $C_1$.
  

 
Comme notre configuration est cohérente, la couleur
d'une arête en $C_1$ sait si ses sommets sont des jumeaux en relation
à $C_2$ (Ex.~\ref{ex:twindet}); donc, s'il y
avait des jumeaux dans $C_1$ en relation à $C_2$, nous aurions, soit qu'une
des couleurs d'arêtes en $C_1$ donne une relation non connexe -- ce qui
contredit le fait que $\mathfrak{X}\lbrack C_1\rbrack$ est uniprimitive -- soit que tous les
éléments de $C_1$ sont des jumeaux en relation en $C_2$.
Dans ce dernier cas, par l'exercice~\ref{ex:gita}, $c|_{C_1\times C_2}$
serait monochrome, ce qui n'est pas le cas. En conclusion, il n'y a pas
de jumeaux dans $C_1$ en relation à $C_2$.

Notre intention est d'appliquer l'exercice \ref{ex:biparcomp} pour obtenir un graphe biparti contracté $C_1\times \{1,2,\dotsc,m\}$ avec $m\leq |C_2|/2$.
Nous devons seulement faire attention à ce
que ce graphe ne soit pas trivial.

Soit $d_k$ le degré de tout $w\in C_2$ dans le graphe biparti
$(C_1,C_2;\mathscr{G}_k)$ pour une couleur~$k$ donnée, où $\mathscr{G}_k$
consiste en les arêtes de couleur $k$.
(Par l'ex.~\ref{ex:biparhom}\ref{it:bar1},
le degré $d_k$ ne dépend pas de $w$.) Si $d_k\leq
|C_1|/2$ pour tout $k$, nous fixons un $w\in C_2$ (non canonique)
et obtenons un $(1/2)$-coloriage de $C_1$ en assignant la couleur
$c(x,w)$ au sommet $x\in C_1$.  Supposons
donc qu'il y a une couleur -- que nous appellerons {\em violet} -- telle que
$d_{\text{violet}}> |C_1|/2$. S'il y a un $1\leq i\leq m$ tel qu'il n'y a aucune
classe de plus que $|C_1|/2$ jumeaux dans $C_1$ en relation à $B_i$,
nous fixons un élément $y\in B_i$ d'un tel $i$ (non canoniquement), fixant ainsi
cet $i$. De cette façon, nous obtenons une réduction
au graphe biparti $(C_1,B_i;\mathscr{G}_{\text{violet}}\cap (C_1\times B_i))$.

Supposons que cela n'est pas le cas. Donc, pour chaque $i$, il existe une
classe $T_i\subset C_1$ de jumeaux en relation à $B_i$ telle que
$|T_i|>|C_1|/2$. Pour chaque $w\in B_i$, les arêtes de $w$ à tout $v\in T_i$
sont de la même couleur; alors, elles doivent être violettes. Soit {\em vert}
une couleur d'arêtes en $C_1\times C_2$ qui ne soit pas violet. Alors, le graphe
$X=(C_1,\{1,\dotsc,m\};D)$ dans l'exercice \ref{ex:biparcomp}
n'est pas vide; comme $(v_i,i)$ est violet pour tout $v\in T_i$, $X$ n'est
pas complet non plus. Comme $X$ est birégulier, il n'y a aucune classe
de jumeaux en $C_1$ en relation à $\{1,\dotsc,m\}$ avec $> |C_1|/2$
éléments (ex. \ref{ex:gita}). Nous avons donc tout réduit à
un graphe biparti $X$ du type que nous désirions.

Considérons maintenant le cas de $\mathfrak{X}\lbrack C_2\rbrack$ primitive\footnote{Le problème dans la preuve originale de Babai était à ce point précis. Ce qui suit est un argument alternatif proposé par lui ({\em col rumore sordo di un galoppo})
  lorsque cet article était en train d'être édité. Il est plus concis et
  élégant que l'argument d'origine, en plus d'être correct. Avant, la preuve
  faisait  deux fois (ou plus) recours à la proposition elle-même, ce qui faisait 
 croître  l'indice $\lbrack G:H\rbrack$ de fa\c{c}on catastrophique.}.
Fixons un $y\in C_2$ arbitraire
(non canoniquement).
Nous pouvons supposer qu'il y a une couleur -- disons, violet -- telle que
$d_{\text{violet}}> |C_1|/2$, puisque, sinon, les couleurs
des arêtes qui connectent les éléments de $C_1$ avec $y$ nous
donneraient un $(1/2)$-coloriage de $C_1$. \'Ecrivons
$V_1 = L_{\text{violet}}(y) = \{x\in C_1: c(x,y) = \text{violet}\}$.
Donc $|V_1|> |C_1|/2$.
Soit {\em bleu} une couleur d'arêtes en $\mathfrak{X}\lbrack C_2\rbrack$
telle que le degré de $\mathscr{G}_{\text{bleu}}$ est (positif et)
$< |C_2|/2$; une telle couleur existe parce que
$\mathfrak{X}\lbrack C_2\rbrack$ n'est pas une clique. (S'il y a
plusieurs couleurs comme cela, nous choisissons la plus bleue d'entre elles.)
Alors, $V_2 = L_{\text{bleu}}(y)\subset C_2$ satisfait  $1\leq |V_2|< |C_2|/2$.

Le graphe biparti $(V_1,V_2;\mathscr{G}_{\text{violet}}\cap (W\times U))$
est semirégulier par l'exercice \ref{ex:biparhom}\ref{it:bara2}. Il est
non vide parce que, pour tout $u\in V_2$,
$|L_{\text{violet}}(u)|> |C_1|/2$, et donc \mbox{$L_{\text{violet}}(u)\cap V_1 \ne \emptyset$.}
S'il était complet, nous aurions $V_1\subset L_\text{violet}(u)$ pour tout
$u\in V_2$; comme $|V_1| = |L_{\text{violet}}(y)| = |L_{\text{violet}}(u)|$, ceci
impliquerait que $V_1=L_{\text{violet}}(u)$. Or, cela voudrait dire que
$y$ et $u$ sont des jumeaux dans le graphe $(C_1,C_2;\mathscr{G}_{\text{violet}})$.
Par le même argument qu'avant (basé sur l'exercice~\ref{ex:twindet}), la primitivité
de $\mathfrak{X}\lbrack C_2\rbrack$ et le fait que $c|_{C_1\times C_2}$ ne soit pas monochrome
impliquent qu'il n'y a pas de jumeaux dans $C_2$ en relation au graphe
$(C_1,C_2;\mathscr{G}_{\text{violet}})$.
Donc, $(V_1,V_2;\mathscr{G}_{\text{violet}}\cap (V_1\times V_2))$ n'est pas complet.
Par l'exercice \ref{ex:gita}, nous obtenons qu'aucune classe de jumeaux
dans $(V_1,V_2;\mathscr{G}_{\text{violet}}\cap (V_1\times V_2))$ n'a plus de $|V_2|/2$
éléments. Nous avons donc terminé.
\qed
 
\subsection{Récursion et réduction}\label{subs:recred}

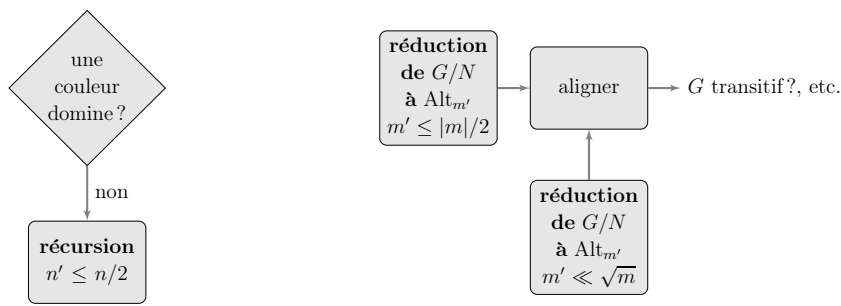
\begin{center}
      \scalebox{0.66}{
\begin{tikzpicture}[scale=2, node distance = 2cm, auto]
    \node [below of=init] (trans) {$G$ transitif?, etc.};
    \node [block, left of = trans, node distance = 3.5cm] (align) {aligner};
    \node [block, left of = align, node distance = 3cm] (redhalf)
                {\bf réduction de $G/N$ à $\Alt_{m'}$\\ $m'\leq |m|/2$};
    \node [block, below of = align, node distance = 3cm] (redsqrt)
          {\bf réduction de $G/N$ à $\Alt_{m'}$\\ $m'\ll \sqrt{m}$};
    \node [decision, left of = redhalf, node distance = 7cm] (coldom) {une couleur domine?};
    \node [block, below of = coldom, node distance = 3.5cm] (luks) {\bf récursion $n'\leq n/2$};       
      \path [line] (align) -- (trans);
      \path [line] (redhalf) -- (align);
      \path [line] (redsqrt) -- (align);
      \path [line] (coldom) -- node [, color = black] {non} (luks);
  \end{tikzpicture}
    }
\end{center}

\subsubsection{Le cas sans couleurs dominantes}\label{subs:cacoudo}

Nous sommes dans le cas dans lequel un coloriage
$c_{\mathfrak{X}}:\Gamma\to \mathscr{C}$ n'a pas de couleur
dominante. Ici $c_{\mathfrak{X}}$ est l'image d'une
structure $\mathfrak{X}$ sous un foncteur $F$
qui commute avec l'action de
$H_{(x_1,\dotsc,x_\ell)}$, où $H=\Alt(\Gamma)$, $x_i\in \Gamma$.
Le fait que $c_\mathfrak{X}$ n'a pas de couleur dominante
nous servira pour trouver ou écarter ses isomorphismes
possibles en $H_{\vec{x}}
= H_{(x_1,\dotsc,x_\ell)}$.
Pour trouver ou écarter des isomorphismes en
tout $H=\Alt(\Gamma)$, nous n'avons qu'à travailler avec
un ensemble de représentants $\{\sigma_1,\dots,\sigma_s\}$, $s\leq
|\Gamma|^\ell = m^\ell$, des classes de $H_{\vec{x}}$ dans $H$, et à
faire l'union de $\Iso_{H_{\vec{x}}}\left(c_{\mathfrak{X}},c_{\mathfrak{Y}_i}\right)$
pour $\mathfrak{Y}_i = \mathfrak{Y}^{\sigma_i^{-1}}$:
\begin{equation}\label{eq:kozron}
  \Iso_{H}(\mathfrak{X},\mathfrak{Y}) =
  \bigcup_{1\leq i\leq s} \Iso_{H_{\vec{x}}}\left(
  \mathfrak{X},\mathfrak{Y}_i\right)
\sigma_i,\;\;\;\;\;\;\;\;\;\;\;\;
\Iso_{H_{\vec{x}}}(\mathfrak{X},\mathfrak{Y}_i) \subset
\Iso_{H_{\vec{x}}}(c_{\mathfrak{X}},c_{\mathfrak{Y}_i}). 
\end{equation}
Ceci est similaire à l'équation (\ref{eq:rulu}), en \S \ref{subs:luks}.
Le coût de la procédure
est multiplié par $s\leq m^\ell$.


Si le coloriage $c_{\mathfrak{X}}$
n'est pas une permutation (en $H_{\vec{x}}$)
    du coloriage 
    $c_{\mathfrak{Y}_i}$, alors 
    $\Iso_{H_{\vec{x}}}(c_{\mathfrak{X}},c_{\mathfrak{Y}_i}) = \emptyset$.
    Supposons, par contre,
    qu'il y a au moins un $\tau_i\in H_{\vec{x}}$ 
    tel que 
    $c_{\mathfrak{X}} = c_{\mathfrak{Y}_i}^{\tau_i}$.
(Nous disons que $\tau_i$ {\em aligne}
    $c_{\mathfrak{X}}$ et $c_{\mathfrak{Y}_i}$.)
    Il est trivial de 
    trouver $\tau_i$. 
    Or
    \[
    \Iso_{H_{\vec{x}}}\left(\mathfrak{X},\mathfrak{Y}_i\right) =
    \Iso_{H_{\vec{x}}}\left(\mathfrak{X},\mathfrak{Y}_i^{\tau_i}\right) \tau_i^{-1}
    \subset \Aut_{H_{\vec{x}}}(c_{\mathfrak{X}}) \tau_i^{-1}. 
    \]
    Comme $c_{\mathfrak{X}}$ n'a pas de couleur dominante, ceci est assez contraignant,
    ce que nous voulions.
    
 Appliquons cette procédure générale au cas de $G$ primitif que nous sommes en
train de discuter. Il y a une bijection $\iota:\Omega\to
\{S\subset \Gamma: |S|=k\}$; donc, $c_{\mathfrak{X}}$ induit un coloriage
$c':\Omega \to \{(k_i)_{i\in \mathscr{C}}: k_i\geq 0, \sum_i k_i =k\}$.
Nous sommes dans une situation similaire à celle de la fin du
\S \ref{sec:grasym}, mais en mieux: il est facile de montrer que,
comme aucune classe de couleur de $c$ possède plus de $\alpha |\Gamma|$
    éléments,
    aucune classe de couleur de $c'$ possède plus de $\alpha |\Omega|$ éléments.

    Nous procédons alors comme dans le
cas intransitif de la preuve de Luks (Thm.~\ref{thm:luxor}), ce qui réduit
le problème à $\leq n$ problèmes d'isomorphisme de chaînes pour des chaînes de
longueur $\leq \alpha n$ et de longueur totale $\leq n$.
Le dernier pas ({\em lifting}, {\og relèvement\fg}) consiste à trouver des éléments de $G$
qui induisent $\tau_i$. \'Etant donnée une bijection $\iota$, ceci est trivial.

\subsubsection{Le cas du découpage}

Considérons maintenant
un $\alpha$-découpage (fin de \S \ref{subs:secf1}) d'un
ensemble de sommets $\Gamma$. Ce découpage sera donné canoniquement, à savoir,
en tant que l'image d'une structure $\mathfrak{X}$ sous un foncteur, tout
comme le coloriage au \S \ref{subs:cacoudo}. Nous pouvons supposer 
que le découpage a une classe de couleurs $C$ dominante
($|C| > \alpha |\Gamma|$, $\alpha>1/2$), puisque, dans le cas contraire, nous
pouvons passer au \S \ref{subs:cacoudo}.

Nous voulons savoir quels éléments de $\Alt(\Gamma)$ respectent le
$\alpha$-découpage; ceci nous aidera à contraindre les isomorphismes de
$\mathfrak{X}$, tout comme en (\ref{eq:kozron}).
Par la définition de $\alpha$-découpage, $C$ est partitionné en
$\ell\geq 2$ ensembles de la même taille $\geq 2$. Les seules permutations en
$\Alt_{m_0}$, $m_0 = C$, qui sont permises sont celles qui respectent la
partition. Le groupe qui respecte la partition est isomorphe à
$\Alt_{m_0/\ell}$.

Nous avons donc  réduit notre problème à un problème
avec $m'=m_0/\ell\leq m/2$. Après avoir r\'esolu ce problème, nous travaillons
-- comme dans le  \S \ref{subs:cacoudo} --
sur les autres classes de couleurs.


\'Etant donnés deux $\alpha$-découpages, nous vérifions si
les partitions des deux découpages ont le même nombre d'ensembles de la même
taille pour chaque couleur, puis nous alignons les deux découpages,
et  procédons exactement comme pour le
problème de l'automorphisme.



\subsubsection{Le cas du schéma de Johnson}
Soit donné un schéma de Johnson sur un ensemble de sommets $\Gamma$,
ou plutôt deux schémas de Johnson $\mathscr{J}(m_i,k_i)$,
$2\leq k_i\leq m_i/2$,
sur des ensembles de sommets $\Gamma_1$,
$\Gamma_2$ de la même taille. Nous avons vu au \S \ref{subs:idgroup}
comment identifier $\Gamma_i$ (là, $\Omega$) explicitement avec
les ensembles de taille $k_i$ d'un ensemble $\Lambda_i$ (là, $\Gamma$)
de taille $m_i$. Si
$k_1\ne k_2$ et $m_1\ne m_2$, nos structures ne sont pas isomorphes.
Si $k_1=k_2$ et $m_1=m_2$, nous établissons une bijection entre
$\Lambda_1$ et $\Lambda_2$ et alignons les deux structures. Nous avons
réduit notre problème à un problème avec $m'\ll \sqrt{m}$ à la place
de $m$.

La situation nous est donc même plus favorable que dans le cas du découpage.
\`A nouveau, nous laissons la comptabilité au lecteur.

\begin{center}
  * * *
\end{center}

Une petite confession: le cas de $G$ primitif, que nous venons de finir de
traiter, pourrait être traité exactement comme le cas de $G$ imprimitif, que nous examinerons maintenant. La motivation du traitement séparé pour
$G$ primitif est pédagogique. Aucune peine n'est perdue, puisque toutes les
techniques que nous avons étudiées nous seront essentielles dans le cas
imprimitif.

\section{Le cas imprimitif}\label{sec:casimp}

Nous avons une application surjective explicite
\[\phi:G\to \Alt(\Gamma),\]
où $G<\Sym(\Omega)$ est un groupe de permutation, $|\Gamma|=m$,
$|\Omega|=n$. Nous pouvons supposer que $|\Gamma|\geq C \log n$, $C$ arbitraire. L'application $\phi$ se factorise comme suit
\[G\to G/N\to \Alt(\Gamma),\]
où $N$ est le stabilisateur d'un système de blocs, et $G/N\to \Alt(\Gamma)$
est un isomorphisme.

Nous devons déterminer $\Iso_G(\mathbf{x},\mathbf{y})$, où
$\mathbf{x}$, $\mathbf{y}$ sont des chaînes. Nous avons déjà
résolu le cas $N = \{e\}$.

Nous attaquerons le problème de façon locale: pour $T\subset \Gamma$,
nous arriverons à obtenir un certificat, soit du fait que 
$\phi(\Aut_{G_T}(\mathbf{x}))|_T$ contient $\Alt(T)$ ({\og certificat
de plénitude\fg}), soit du contraire. (Ici $G_T$ désigne
le groupe $\{g\in G: T^{\phi(g)} = T\}$.)
Nous calculerons tous ces certificats pour
$T$ d'une taille $k$ modérée. Si le nombre de certificats de plénitude
est très grand, nous aurons prouvé que $\phi(\Aut_G(\mathbf{x}))$ contient
un grand groupe alternant; ce qui restera à faire sera une version de
la procédure du \S \ref{sec:grasym} ({\og pull-back\fg}).

Dans le cas contraire, les certificats formeront une structure
$k$-aire dont la symétrie est bornée. Nous pourrons donc appliquer le Lemme
des designs, suivi de Coupe-ou-Johnson, comme avant. Il y a aussi quelques
autres cas particuliers, mais ils nous amènent à des $\alpha$-découpages,
$\alpha<1$, ce qui est aussi bien.
\subsection{Les certificats locaux}\label{sec:certloc}
\subsubsection{Certificats d'automorphismes}\label{subs:certaut}
Un certificat local\footnote{
Ou {\og local-global\fg}, dans la nomenclature de Babai. {\og Global\fg}
fait référence à 
$\Aut_{G_T}(\mathbf{x}) <\Sym(\Omega)$.}
  pour $T\subset \Gamma$ est
\begin{itemize}
\item  soit une paire $(\text{{\og pas plein\fg}},W,M(T))$, où
  $W\subset \Omega$, $M(T)<\Sym(T)$,
  $M(T)\neq \Alt(T)$ (donc {\og pas plein\fg}) et
  $\phi\left(\Aut_{G_T}^W(\mathbf{x})\right)|_T
  < M(T)$,
\item  soit une paire $(\text{{\og plein\fg}},K(T))$, où $K(T)<\Aut_{G_T}(\mathbf{x})$,
  et $\phi(K(T))|_T = \Alt(T)$.
\end{itemize}
Le certificat local dépend de $\mathbf{x}$ de façon canonique.
Il est clair qu'un certificat plein, voire pas plein, 
garantit que $\phi(\Aut_{G_T}(\mathbf{x}))|_T$ est $\Alt(T)$, voire
ne l'est pas.

Si $T$ est donné en tant que tuple ordonné, son certificat dépend de l'ordre
de $T$ seulement dans le sens de ne pas en dépendre: le même groupe
$\{(2 3), e\} < \Sym(\{1,2,3\})$ (disons) a une apparence différente si nous
le regardons du point de vue de l'ordre $(1,2,3)$ ou de l'ordre $(2,1,3)$.

Nous construisons le certificat par une procédure itérative.
Au début de chaque pas, $W\subset \Omega$ et
$A(W)$ est le groupe $\Aut_{G_T}^W(\mathbf{x})$; la fenêtre $W$ sera invariante
sous $A(W)$. Au tout début de la procédure, $W=\emptyset$ et
$A(W) = G_T$.
(Nous pouvons calculer
$G_T$ comme dans l'exercice \ref{ex:fhl}\ref{it:richt}
en temps $|\Omega|^{O(k)}$, où $k = |T|$.)
\`A chaque pas, nous ajoutons à $W$
tous les éléments {\em atteints} par $A(W)$
(voir \S \ref{subs:stab}), puis nous mettons $A(W)$ à jour,
selon le nouveau $W$. Nous nous arrêtons si
$\phi(A(W))|_T \ne \Alt(T)$ (non-plénitude) ou si $W$
ne croît plus, ce qui veut dire qu'aucun élément de $\Omega\setminus
W$ n'est atteint par $A(W)$.

Il est clair qu'il y aura $\leq |\Omega|$ itérations. \`A la fin, dans
le cas de non-plénitude, nous retournons $(\text{{\og pas plein\fg}},
W,\phi(A(W)))$; dans le cas de plénitude, nous retournons
$\left(\text{{\og plein\fg}},A(W)_{(\Omega\setminus W)}\right)$.
Il est clair
  que le stabilisateur des points $A(W)_{(\Omega\setminus W)}$
  est contenu non pas seulement dans $\Aut_{G_T}^W(\mathbf{x})$, mais aussi dans $\Aut_{G_T}(\mathbf{x})$,
  puisqu'il fixe tous les points de $\Omega\setminus W$.
  Nous savons que $\phi\left(A(W)_{(\Omega\setminus W)}\right) = \Alt(T)$
  par la Proposition \ref{prop:atinl}\ref{it:unaffstab}, sous la condition
  que $|T|\geq \max(8,2+\log_2 |\Omega|)$. 
  
  Vérifier si $\phi(A(W))|_T = \Alt(T)$ est facile: nous n'avons qu'à vérifier,
  en utilisant Schreier-Sims, si deux générateurs arbitraires de $\Alt(T)$ sont
  en $\phi(A(W))|_T$. De la même façon, il est simple de déterminer
  quels éléments sont atteints par $A(W)$: nous calculons $A(W)_x$
  pour chaque $x\in \Omega$ (par Schreier-Sims) et, toujours par Schreier-Sims,
  vérifions si $\phi(A(W)_x)|_T = \Alt(T)$. Ceci prend du temps polynomial en
  $|\Omega|$.

  Il reste à voir comment mettre à jour $A(W)$, étant donné
  $A\left(W^-\right)$, où nous écrivons $W^-$ pour l'ancienne valeur de $W$.
  Tout élément de $A(W)$ est dans $A(W^-)$, et
  donc $A\left(W\right) = \Aut_{A\left(W^-\right)}^{W}(\mathbf{x})$.
  Comme dans l'équation (\ref{eq:rulu}),
  \begin{equation}\label{eq:ruyur}
    \Aut_{A\left(W^-\right)}^{W}(\mathbf{x}) =
  \bigcup_\sigma \Aut^{W}_{N\sigma}(\mathbf{x}) =
  \bigcup_\sigma \Iso_N^{W}\left(\mathbf{x},\mathbf{x}^{\sigma^{-1}}\right),
  \end{equation}
  où $N$ est le noyau de $\phi|_{A(W_-)}$ et $\sigma$ parcourt des
  représentants des $k!/2$ classes de $N$ en $A(W)$.
  Nous pouvons trouver rapidement
  un $\sigma\in A(W^-) \cap \phi^{-1}(\{\tau\})$ pour tout $\tau\in
  \Sym(\Gamma)$, par Schreier-Sims.

  La Proposition \ref{prop:atinl}\ref{it:afforb}
  nous donne que toute orbite de $N$ contenue en $W$ (l'ensemble
  d'éléments atteints par $A(W^-)$) est de longueur $\leq |W|/k \leq
  |\Omega|/k$. En conséquence, par la règle (\ref{eq:udu3}), mettre
  $A(W)$ à jour se réduit à $|\Omega|\cdot (k!/2)$ problèmes
  de détermination de $\Iso$ pour des chaînes de longueur $\leq |\Omega|/k$.

  Comme le nombre d'itérations est $\leq |\Omega|$, la procédure fait
  appel à Isomorphisme-de-Chaînes $\leq |\Omega|^2 \cdot (k!/2)$ fois
  pour des chaînes de longueur $\leq |\Omega|/k$. Ceci -- comme
  la routine qui prenait $|\Omega|^{O(k)}$ de temps -- est acceptable pour
  $k\ll (\log |\Omega|)^\kappa$. Nous choisirons $\kappa=1$.
  
  \subsubsection{Comparaisons de certificats}
  Une légère modification de la procédure ci-dessus nous permet
  d'élucider la relation entre deux certificats locaux pour deux
  chaînes. Soient
  $\mathbf{x}, \mathbf{x}':\Omega\to \Sigma$, $T,T'\subset\Sigma$,
  $|T|=|T'|=k$. Pour $S\supset T$, soit $\mathbf{x}^S$ la chaîne
  \[\mathbf{x}^S(i) = \begin{cases} \mathbf{x}(i) & \text{si $i\in S$,}\\
    \text{glauque} & \text{si $i\notin S$}\end{cases}\]
  où $\text{glauque} \notin \Sigma$. Nous voulons calculer
  \begin{equation}\label{eq:headtail}
    \Iso_{G_T \cdot \tau_{T,T'}}\left(\mathbf{x}^W,\mathbf{x}^{W'}\right),\end{equation}
    où $G_T \cdot \tau_{T,T'}$ est la classe des éléments de $G$ qui envoient
    l'ensemble $T$ à $T'$, et $W'$ est la valeur de $W$ retournée
    quand la donnée est $T'$ à la place de $T$.

    Pour déterminer (\ref{eq:headtail}), nous suivons la procédure
    (\S \ref{subs:certaut}), modifiée de la façon suivante:
    nous mettrons à jour,
    dans chaque itération, non pas seulement $A(W)$, mais aussi la classe
    $A(W) \tau$ d'isomorphismes en $G_T \cdot \tau_{T,T'}$ de $\mathbf{x}^W$
    à $(\mathbf{x}')^{W'}$. Voilà comment le faire, de façon analogue à
    (\ref{eq:ruyur}):
    \begin{equation}\label{eq:rurik}
    \bigcup_\sigma \Iso_{N\sigma}\left(\mathbf{x}^{W},
    \left(\mathbf{x}'\right)^{W'}\right) =
    \bigcup_\sigma \Iso_N\left(\mathbf{x}^{W},
    \left(\left(\mathbf{x}'\right)^{W'}\right)^{\sigma^{-1}}\right),\end{equation}
      où $N$ est le noyau de $\phi|_{A(W_-)}$ et $\sigma$ parcourt des
      représentants des $k!/2$ classes de $N$ contenues en $A(W^-) \tau^-$.
      Comme $W$ est stabilisé par $A(W_-)$ (et donc par $N$),
      le fait que $\sigma$ envoie $W$ sur $W'$ ou non dépend seulement
      de la classe de $N$ à laquelle $\sigma$ appartient. (La classe
      $\Iso$ dans la dernière expression de (\ref{eq:rurik}) est vide
      si $W^\sigma \ne W'$.)

      Comme avant, toute orbite de $N$ contenue en $W$ est de longueur
      $\leq |W|/k$, et le problème se réduit à $|\Omega|\cdot (k!/2)$ appels
      par itération
      à Isomorphisme-de-Chaînes pour des chaînes de longueur $\leq |W|/k\leq
      |\Omega|/k$.

      Par ailleurs, si $T$ et $T'$ nous sont données comme tuples ordonnés
      $(T)$, $(T')$, il est facile de déterminer
      \begin{equation}\label{eq:isoma}I(\mathbf{x},\mathbf{x}',T,T') =
        \Iso_{G_{(T)} \cdot \tau_{(T),(T')}}\left(\mathbf{x}^W,
      \left(\mathbf{x}'\right)^{W'}\right),\end{equation}
      où $G_{(T)} \cdot \tau_{(T),(T')}$ est la classe des éléments de $G$ qui
      envoient le tuple ordonné $(T)$ à $(T')$. En effet, nous n'avons qu'à
      déterminer (\ref{eq:headtail}), puis utiliser Schreier-Sims pour déceler
      les éléments de (\ref{eq:headtail}) qui envoient $(T)$ à $(T')$ dans
      le bon ordre.

\subsection{L'agrégation des certificats}\label{subs:agrecert}

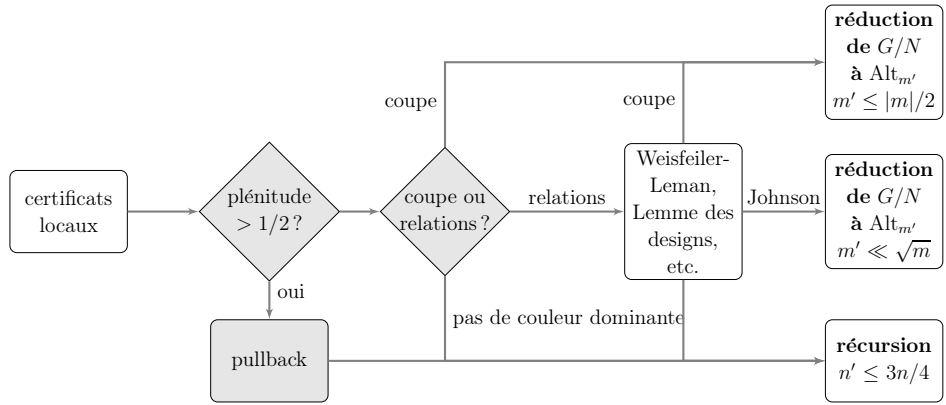
\begin{center}
\scalebox{0.66}{
  \begin{tikzpicture}[scale=2, node distance = 2cm, auto]
    \node [noblock] (cert) {certificats locaux};
    \node [decision , right of = cert, node distance = 4cm] (certsym)
       {plénitude $>1/2$?};
    \node [block, below of = certsym, node distance = 3cm] (pull)
             {pullback};
    \node [decision, right of = certsym, node distance = 3.5cm] (couprel)
           {coupe ou relations?};
    \node [noblock, right of = couprel, node distance = 4.75cm] (designs)
          {Weisfeiler-Leman, Lemme des designs, etc.};
    \node [noblock, right of = designs, node distance = 4cm] (redsqrt) {\bf réduction de $G/N$ à $\Alt_{m'}$\\ $m'\ll \sqrt{m}$};
              \node [noblock, below of = redsqrt, node distance = 3cm] (luks) {\bf récursion
      $n'\leq 3n/4$};
   \node [noblock, above of = redsqrt, node distance = 3cm] (redhalf)
           {\bf réduction de $G/N$ à $\Alt_{m'}$\\ $m'\leq |m|/2$};
    \path[line] (couprel) -- node [color=black] {relations} (designs);
    \path[line] (cert) -- (certsym);
    \path [line] (certsym) -- node [near start, color=black] {oui} (pull);
    \path[line] (pull) -- (luks);
    \path [line] (certsym) -- (couprel);
    \path [line] (couprel) |- node [near start, color=black] {coupe} (redhalf);
   \path [line] (designs) |- node [near start, color=black] {coupe} (redhalf);
   \path [line] (couprel) |- node [near start, color=black]
         {pas de couleur dominante} (luks);
         \path [line] (designs) -- node [color=black] {Johnson} (redsqrt);
            \path [line] (designs) |- (luks);

\end{tikzpicture}
  }
\end{center}

En suivant la procédure du
\S \ref{subs:certaut} pour une chaîne $\mathbf{x}$, nous trouvons des certificats locaux pour chaque
$T\subset \Gamma$ de taille $k$, où $k$ est une constante
$\sim C \log |\Omega|$ ($C>1/\log 2$) et
$k<|\Gamma|/10$. Soit $F<\Aut_G(\mathbf{x})$
le groupe engendré par les certificats pleins $K(T)$. 
Soit $S\subset \Gamma$ le support de $\phi(F)$, c'est-à-dire
l'ensemble des éléments de $\Gamma$ qui ne sont
pas fixés par tout \'el\'ement de $\phi(F)$.

Notre objectif est de déterminer les isomorphismes $\Iso_G(\mathbf{x},
\mathbf{x}')$ de $\mathbf{x}$ à une autre chaîne~$\mathbf{x}'$.
Puisque les certificats sont canoniques, l'assignation de $F$ et $S$ à
une chaîne l'est aussi. Donc, si nous arrivons à deux cas différents
ci-dessous en suivant la procédure pour $\mathbf{x}$ et pour $\mathbf{x}'$, les deux chaînes
ne sont pas isomorphes.
  
{\noindent Cas 1: $|S|\geq |\Gamma|/2$, mais aucune orbite de $\phi(F)$
  n'est de longueur $>|\Gamma|/2$.}

Alors, nous colorions chaque élément de $\Gamma$ par la longueur de l'orbite
qui le contient. Ceci est un coloriage canonique.
Soit aucune classe de couleurs n'est de taille $>|\Gamma|/2$,
soit une classe de couleurs de taille $>|\Gamma|/2$ est découpée en $\geq 2$ ensembles de la même taille $\geq 2$. Dans un cas comme dans l'autre,
nous passons à une réduction/récursion.

{\noindent Cas 2: $|S|\geq |\Gamma|/2$ et une orbite $\Phi$ de $\phi(F)$ est
  de longueur $>|\Gamma|/2$.}

{\em Cas 2a: $\Alt(\Phi) < \phi(F)|_\Phi$.} Nous sommes dans
le cas de grande symétrie.
Nous procédons comme au \S \ref{sec:grasym}, jusqu'au point où nous devons
déterminer $\Iso_H(\mathbf{x},\mathbf{y})$ (où
$\mathbf{y}$ est $\left(\mathbf{x}'\right)^{\sigma'}$, $\sigma'\in G$,
et $H = \phi^{-1}\left(\Alt(\Gamma)_\Phi\right)$). Définissons
$K = \phi^{-1}\left(\Alt(\Gamma)_{(\Phi)}\right)$, et soient
$\sigma_1,\sigma_2\in G$ des préimages (arbitraires) sous $\phi$
de deux générateurs de $\Alt(\Phi)<\Alt(G)$, trouvées par Schreier-Sims.
Nous savons que les classes
$\Aut_{K \sigma_i}(\mathbf{x})$, $i=1,2$, sont non vides, puisque
$\Alt(\Phi)<\phi(F)|_\Phi$.
Comme $K$ n'a pas d'orbites de longueur $>|\Omega|/2$, nous pouvons
déterminer ces deux classes par des appels à 
Isomorphisme-de-Chaînes pour des chaînes de longueur $\leq |\Omega|/2$
et longueur totale $\leq 2 |\Omega|$. Elles engendrent
$\Aut_H(\mathbf{x})$.

Encore par le fait que $\Alt(\Phi)<\phi(F)|_\Phi$,
la classe $\Iso_H(\mathbf{x},\mathbf{y})$
sera non vide ssi $\Iso_K(\mathbf{x},\mathbf{y})$ est non vide.
Nous pouvons déterminer
cette dernière classe par des appels à Isomorphisme-de-Chaînes comme ci-dessus,
puisque $K$ n'a pas d'orbites de longueur $>|\Omega|/2$. Si elle est non vide,
nous obtenons la réponse
\[\Iso_H(\mathbf{x},\mathbf{y}) = 
\Aut_H(\mathbf{x}) \Iso_K(\mathbf{x},\mathbf{y}).\]


{\em Cas 2b: $\Alt(\Phi) \nless \phi(F)|_\Phi$.} Soit $d\geq 1$ l'entier maximal
avec la propriété que $\phi(F)|_\Phi$ est $d$-transitif, c'est-à-dire,
$\phi(F)|_\Phi$ agit transitivement sur l'ensemble des $d$-tuples d'éléments
distincts de $\Phi$. Par CGFS, $d\leq 5$; si nous ne voulons pas utiliser
CGFS, nous avons la borne classique $d\ll \log |\Gamma|$.

Choisissons $x_1,\dotsc,x_{d-1}\in \Phi$ arbitrairement.
Le reste de notre traitement de ce cas sera donc seulement canonique  en relation à
\[G_{(x_1,\dotsc,x_{d-1})} = \{g\in G: x_i^{\phi(g)} = x_i \; \forall 1\leq i\leq d-1\},\]
ce qui, comme
nous le savons, n'est pas un problème; voir
le début du \S \ref{subs:cacoudo}. 

La restriction du groupe $\phi(F)_{(x_1,\dotsc,x_{d-1})}$ à $\Phi' =
\Phi\setminus \{x_1,\dotsc,x_{d-1}\}$
est transitive sur $\Phi'$, mais elle n'est pas doublement transitive.
Donc, la configuration cohérente schurienne qui lui correspond n'est pas une
clique. Nous livrons cette configuration à Coupe-ou-Johnson (\S \ref{sec:coujoh}), tout comme à la fin du \S \ref{sec:coujoh}.

Pour comparer les configurations qui correspondent à deux chaînes
$\mathbf{x}$, $\mathbf{x}'$, nous alignons leurs classes $\Phi$ d'abord.
(Si elles ne sont pas de la même taille, ou si une chaîne nous donne le
cas 2a et l'autre pas, les chaînes ne sont pas isomorphes.) Les isomorphismes
seront donc contenus dans le stabilisateur $H<G$ de l'ensemble $\Phi$ (facile à
déterminer, comme vers
la fin du \S \ref{sec:grasym}, puisque $\phi$ est surjective). Nous pouvons
remplacer $\phi$ par l'application $g\mapsto \phi(g)|_\Phi$ de $H$ à $\Alt(\Phi)$.
Puis nous construisons les configurations comme ci-dessus, et comparons ce
que Coupe-ou-Johnson nous donne.

Tout à la fin, nous nous occupons du complément de $C$. Il s'agit, comme d'habitude,
d'appels à Isomorphisme-de-Chaînes pour des chaînes de longueur $\leq |\Omega|/2$
et longueur totale $<|\Omega|$.

{\noindent Cas 3: $|S|< |\Gamma|/2$.}
Nous commençons en alignant les supports $S$ pour les
chaînes $\mathbf{x}$, $\mathbf{x}'$,
et en remplaçant $\phi$ par $g\mapsto \phi(g)|_{\Gamma \setminus S}$,
tout comme dans le cas 2(b).

Nous allons définir une relation $k$-aire avec très peu
de jumeaux, pour la donner après au Lemme des designs.

Regardons la catégorie de toutes les chaînes $\Omega \to \Sigma$,
où $\Omega$ et $\Sigma$ sont fixes, une action de $G$ sur $\Omega$ est
donnée, et $\phi:G\to \Gamma$ est aussi donnée. Nous la regardons depuis
longtemps, puisque nous devons comparer les couleurs sur des configurations
induites par des chaînes différentes pour décider si ces dernières sont
isomorphes.

Cette fois-ci, nous définirons des couleurs par des classes d'équivalence:
deux paires $(\mathbf{x},(T))$, $(\mathbf{x}',(T'))$
($T,T'\subset \Gamma\setminus S$, $|T|=|T'|=k$) sont équivalentes si l'ensemble
des isomorphismes en (\ref{eq:isoma}) est non vide.
Nous colorions $(T)$ -- dans le coloriage de $(\Gamma\setminus S)^k$
correspondant à
$\mathbf{x}$ -- par la classe d'équivalence de $(\mathbf{x},(T))$.
Ici, $(T)$ est un
$k$-tuple ordonné sans répétitions; si $(T)$ a des répétitions,
elle est coloriée en gris.

Pour $\mathbf{x}$ donné, aucune classe de jumeaux en $\Gamma$ ne peut avoir
$\geq k$ éléments: s'il existait un tel ensemble avec $\geq k$
éléments, il contiendrait un ensemble $T$ avec $k$ éléments,
et tous les ordres $(T)$ de $T$ auraient la même couleur. Ceci voudrait
dire que l'ensemble des isomorphismes en (\ref{eq:isoma}) serait non vide
pour n'importe quels ordres $(T)$, $(T')$ de $T$. En conséquence,
$\Aut_{G_T}(\mathbf{x}^W)$ contiendrait des éléments donnant toutes les
permutations possibles de $T$. Ceci nous donnerait une contradiction, puisque
$T$, étant contenu en $\Gamma\setminus S$, n'est pas plein.

Alors, pourvu que $k\leq |\Gamma|/4$, nous avons un coloriage de
$(\Gamma\setminus S)^k$ sans aucune classe
de jumeaux avec $\geq |\Gamma\setminus S|$ éléments.
Nous pourrons donc appliquer le Lemme des designs, après une application
de raffinements habituels $F_2$, $F_3$ (Weisfeiler-Leman).

Mais -- pouvons-nous calculer ces coloriages? Les classes d'équivalence sont
énormes. Par contre, il n'y a aucun besoin de les calculer. Tout ce dont nous aurons
besoin, pour comparer des structures qui viennent
de chaînes $\mathbf{x}$, $\mathbf{y}$, sera d'être capables de
comparer deux tuples $(T)$ (sur la configuration donnée par $\mathbf{x}$
ou $\mathbf{y}$) et
$(T')$ (sur la configuration donnée par $\mathbf{x}'=\mathbf{x}$ ou
$\mathbf{x}'=\mathbf{y}$) et dire si elles sont de la même couleur.
En d'autres termes, nous devrons calculer -- au tout début de la procédure,
pour toute paire $((T),(T'))$, $|T|=|T'|=k$, et pour les paires de
chaînes $(\mathbf{x},\mathbf{x})$, $(\mathbf{x},\mathbf{y})$, $(\mathbf{y},\mathbf{y})$ --
l'ensemble d'isomorphismes en
(\ref{eq:isoma}), ce que nous savons déjà faire. Les couleurs
sont donc, dans la pratique, des entrées dans un index que nous
enrichissons et auquel nous faisons référence durant nos procédures.

Nous invoquons donc le Lemme des Designs, suivi par Coupe-ou-Johnson, et le
reste de la procédure.
\begin{center}
  {--- \textsc{Fine dell'opera} ---}\\ \vskip 7pt
\end{center}

Le lecteur peut vérifier que les informations précisées jusqu'à ici
(temps pris par des
procédures, type de récursion) 
sont assez pour donner
une borne du type $\exp(O(\log |\Omega|)^c)$ pour le temps
de l'algorithme qui résout le problème
de l'isomorphisme de chaînes. Ceci donne une borne
$\exp(O(\log n)^c)$ pour le problème de l'isomorphisme de graphes avec $n$ sommets. Avec un peu plus de travail, il devient clair que, dans un cas comme dans l'autre, $c=3$. Nous donnons
les détails dans l'appendice. L'exposant $c=3$ est plus petit que celui
d'origine; il est devenu possible grâce à quelques améliorations et simplifications que j'ai été capable d'apporter.


Remerciements .--- Je remercie vivement L. Babai, J. Bajpai,
L. Bartholdi, D. Dona, E. Kowalski,
W.~Kantor, G. Puccini,
L. Pyber, A. Rimbaud et C. Roney-Dougal pour des corrections et suggestions.
En particulier, L. Babai a répondu à beaucoup de mes questions, et m'a aussi
fourni des versions corrigées ou améliorées de plusieurs sections de \cite{Ba}.
En particulier, les \S \ref{subs:secf1}--\ref{subs:concoh} et \S \ref{sec:design}
sont basés sur ces nouvelles versions. Je voudrais aussi remercier
V.~Ladret et V. Le Dret pour un grand nombre de corrections d'ordre typographique et linguistique.



\appendix

\section{Analyse du temps d'exécution}

\subsection{Quelques précisions sur la procédure principale}

\`A tout moment donné, nous travaillons avec un groupe transitif
$G<\Sym(\Omega)$ qui agit sur un système de blocs $\mathscr{B} = \{B_i\}$,
$\Omega = \bigcup_i B_i$, $B_i$ disjoints;
nous notons $N$ le noyau de l'action sur $\mathscr{B}$.
\`A vrai dire, nous aurons
toute une tour de systèmes de blocs $\mathscr{B}_1,\mathscr{B}_2,\dots,
\mathscr{B}_k$,
où $B_i$ est un raffinement de $B_{i+1}$; $\mathscr{B}$
signifiera $\mathscr{B}_k$, le système le moins fin.
Au début, il n'y a qu'un système, $\mathscr{B}_1$, dont les blocs $B_i$ sont
tous de taille $1$, et dont
le noyau $N$ est trivial. 



Nous voudrions que l'action de $G$ sur $\mathscr{B}$ soit primitive.
Donc, si elle ne l'est pas, nous ajoutons à la tour un système minimal
$\mathscr{B}_{k+1}$ tel que
$\mathscr{B}_k$ soit un raffinement de $\mathscr{B}_{k+1}$. Nous
redéfinissons $\mathscr{B} = \mathscr{B}_{k+1}$;
$N$ sera le noyau du nouveau $\mathscr{B}$.

Si $G/N$ est petit ($\leq b^{O(\log b)}$, où $b=|\mathscr{B}|$;
cas (a) du Th\'eor\`eme~\ref{thm:cam} (Cameron)), nous
réduisons notre problème à plusieurs instances du problème avec $N$ à la place
de $G$. Chacune de ces instances se décompose en plusieurs instances -- une pour
chaque orbite de $N$. Chaque orbite $\Omega'$ de $N$ est contenue dans un
bloc de $\mathscr{B}$. Les intersections de $\Omega'$ avec les blocs de
$\mathscr{B}_1,\mathscr{B}_2,\dotsc$ nous donnent une tour
de systèmes de blocs pour $N|_{\Omega'}$. 

Si nous sommes dans le cas (b) du Th\'eor\`eme~\ref{thm:cam}, nous passons à
$\leq b$ instances du problème avec $M\triangleleft G$ (où $\lbrack G:M\rbrack
\leq b$) à la place de $G$.
Nous passons à un nouveau système\footnote{Ce système peut être égal à $\mathscr{B}$ seulement si
  $M=G$; voir la deuxième note de pied de page dans l'énoncé du Th\'eor\`eme~3.1.
Dans ce cas-là, le passage de $G$ à $M$ est bien sûr gratuit.}
$\mathscr{B}'$
de $m' = \binom{m}{k} \leq b$ blocs, et l'ajoutons à la tour comme
son nouveau dernier niveau.
Nous notons $N'$ le noyau de l'action de $M$ sur $\mathscr{B}'$.
Alors, $M/N' = \Alt_m^{(k)}$. Nous remplaçons $G$ par $M$ et
redéfinissons $\mathscr{B} = \mathscr{B}'$, $N=N'$.

Donc, nous avons un isomorphisme de $G/N$ à $\Alt_m$.
Nous sommes dans le cas principal que Babai attaque.
Ses méthodes amènent à une réduction de $\Alt_m$, soit à un groupe
intransitif sans grandes orbites, soit à un produit $\Alt_{s_1}\wr \Alt_{s_2}$,
$s_1,s_2>1$, $s_1 s_2\leq m$, 
soit à un groupe $\Alt_{m'}$, $m'\ll \sqrt{m}$. (Nous simplifions
quelque peu. Nous pourrions avoir, disons, un produit $\Alt_{s_1}\wr \Alt_{s_2}$,
agissant sur une orbite de grande
taille $s_1 s_2\leq m$, et d'autres groupes sur des
petites orbites, ou plusieurs produits agissant sur des petites orbites.)

Dans le cas intransitif
sans grandes orbites, nous procédons comme dans la preuve de Luks.
(La procédure aura été plus coûteuse que dans Luks, mais grâce au
manque de grandes orbites, le gain dans la récursion est aussi plus grand.)
Dans le cas de $\Alt_{m'}$, $m'\ll \sqrt{m}$, nous itérons la procédure.
Dans le cas de $\Alt_{s_1}\wr \Alt_{s_2}$ -- qui correspond à un découpage dans
des ensembles de taille $r$ de la même couleur
-- nous avons une action primitive de
$\Alt_s$ sur un système de $s$ blocs de taille $r$.
Nous passons, alors, à cette action
et à ces blocs, sans oublier les blocs $\mathscr{B}'$, auxquels nous retournons
plus tard, après avoir fini de travailler sur $\Alt_s$.

Il est clair que ce type de procédure réduit complètement $\Alt_k$ 
en un nombre d'it\'erations qui n'est pas sup\'erieur \`a $\log_2 m$. 

\subsection{Récursion et temps}

Examinons le temps total d'exécution de l'algorithme qui trouve les
isomorphismes entre deux chaînes.
Les pas individuels sont peu onéreux; aucun ne précise plus de $n^{O(\log n)}$
de temps. Notre attention doit se porter avant tout sur la récursion.

Dans la procédure générale, une récursion 
est toujours d'une descente, soit vers des chaînes plus
courtes, soit vers un groupe plus petit, ou au moins coupé dans des tranches
plus fines par une tour de systèmes de blocs ayant plus de niveaux. Dans le
premier type de descente, le groupe
reste le même ou, plutôt, est remplacé par une restriction de lui-même.
Dans le deuxième cas, la longueur des chaînes reste la même. (Nous pouvons
aussi avoir un mélange des deux cas -- tant mieux:
le groupe devient plus petit et les chaînes se raccourcissent aussi.)

La descente la moins coûteuse, et moins avantageuse, est celle du
cas intransitif de la procédure de Luks. Il pourrait arriver que $G$
ait deux orbites sur $\Omega$ ($|\Omega|=n$), une de longueur
$n-1$ et une de longueur
$1$. Ceci serait même compatible  avec une borne polynomiale sur le temps,
 pourvu que le temps pris avant la descente soit lui-même polynomial:
$n^{c+1} \leq (n-1)^{c+1} + 1^{c+1} + n^c$ pour $c\geq 1$.

D'autres types de descente sont plus coûteux, mais aussi plus avantageux:
nous descendons à des chaînes de longueur $\leq n/2$ (ou $\leq 2 n/3$),
ou de $\Alt_m$
à $\Alt_{s_1}\wr \Alt_{s_2}$, $s_1 s_2 \leq m$, $s_1,s_2\leq m/2$, par exemple. Il est
clair qu'il est impossible de descendre plus qu'un nombre logarithmique
de fois de cette façon.

Il est crucial de ne pas oublier qu'un coût (considérable) peut être caché
dans une perte de canonicité. 
Si nos choix ne sont canoniques qu'en relation
à un sous-groupe $H$ de notre groupe $G$, le coût de leur application sera
multiplié par $\lbrack G:H\rbrack$. (Voir \S \ref{subs:cacoudo}.)

\begin{center}
  * * *
\end{center}

Considérons, alors, le coût de chaque procédure.
Le cas intransitif de Luks est, comme nous l'avons déjà vu, compatible même
avec une borne polynomiale. Concentrons-nous alors sur le cas où $G$ agit
de façon primitive sur un système de blocs; soit $N$ le noyau.

Si nous sommes dans le
cas (\ref{it:rila}) du Th\'eor\`eme~\ref{thm:cam}, ou dans le cas (\ref{it:rilb}), mais avec $m\leq C \log n$,
nous faisons appel \`a $(m')^{O(\log n)}$ instances de la procédure principale
pour des chaînes de longueur
$n/m'$ (où $m'\geq m$). Ceci est consistant avec une borne totale du type
$\exp(O((\log n)^c))$, $c\geq 2$. Nous pouvons, donc, nous concentrer sur le cas
où il existe un isomorphisme $\phi:G/N\to \Alt(\Gamma)$, $|\Gamma| = m
> C\log n$. (La procédure du \S \ref{subs:idgroup} rend cet isomorphisme explicite.)

Le premier pas à considérer est la création de certificats
locaux, avec, comme objectif, la création d'une relation $k$-aire sur 
$\Gamma$. (Si $G$ est primitif, créer une telle relation est trivial; voir le
début du \S \ref{subs:chasche}.)
Il y a $n^k$ certificats locaux, où $k = 2 \log n$
(disons); nous devons les calculer et aussi
comparer toute paire de certificats.
Déjà le premier pas du calcul d'un certificat, à savoir
le calcul de $G_T$, prend un temps $n^{O(k)}$ (plus précisément,
$O((n/k)^{O(k)})$). D'autres calculs prennent moins de temps. 
L'usage de la récursion, par contre, est relativement lourd:
nous faisons appel \`a  la procédure principale $\leq n^2 \cdot k!$ fois
pour des chaînes de longueur $\leq n/k$. Ceci se passe pour chaque
ensemble $T$ de taille~$k$, c'est-à-dire $\leq n^k/k!$ fois.
La procédure pour comparer des paires de certificats est analogue.

Nous faisons donc appel à la procédure principale $O(n^{2 k+1})$ fois
pour des chaînes de longueur $\leq n/k$.
Dans chacun de ces appels, notre tour de stabilisateurs est héritée:
notre groupe est un groupe transitif, égal à la restriction de $N^-$
à une de ses orbites, où $N^-$ (noté $N$ au \S \ref{sec:casimp})
est un sous-groupe d'un sous-groupe $A(W^-)$ de $G$.

Pour deux systèmes de blocs consécutifs $\mathscr{B}_i$,
$\mathscr{B}_{i+1}$, notons $r_i$ le nombre de blocs de $\mathscr{B}_i$
dans chaque bloc de $\mathscr{B}_{i+1}$. Il est clair que
ce nombre n'augmente pas quand nous passons à la restriction d'un sous-groupe
de $G$ (par exemple, $N^-$) à une de ses orbites.

Examinons maintenant l'agrégation des certificats locaux (\S \ref{subs:agrecert}). Il y
a trois cas. Dans le premier, le temps de calcul additionnel
est à peu près trivial, et nous obtenons une réduction, soit à un groupe
intransitif sans grandes orbites, soit à 
un produit $\Alt_{s_1}\wr \Alt_{s_2}$ sur une grande orbite et éventuellement
d'autres groupes sur des orbites plus petites.

Ici, déjà, l'analyse devient délicate. Nous devons prendre en considération
non seulement la taille du domaine mais aussi le groupe qui agit sur lui.
Plus précisément, nous devons borner le nombre de fois que notre tour
$\mathscr{B}_1, \mathscr{B}_2,\dotsc,\mathscr{B}_k$ pourrait être
raffiné ou raccourci encore. Ceci sera mesuré par
\[\rho = \sum_{1\leq i\leq k-1} (2 \lfloor \log_2 r_i\rfloor - 1),\]
où nous supposons que nous avons enlevé des systèmes répétés
de la tour (donc $r_i>1$).

Notons $F(n,r)$ le temps d'exécution de la procédure principale
pour des chaînes de longueur $n$ et pour une tour de systèmes de blocs pour $G$
telle que le paramètre $\rho$ est $\leq r$.
Une réduction de $G/N$ fait décroître $r$ par au moins $1$; un coloriage
sans aucune grande classe de couleurs assure une descente vers des chaînes
de longueur $\leq n/2$. Nous devrons aussi inclure un facteur de $\log n^{2k}$, prenant en considération le temps requis pour accéder à nos comparaisons
de paires de certificats locaux\footnote{Faire ce type de comparaisons
  à l'avance nous aide, mais ne pas les faire à l'avance ne changerait
  pas l'ordre du temps utilisé, asymptotiquement.}.
Donc, dans le cas que nous examinons, $F(n,r)$ est borné par
\[n^{O(k)} + \left(n^{2 k +1} F(n/k,r) + F(n_1,r-1) + \sum_{i\geq 2} F(n_i,r)\right)\cdot O(k \log n),\]
où $\sum n_i = n$ et $n_i\leq n/2$ pour $i\geq 2$, ou
\[n^{O(k)} + \left(n^{2 k +1} F(n/k,r) + \sum_i F(n_i,r)\right) \cdot O(k \log n),
\]
où $\sum n_i = n$ et $n_i\leq n/2$ pour $i\geq 1$. 
Puisque $k\ll \log n$, ceci est
consistant avec
$F(n,r) = \exp\left(O\left(r+\log n\right)^{c}\right)$ pour $c\geq 3$, ou
même avec
$F(n,r) = \exp\left(O\left((\log n)^{c_1} + (\log r)^{c_2}\right)\right)$
pour $c_1\geq 3$ et $c_2\geq 1$, par exemple. 

Le cas 2a a un coût très similaire, à un facteur constant près.
Le cas 2b et 3 sont différents. Dans les deux cas, nous arrivons à
construire une relation $d$-aire, avec
$d\leq 5$, dans le cas 2b, et $d=k\ll \log n$ dans le cas 3.
Puis, nous appelons Weisfeiler-Leman,
suivi du Lemme des Designs pour des configurations $d$-aires, et,
finalement, Coupe-ou-Johnson.

Weisfeiler-Leman prend un temps $|\Gamma|^{O(d)} = m^{O(d)}$.
Le Lemme des designs
garantit l'existence d'un tuple $(x_1,\dotsc,x_\ell)\in \Gamma$,
$\ell\leq d-1$, avec certaines propriétés. Nous cherchons
un tel tuple
par force brute, ce qui prend un temps $O(m^d)$. Ce qui est plus important est
que ce choix n'est pas canonique. Donc, le temps d'exécution de
tout ce qui reste est multiplié par $m^d = m^{O(\log n)}$.

Coupe-ou-Johnson prend un temps $O(m^d)$. Ici, à nouveau, nous faisons
des choix qui ne sont pas
complètement canoniques; ils imposent un facteur de $m^{O(\log m)}$ sur tout
ce qui suit.
Le résultat de Coupe-ou-Johnson est soit un $\beta$-découpage, ce qui
implique une réduction à un produit du type
$\Alt_{s_1}\wr \Alt_{s_2}$ et/ou à des chaînes plus courtes,
soit un schéma de Johnson, ce qui implique une réduction à $\Alt_{m'}$,
$m'\ll \sqrt{m}$. Donc, soit 
\begin{equation}\label{eq:urg1}
  F(n,r)\leq n^{O(k)} + O\left(k n^{2k+2} F(n/k,r)\right) +
  m^{O(\log n)} \left(1+ F(n_1,r-1) + \sum_{i\geq 2} F(n_i,r)\right),
  \end{equation}
où $\sum n_i = n$,
et $n_i\leq n/2$ pour $i\geq 2$, ou
\begin{equation}\label{eq:urg2}
  F(n,r)\leq n^{O(k)} + O\left(k n^{2k+2} F(n/k,r)\right) +
m^{O(\log n)} \left(1 + \sum_i F(n_i,r)\right),\end{equation}
où $\sum n_i = n$ et $n_i\leq n/2$ pour $i\geq 1$.

Ici $m\leq n$. (Nous pourrions travailler avec une borne moins grossière,
mais cela nous servirait peu.) Donc, les inégalités
(\ref{eq:urg1}) et (\ref{eq:urg2}) sont consistantes avec
$F(n,r) = \exp\left(O\left(r+\log n\right)^{c}\right)$ pour $c\geq 3$.

Comme $r\leq 2 \log_2 n$, nous concluons que le temps total d'exécution de
la procédure pour déterminer les isomorphismes entre deux chaînes de longueur
$n$ est
\[F(n,r) \leq e^{O\left(\log n\right)^3}.\]


\end{document}